\documentclass[10pt]{amsart}
\usepackage{stmaryrd}
\usepackage{amscd,amsxtra,amssymb,mathrsfs}
\usepackage[dvips]{graphics}
\usepackage{epsfig,graphicx}
\usepackage{xypic}
\usepackage[all]{xy}
\xyoption{arc}

\newtheorem{theorem}{Theorem}[section]
\newtheorem{corollary}[theorem]{Corollary}
\newtheorem{lemma}[theorem]{Lemma}
\newtheorem{proposition}[theorem]{Proposition}
\newtheorem{conjecture}[theorem]{Conjecture}
\newtheorem{question}[theorem]{Question}

\theoremstyle{definition}
\newtheorem{definition}[theorem]{Definition}
\newtheorem{remark}[theorem]{Remark}
\newtheorem{example}[theorem]{Example}

\DeclareMathAlphabet{\mathpzc}{OT1}{pzc}{m}{it}

\DeclareMathOperator{\Perf}{\mathsf{Perf}}

\DeclareMathOperator{\SL}{\mathsf{SL}}

\DeclareMathOperator{\rad}{\mathsf{rad}}
\DeclareMathOperator{\coker}{\mathsf{coker}}
\DeclareMathOperator{\chr}{\mathsf{char}}
\renewcommand{\ker}{\mathsf{ker}}
\newcommand{\im}{\mathsf{im}}
\renewcommand{\dim}{\mathsf{dim}}
\DeclareMathOperator{\tor}{\mathsf{tor}}
\DeclareMathOperator{\Coh}{\mathsf{Coh}}

\DeclareMathOperator{\VB}{\mathsf{VB}}

\DeclareMathOperator{\Tri}{\mathsf{Tri}}
\DeclareMathOperator{\Pic}{Pic}

\DeclareMathOperator{\art}{\mathsf{art}}
\DeclareMathOperator{\Art}{\mathsf{Art}}
\DeclareMathOperator{\CM}{\mathsf{CM}}
\DeclareMathOperator{\Hot}{\mathsf{Hot}}
\DeclareMathOperator{\Pro}{\mathsf{pro}}
\DeclareMathOperator{\Rep}{\mathsf{Rep}}

\DeclareMathOperator{\Ob}{\mathsf{Ob}}
\DeclareMathOperator{\Tr}{\mathsf{Tr}}

\DeclareMathOperator{\krdim}{\mathsf{kr.dim}}
\DeclareMathOperator{\prdim}{\mathsf{pr.dim}}
\DeclareMathOperator{\injdim}{\mathsf{inj.dim}}
\DeclareMathOperator{\gldim}{\mathsf{gl.dim}}
\DeclareMathOperator{\syz}{\mathsf{syz}}
\DeclareMathOperator{\add}{\mathsf{add}}
\DeclareMathOperator{\rank}{\mathsf{rank}}

\DeclareMathOperator{\depth}{\mathsf{depth}}
\DeclareMathOperator{\Mod}{\mathsf{Mod}}

\DeclareMathOperator{\edim}{\mathsf{edim}}

\DeclareMathOperator{\Hom}{\mathsf{Hom}}
\DeclareMathOperator{\RHom}{\mathsf{RHom}}

\DeclareMathOperator{\Ext}{\mathsf{Ext}}

\DeclareMathOperator{\GL}{\mathsf{GL}}
\DeclareMathOperator{\Aut}{\mathsf{Aut}}
\DeclareMathOperator{\Ann}{\mathsf{ann}}
\DeclareMathOperator{\End}{\mathsf{End}}
\DeclareMathOperator{\Mat}{\mathsf{Mat}}
\DeclareMathOperator{\Spec}{\mathsf{Spec}}
\DeclareMathOperator{\Specan}{\mathsf{Specan}}
\DeclareMathOperator{\Det}{\mathsf{det}}

\input xy
\xyoption{all}

\newcommand{\kk}{k}

\renewcommand{\mod}{\mathsf{mod}}

\newcommand{\kA}{\mathcal{A}}

\newcommand{\kE}{\mathcal{E}}
\newcommand{\kF}{\mathcal{F}}
\newcommand{\kG}{\mathcal{G}}

\newcommand{\kO}{\mathcal{O}}
\newcommand{\kL}{\mathcal{L}}
\newcommand{\kP}{\mathcal{P}}
\newcommand{\kK}{\mathcal{K}}
\newcommand{\kM}{\mathcal{M}}

\newcommand{\kV}{\mathcal{V}}

\newcommand{\sB}{B}

\newcommand{\sP}{\mathcal{P}}

\newcommand{\lar}{\longrightarrow}

\newcommand{\rA}{A}
\newcommand{\rB}{B}
\newcommand{\rR}{R}
\newcommand{\rS}{S}
\newcommand{\rQ}{Q}
\newcommand{\rK}{K}
\newcommand{\rL}{L}

\newcommand{\idm}{\mathfrak{m}}
\newcommand{\idn}{\mathfrak{n}}
\newcommand{\idp}{\mathfrak{p}}
\newcommand{\idq}{\mathfrak{q}}

\newcommand{\mM}{M}
\newcommand{\mN}{N}
\newcommand{\mK}{K}
\newcommand{\mT}{T}
\newcommand{\mE}{E}
\newcommand{\mF}{F}
\newcommand{\mG}{G}
\newcommand{\mX}{X}
\newcommand{\mD}{D}
\newcommand{\mV}{V}
\newcommand{\mP}{P}

\newcommand{\sX}{X}
\newcommand{\sY}{Y}
\newcommand{\sZ}{Z}
\newcommand{\sE}{E}
\newcommand{\sC}{C}
\newcommand{\sD}{D}
\newcommand{\sT}{T}
\newcommand{\sU}{U}

\newcommand{\DD}{\mathbb{D}}
\renewcommand{\SS}{\mathbb{S}}
\newcommand{\TT}{\mathbb{T}}
\newcommand{\II}{\mathbb{I}}
\newcommand{\OO}{\mathbb{O}}

\newcommand{\bX}{X}

\newcommand{\bU}{U}

\newcommand{\cA}{\mathcal{A}}

\title[Cohen-Macaulay modules over surface  singularities]{Maximal Cohen-Macaulay
modules over \linebreak surface singularities}

\author{Igor Burban}
\address{
Mathematisches Institut
Universit\"at Bonn,
Beringstr. 1,
D-53115 Bonn,
Germany
}
\email{burban@math.uni-bonn.de}

\author{Yuriy Drozd}
\address{%
 Institute of Mathematics
National Academy of Sciences,
Tereschenkivska str. 3,
01004 Kyiv, Ukraine
}
\email{drozd@imath.kiev.ua}

\begin{document}


%


\begin{abstract}
This is a survey article about
properties of Cohen-Macaulay modules over surface singularities.
We discuss properties of the Macaulayfication functor, reflexive
modules over simple, quotient and minimally elliptic singularities,
geometric and algebraic McKay Correspondence. Finally, we
describe matrix factorizations corresponding to indecomposable
Cohen-Macaulay modules over the non-isolated singularities
 $A_\infty$ and  $D_\infty$.
\end{abstract}



\maketitle

\section{Introduction and historical remarks}
The study  of Cohen-Macaulay modules over Noetherian local rings
originates from the theory of integral representations of finite groups
(which, on its side, grew up from a very classical problem of
classification of crystallographic  groups, considered  by Bieberbach, Fedorov,
Schoenflies and others at the end of the 19th century). Namely, let
$G$ be a finite group, $p$ a prime number  and $\rA = \mathbb{Z}_{(p)}$ be the ring of p-adic
integers. Then the category of representations of the group  $G$ over
the ring $\rA$ is equivalent to the category of finitely generated
$\rA[G]$--modules, which are free as $\rA$--modules. This category is not abelian.
However,
it  is extension-closed  in $\rA[G]-\mod$ and has Krull-Remak-Schmidt property. Although
the ring $\rA[G]$ is in general not commutative, it is finite as a module
over its center, i.e.~it is a so-called \emph{order}. A
finitely-generated $\rA[G]$--module which is  free as $\rA$--module,
is  called a \emph{lattice}.

If $G$ is a finite commutative group, then the ring
$\rA[G]$ is commutative of Krull dimension one. In this case, the
category of lattices coincides with the category of
\emph{maximal Cohen-Macaulay modules} this article deals with.

In algebraic geometry, the property of an algebraic variety to be
Cohen-Macaulay (respectively, Gorenstein) describes a restricted class of singular varieties, for
which one has a duality theory in the  way it exists for smooth varieties.
For certain classes of local Noetherian rings of small  dimensions the
category of Cohen-Macaulay modules can be defined in a very simple way.

For example, if $(\rA, \idm)$ is a reduced Noetherian local ring of Krull
dimension one then the category of Cohen-Macaulay modules coincides
with the category of torsion free modules. Note that
 the study of
torsion free modules over curve singularities is important from
the point of view of moduli problems (for example, to investigate  the singularities
of the compactified Jacobian of an irreducible projective curve,
see \cite{AltKle} and \cite{DucoandK}).
If $(\rA, \idm)$ is a \emph{normal surface singularity}, then
the Cohen-Macaulay modules over $\rA$ are precisely the Noetherian
reflexive modules, i.e. finitely generated $\rA$--modules such that
$\mM \cong \mM^{**}$.

First results about the representation type of the category of
Cohen-Macaulay modules were obtained by Drozd-Roiter \cite{DroRoi}
and Jacobinski \cite{Jacob} in the 60-ies
in the framework of integral representations
of finite groups. A suggestion to study homological  properties
of the category of Cohen-Macaulay modules over a Gorenstein local
ring of any Krull dimension as well as the first fundamental results
going in this direction (existence of almost split sequences,
Serre duality in
the stable category) are due to Auslander \cite{PhilNotes}, see
also \cite{AuslICM86}. These ideas have found its further development in
a work of Buchweitz \cite{Buchweitz}, who suggested to study
the stable category of Cohen-Macaulay modules over a Gorenstein local ring as
a  triangulated category, see also \cite{Orlov1}.

In 1978 Herzog has shown that the two-dimensional quotient
singularities have  finitely many indecomposable Cohen-Macaulay
modules \cite{Herzog1}. This results was a starting point and a motivation
for a whole bunch of interesting results  dedicated to the so-called
McKay correspondence. This theory can be divided (although
a little bit artificially)
into two parts: geometric
(due to Artin, Esnault, Gonzalez-Springer, Kn\"orrer,  Verdier and others)
and algebraic (Auslander, Reiten and others).
A combination of both approaches yields a construction of
explicit bijections between the set of
isomorphy classes of irreducible representations of a finite
subgroup $G \subset \SL(2, \mathbb{C})$, the 
irreducible components 
of the exceptional divisor of a minimal resolution
of $\Spec(\mathbb{C}\llbracket x,y\rrbracket ^G)$ and the set of isomorphy classes of 
indecomposable Cohen-Macaulay modules
over $\mathbb{C}\llbracket x,y\rrbracket ^G$. A discussion of these and related results is
the main goal  of our article.

The geometric approach of the study of Cohen-Macaulay modules over the quotient surface
singularities  was generalized by
Kahn \cite{Kahn} to the case of \emph{minimally
elliptic} surface singularities. In the case of \emph{simply elliptic}
singularities he used a classification of Atiyah \cite{Atiyah}
of indecomposable vector bundles on elliptic curves to prove they have
tame Cohen-Macaulay representation type. His programme was completed
by Drozd, Greuel and Kashuba \cite{DGK} who elaborated and extended Kahn's results
on the case of \emph{cusp singularities}.

Recently, Kapustin and Li \cite{KapLi}  revived   the
interest of mathematicians to a study of Cohen-Macaulay modules
over isolated hypersurface singularities because of their
relations with the so-called \emph{Landau-Ginzburg models} arising in  string theory.
In \cite{Orlov} Orlov established an exact connection between the derived category of coherent
sheaves on a projective variety and the stable category of graded Cohen-Macaulay modules
over its homogeneous coordinate algebra. This gave a rigorous meaning to the physical equivalence
between Landau-Ginzburg models and non-linear  sigma models, predicted
by string theorists a long time ago. Another interesting application of 
Cohen-Macaulay modules
to link invariants was discovered  by Khovanov and Rozansky \cite{Khovanov}.

The plan of this article is the following. First, we recall
some basic results from commutative algebra related
to our study of Cohen-Macaulay modules. Then
we discuss some ``folklore results''  on the Macaulayfication functor
in the case of two-dimensional Cohen-Macaulay surface singularities.
Next,  we describe  both the algebraic and the
geometric approaches
to McKay correspondence for quotient surface singularities as well
as its generalization for simply elliptic and cusp singularities.
Finally, we give a new proof of a result
of Buchweitz, Greuel and Schreyer stating  that the
 surface singularities
$A_\infty$ and $D_\infty$ have countable (also called discrete) Cohen-Macaulay
representation type.

\medskip
\noindent
\emph{Acknowledgement}. The first-named author would like to thank Gert-Martin Greuel and
Duco van Straten  for stimulating discussions of the material of 
this article. The research of the first-named author was supported 
by the DFG project  Bu-1866/1-2 and of the second-named author by
the INTAS  grant no. 06-1000017-9093.

\section{Generalities about Cohen-Macaulay modules}

Let $(\rA,\idm)$ be a Noetherian local ring, $\kk = \rA/\idm$ its residue field and $d = \krdim(\rA)$ its Krull dimension.
Throughout the paper
$\rA-\mod$ denotes the category of Noetherian (i.e. finitely generated) $\rA$-modules,
whereas
$\rA-\Mod$ stands for the category of all $\rA$-modules.
In this section we collect some basic facts about  Cohen-Macaulay modules.

\begin{definition}\label{D:Depth}
For a Noetherian $\rA$-module  $\mM$ its \emph{depth} is defined as
$$
\depth_\rA(\mM) = \inf_{i \ge 0}\bigl\{i \, \big| \, \Ext^{i}_\rA(\kk, \mM) \ne 0
\bigr\}.
$$
\end{definition}

\begin{lemma}\label{L:depdim}
Let $\mM$ be a Noetherian $\rA$-module. Then we have:
$$
\depth_\rA(\mM) \le \dim(\mM) := \krdim\bigl(\rA/\Ann(\mM)\bigr),
$$
in particular, the depth of a Noetherian module is always finite. Moreover,
if $I = \Ann_\rA(M)$ then $\depth_\rA(\mM) = \depth_{\rA/I}(\mM)$. Hence, in
what follows, we shall write $\depth(\mM)$ for the depth of $\mM$, omitting  the
subscript.
\end{lemma}

\noindent For a proof of this lemma we refer to
\cite[Proposition 1.2.12]{BrunsHerzog}.

\medskip

\begin{definition}\label{D:CMM}
A  Noetherian  $\rA$--module $\mM$ is  called \emph{maximal Cohen-Macaulay} if
$\depth(\mM) = d$. In what follows, we simply  call such modules Cohen-Macaulay.
\end{definition}

\medskip
\noindent
The following properties follow immediately  from the definition.

\begin{proposition}\label{P:propCM1}
Let $\CM(\rA)$ denote  the full subcategory of Cohen-Macaulay modules. Then
\begin{itemize}
\item The category $\CM(\rA)$ is closed under extensions and direct summands.
\item $\CM(\rA) = \bigl\{\mM \in \rA-\mod\, | \, 
\Ext^i_\rA(\mT, \mM) = 0 \mbox{\,\,\textrm{for all}\,\,\,}
\mT \in \art(\rA), 0 \le i < d \bigr\}$,
where $\art(\rA)$ is the category of Noetherian modules of finite length.
\item We have an exact functor $\CM(\rA) \lar \CM(\widehat{\rA})$.
\item If $\rA$ is Artinian then $\CM(\rA) = \rA-\mod$.
\end{itemize}
\end{proposition}

The given definition of Cohen-Macaulay modules might look like
a little bit artificial.
The following  lemma explains why this
 notion  is natural in the context
of the commutative algebra.

\begin{lemma}\label{L:CMdim1}
Let $(\rA,\idm)$ be a reduced Noetherian ring of Krull dimension \emph{one}. Then an $\rA$--module $\mM$ is
Cohen-Macaulay if and only if it is torsion free.
\end{lemma}

\noindent
\emph{Proof}. Indeed, since the ring $A$ is one-dimensional, the torsion part of a finitely generated
 $\rA$-module is zero-dimensional, hence of  finite length.
Moreover, a Noetherian $\rA$--module $\mM$ is torsion free if
and only if $\Hom_\rA(\kk, \mM) = 0$, what is equivalent to the condition $\depth(\mM) = 1$.
\qed

\begin{definition}\label{D:CMring}
Let $(\rA, \idm)$ be a Noetherian local ring of Krull dimension $d$.
\begin{itemize}
\item The ring $\rA$ is  Cohen-Macaulay if it is Cohen-Macaulay as an $\rA$-module.
\item The ring $\rA$ is  Gorenstein if it is Cohen-Macaulay and
$\Ext^d_\rA(\kk, \rA) \cong  \kk$.
\end{itemize}
\end{definition}

\medskip

\begin{remark}
A Cohen-Macaulay ring need not be reduced. For example, the local ring
$\kk\llbracket x,y\rrbracket /(y^2)$ is Cohen-Macaulay.
\end{remark}

\medskip

\begin{example}\label{E:CMrings} Consider the following examples of Cohen-Macaulay and Gorenstein rings.
\begin{itemize}
\item A regular local ring is always Gorenstein.
\item Let $(\rA, \idm)$ be regular and  $f \in \idm$,
then $\rA/f$ is Gorenstein, see \cite[Proposition 3.1.19]{BrunsHerzog}.
\item Any reduced curve singularity is Cohen-Macaulay.
\item The curve singularity $\kk\llbracket x,y,z\rrbracket /(xy, xz, yz)$ is Cohen-Macaulay but not Gorenstein.
\item The ring $\kk\llbracket x,y,z\rrbracket /(xy, xz)$ is not Cohen-Macaulay.
\item Let $\kk$ be an algebraically closed field of characteristic zero,
$G \subset \GL_n(\kk)$  a finite
subgroup without pseudo-reflections, then the invariant ring
$\kk\llbracket x_1, \dots, x_n\rrbracket ^G$ is always Cohen-Macaulay, see
Theorem \ref{T:quotbasicprop} below. It is Gorenstein
if any only if $G \subseteq
\SL_n(\kk)$, see \cite{WatanabeI} and \cite{WatanabeII}.
\end{itemize}
\end{example}

\begin{example}
Let $\kk$ be an algebraically closed field of characteristic zero,
 $\rA =  \linebreak \kk\llbracket x,y\rrbracket /(y^2 - x^{n+1})$ a simple curve singularity of type
$A_n$. Then the ideals $I_1 = \langle x, y\rangle, I_2 = \langle x^2, y\rangle, \dots,
I_n = \langle x^n, y\rangle$ are Cohen-Macaulay as $\rA$--modules.
The modules $I_2, \dots, I_n$ are always indecomposable and
$I_1 \cong \rR$, where $\rR$ is the normalization of $\rA$. The module $I_1$ is indecomposable
for even $n$  and   $I_1 \cong \rA/(y - x^{\frac{n+1}{2}}) \oplus
\rA/(y + x^{\frac{n+1}{2}})$ for odd $n$.  Moreover, this
is a complete list of indecomposable $\rA$--modules, see for example
\cite{Yoshino}.
\end{example}

\medskip

\begin{theorem}[Auslander-Buchsbaum formula]\label{T:AuslBuchsb}
Let $(\rA,\idm)$ be a Noetherian ring, $\mM$ a Noetherian  $\rA$--module such that  $\prdim(\mM) < \infty$. Then
$$
\prdim(\mM) + \depth(\mM) = \depth(\rA).
$$
\end{theorem}

\medskip
This theorem implies that the structure of Cohen-Macaulay modules over a regular local ring is particularly easy.

\begin{corollary}\label{C:CMreg}
Let $(\rA,\idm)$ be a regular local ring, then a Cohen-Macaulay module $\mM$ is free.
\end{corollary}

\noindent
\emph{Proof}. Since $\rA$ is regular, 
$\krdim(\rA) = \gldim(\rA) = \depth(\rA)$, see for example 
\cite[Section IV.D]{Serre}. 
In particular, the projective dimension  of a finitely generated module $\mM$ 
is always finite.
From the condition $\depth(\mM) = \krdim(\rA)$ it follows $\prdim(\mM) = 0$, i.e. $\mM$ is projective, hence free.
\qed

\medskip

\begin{theorem}[Depth Lemma]\label{T:depthlemma}
Let $0 \to \mM \to \mN \to \mK \to 0$ be a short exact sequence of $A$-modules. Then we have:
\begin{itemize}
\item $\depth(\mN) \ge \min\bigl(\depth(\mM), \depth(\mK)\bigr)$
\item If $\depth(\mN) > \min\bigl(\depth(\mM), \depth(\mK)\bigr)$ then $\depth(\mM) = \depth(\mK) +1$.
\end{itemize}
\end{theorem}

\noindent A proof of this  statement can be found for example in
\cite[Proposition 1.2.9]{BrunsHerzog}.
As its   immediate corollary we obtain.

\begin{proposition}\label{P:syzisCM}
Let $(\rA, \idm)$ be a
Cohen-Macaulay  local ring of Krull dimension $d$ and
 $\mT$ any $\rA$-module.
Then the $d$-th syzygy module $\syz^d(\mT)$ is always Cohen-Macaulay.
\end{proposition}

\medskip
\noindent
In the literature one can find other definitions of
Cohen-Macaulay modules. One possibility to
introduce them uses  regular sequences.

\begin{proposition}\label{P:regseq}
Let $(\rA, \idm)$ be a Noetherian local ring of Krull dimension $d$.
\begin{itemize}
\item A module $\mM$ is Cohen-Macaulay if and only if there exists an $\mM$-regular sequence, i.e.
a sequence
$f_1, \dots, f_d \in \idm$ such that the homomorphism
$$\mM/(f_1, \dots, f_i)\mM \stackrel{f_{i+1}}\lar  \mM/(f_1, \dots, f_i)\mM$$
is injective
for all $0 \le i < d$.
\item If $\rA$ is Cohen-Macaulay, then any $A$-regular sequence is also $\mM$-regular.
\end{itemize}
\end{proposition}

\noindent For a proof of this proposition see
\cite[Chapter 2]{BrunsHerzog}.

\medskip

The following result  says that the category of Cohen-Macaulay modules behaves well under basic functorial
operations of the commutative algebra.

\begin{proposition}\label{P:CMfunct}
Let $(\rA, \idm)$ be a Noetherian local ring.
\begin{itemize}
\item For any non-zero divisor $f \in \idm$ we have a functor $\CM(\rA) \to \CM(\rA/f)$, $\mM \mapsto \mM/f\mM$.
\item For any prime ideal $\mathfrak{p} \in \Spec(\rA)$ we have a functor
$\CM(\rA) \to  \CM(\rA_\mathfrak{p})$, $\mM \mapsto \mM_{\mathfrak{p}}$.
\end{itemize}
\end{proposition}

\noindent For a proof of this proposition,  see
\cite[Theorem 2.1.3]{BrunsHerzog}.

\medskip

\begin{definition}\label{D:canmod}
Let $(\rA,\idm)$ be a Cohen-Macaulay ring of Krull dimension $d$.
A Cohen-Macaulay module $\mK$ is called canonical if
$\Ext^d_\rA(\kk, \mK) \cong \kk$ and
$\injdim(\mK) < \infty$.
\end{definition}

\medskip

\begin{theorem}\label{T:canmod}
Let $(\rA,\idm)$ be a Cohen-Macaulay ring of Krull dimension $d$.
\begin{itemize}
\item If Cohen-Macaulay modules $\mK$ and $\mK'$ are canonical then $\mK \cong \mK'$. Hence, we may use the notation
$\mK = \mK_\rA$.
\item A Noetherian  $\rA$-module $\mK$ is canonical if and only if $\dim_{\kk}\Ext^i_\rA(\kk, \mK) = \delta_{i,d}$.
\item If $\rA$ is Gorenstein then the regular module $\rA$ is canonical.
\item A canonical module exists in $\CM(\rA)$ if and
only if there exists a local Gorenstein
ring $(\rB, \idn)$ and a surjective ring homomorphism  $\rB \to \rA$. In this case
$$\mK_\rA  \cong  \Ext^t_\rB(\rA, \rB),$$
 where $t = \krdim(\rB) - \krdim(\rA)$.
\item If $\rB \to \rA$ is a homomorphism of local Cohen-Macaulay rings such that $\rA$ is
finite as $\rB$-module and $\rB$ has a
canonical module,  then
$\mK_\rA \cong  \Ext^t_\rB(\rA,\mK_\rB)$, where $t = \krdim(\rB) - \krdim(\rA)$.
In particular, if $\kk$ is a field, $\rA =  \kk\llbracket x_1, x_2, \dots, x_n\rrbracket /I$   
a complete $\kk$-algebra
and $\rB \to \rA$ its Noether normalization, then $\mK_\rA = \Hom_\rB(\rA, \rB)$.
\end{itemize}
\end{theorem}

\noindent For a proof of these results we refer to
\cite[Section 3.3]{BrunsHerzog}.

\medskip
\noindent
The importance of the canonical module  becomes clear after   the following theorem.

\begin{theorem}\label{T:CMdual}
Let $(\rA, \idm)$ be a Cohen-Macaulay ring having
a canonical module  $\mK$.
\begin{itemize}
\item For any Cohen-Macaulay module $\mM$ and any integer $t>0$ we have:
$
\Ext^t_\rA(\mM, \mK) = 0.
$ In particular, $\mK$ is an injective object in the exact category of Cohen-Macaulay modules.
\item For any Cohen-Macaulay module $\mM$ the dual module  $\mM^\vee = \Hom_\rA(\mM, \mK)$
is again Cohen-Macaulay.
Moreover, the canonical morphism $\mM \to \mM^{\vee\vee}$ is an isomorphism.
\item The canonical module $\mK$ behaves well under basic functorial operations of the commutative algebra.
In particular,
\begin{itemize}
\item $\widehat{\mK_\rA} \cong  \mK_{\widehat{\rA}}$,
\item for any $\mathfrak{p} \in \Spec(\rA)$ we have: $(\mK_\rA)_{\mathfrak{p}} \cong
 K_{\rA_\mathfrak{p}}$
\item for any non-zero divisor $f \in \idm$ we have: $\mK_\rA/ f \mK_\rA \cong  \mK_{\rA/f}$.
\end{itemize}
\end{itemize}
\end{theorem}

\noindent For a proof of this theorem, see
\cite[Section 3.3]{BrunsHerzog}.

\medskip
By Theorem \ref{T:canmod} the canonical module of a Cohen-Macaulay
ring exists under 
very general assumptions. Nevertheless, there are examples of Noetherian 
local Cohen-Macaulay rings without the canonical module.

\begin{example}\label{E:strangering}
In \cite[Proposition 3.1]{FerrandRaynaud} the authors construct 
an integral local ring $(\rA, \idm)$ of Krull dimension one such that 
its completion $\widehat{\rA}$ contains a minimal prime ideal 
$\idp$ such that $\idp^2 = 0$ and $\idp \cong (\rA/\idp)^n$ for 
$n \ge 2$. Since $\rA$ is a one-dimensional integral domain, it is 
Cohen-Macaulay. Our goal is to show that $\rA$ does not 
have the canonical module.

Indeed, assume $\mK$ is the canonical module of $\rA$. 
Since the quotient ring $\rQ = \rQ(\rA)$ is a field, by Theorem 
\ref{T:CMdual} we have:
$\mK \otimes_\rA \rQ \cong \mK_\rQ = \rQ$. Hence, $\mK$ is isomorphic to 
an ideal in $\rA$. Moreover, the completion $\widehat{\rA}$ is also 
Cohen-Macaulay and its canonical module exists and is an ideal.
By \cite[Korollar 6.7]{HerzogKunz} the total ring of quotients
$\widehat{\rQ} = \rQ(\widehat{\rA})$ is Gorenstein. 

Let $\widetilde\idp = \widehat{\rQ} \otimes_{\widehat\rA} \idp$.
Since localization is an exact functor, we have: 
$\widetilde\idp$ is a minimal prime ideal in $\widehat{\rQ}$, 
$\widetilde\idp^2 = 0$ and $\widetilde\idp \cong 
(\widehat{\rQ}/\widetilde\idp)^n$. By 
\cite[Lemma 1.2.19]{BrunsHerzog} the top of a Gorenstein
Artinian ring  is isomorphic to its socle. However, the top
of $\widehat\rQ$ contains the simple module  
$\widehat\rQ/\widetilde\idp$ with multiplicity one. 
However, the socle of $\widehat\rQ$ contains a semi-simple 
submodule $\widetilde\idp = (\widehat{\rQ}/\widetilde\idp)^n$.
Hence, the ring $\widehat\rQ$ can not be Gorenstein. Contradiction.
\qed
\end{example}

\medskip
\noindent
In what follows, we shall need the notion of the local cohomologies of a Noetherian module $\mM$.

\begin{definition}\label{D:loccohom}
Let $(\rA, \idm)$ be a Noetherian local ring, then the  functor
$\Gamma_{\idm}: \rA-\mod \lar
\rA-\Mod$ is left exact:
$$
\Gamma_{\idm}(\mM) = \varinjlim \Hom_\rA(\rA/\idm^t, \mM) =
\bigl\{x \in \mM\, \big| \, \idm^t x = 0 
\mbox{\textrm{\, \, for some\,\,}} t > 0\bigr\}.
$$
By the definition, $H^i_{\idm}(\mM) := R^i \Gamma_{\idm}(\mM)$.
\end{definition}

\medskip

\begin{remark}
Since taking a direct limit preserves exactness, we have a functorial
isomorphism:
$$
H^i_{\idm}(\mM) \cong  \varinjlim \Ext^i_\rA(\rA/\idm^t, \mM).
$$
\end{remark}

\vspace{0.1cm}

\begin{theorem}\label{T:Groth}
Let $(\rA,\idm)$ be a Noetherian local ring, $\mM$ a Noetherian $\rA$-module, $t = \depth(\mM)$ and
$d = \dim(\mM)$. Then we have the  following properties:
\begin{itemize}
\item All local cohomologies $H^i_{\idm}(\mM)$ are Artinian $\rA$-modules, $i \ge 0$.
\item $H^i_{\idm}(\mM) = 0$ for $ i < t$ and $i > d$. In particular,
$H^i_{\idm}(\mM) = 0$ for $i > \krdim(\rA)$ (Grothendieck's vanishing).
\item $H^t_{\idm}(\mM) \ne 0$ and $H^d_{\idm}(\mM) \ne 0$ (Grothendieck's non-vanishing).
\end{itemize}
\end{theorem}

\noindent For a proof of this theorem  we refer to
\cite[Section 3.5]{BrunsHerzog} and \cite[Chapter 6]{BS}.

\medskip
\noindent
In particular, we obtain an alternative characterization of Cohen-Macaulay modules.

\begin{corollary}\label{C:CMandloccoh}
Let $(\rA,\idm)$ be a Noetherian local ring of Krull dimension $d$, $\mM$ a Noetherian $\rA$-module.
Then $\mM$ is Cohen-Macaulay if and only if $H^i_{\idm}(\mM) = 0$ for $i \ne d$.
\end{corollary}

\noindent
Using this characterization, the following lemma is easy to show.

\begin{lemma}\label{C:CMringext}
Let $(\rA, \idm) \subset (\rB, \idn)$ be a finite extension of local 
Noetherian rings,
$\mM$ a Noetherian $\rB$-module. Then
$\mM$ is Cohen-Macaulay over $\rB$ if and only if it is Cohen-Macaulay over $\rA$.
\end{lemma}

\noindent
\emph{Proof}. Note that the forgetful functor $\rB-\mod \to \rA-\mod$ is exact and
$\Gamma_{\idn}(\mM) \cong  \Gamma_{\idm}(\mM)$ as $\rA$-modules. Hence,
$H^i_{\idn}(\mM) \cong
H^i_{\idm}(\,M)$ as $\rA$-modules, what implies the claim.
\qed

\medskip
\begin{corollary}
Let $\rA = \kk\llbracket x_1,x_2,\dots,x_n\rrbracket /I$ be  Cohen-Macaulay
of Krull dimension $d$ and $\rB = \kk\llbracket y_1, y_2, \dots,y_d\rrbracket 
\to \rA$ be its Noether normalization. Then a Noetherian
$\rA$--module $\mM$ is Cohen-Macaulay if and only if it is free
as $\rB$--module.
\end{corollary}

\noindent
\emph{Proof}. Indeed, by Corollary \ref{C:CMreg} a Cohen-Macaulay
$\rB$--module is free. It remains to apply Lemma \ref{C:CMringext}.
\qed

\medskip
\noindent
Let $\mE = \mE(\kk)$ be the injective envelope of the residue
field $\kk = \rA/\idm$,
$\DD = \Hom_\rA\bigl(-\,\,,\mE(\kk)\bigr): \rA-\Mod \to \rA-\Mod$.

\begin{theorem}[Matlis Duality]\label{T:MatDual}
Let $(\rA,\idm)$ be a local Noetherian ring, $\art(\rA)$ the category of Noetherian  $\rA$-modules
of finite length and $\Art(\rA)$ the category of Artinian $\rA$--modules.
\begin{itemize}
\item The functor
$\DD: \art(\rA) \to \art(\rA)$ is exact and fully faithful, moreover, $\DD^2 \cong Id$.
\item If the ring $\rA$ is complete, then $\DD: \rA-\mod \to \Art(\rA)$ and
$\DD: \Art(\rA) \to \rA-\mod$ are exact fully faithful functors and $\DD^2 \cong Id$.
\end{itemize}
\end{theorem}

\noindent For a proof of this theorem  we refer to
\cite[Proposition 3.2.12 and  Theorem 3.2.13]{BrunsHerzog} and
\cite[Section 10.2]{BS}.

\medskip
The duality functor $\DD$ enters in the  formulation of  the following fundamental result of the
commutative algebra.

\begin{theorem}[Grothendieck's Local Duality]\label{T:locdual}
Let $(\rA,\idm)$ be a Noetherian local ring with a canonical module $\mK_\rA$. Then there exists
an isomorphism of $\delta$-functors
$$
\phi_i: H^i_{\idm} \cong \DD \Ext^{d-i}_\rA(-\,\,, \mK_\rA):
\rA-\mod \to \Art(\rA), \quad i\ge 0.
$$
\end{theorem}

\noindent A proof of this theorem can be found in
\cite[Section 11.2.8]{BS}.

\begin{corollary}
Let $(\rA, \idm)$ be a Cohen-Macaulay  local ring having a canonical module $\mK$.
Then the Cohen-Macaulay $\rA$-modules  are precisely those Noetherian modules
$\mM$ for which
the complex $\RHom_\rA(\mM, \mK) \in \mathsf{Ob}\bigl(D^+(\rA-\mod)\bigr)$
has exactly one non-vanishing cohomology.
\end{corollary}

\medskip
It turns out that in the case of \emph{hypersurface} singularities
the Cohen-Macaulay modules  have the following convenient description.

\begin{proposition}[Eisenbud \cite{Eisenbud}]
Let $(\rS, \idn)$ be a regular local ring of Krull dimension $d \ge 2$, $f \in \idn^2$,  $\rA = \rS/f$ and
$\mM$ a Cohen-Macaulay $\rA$-module without  free direct summands.
Then $\mM$ considered as an $\rA$--module, has a $2$-periodic minimal free resolution.
\end{proposition}

\noindent
\emph{Proof}. Since the ring $\rS$ is regular and $\depth(\mM) = \krdim(\rA) = d -1$, by the Auslander-Buchsbaum formula,
we get:
  $\prdim_\rS(\mM) = 1$. Hence, $\mM$ viewed as an $\rS$-module has a free resolution
$$
0 \to \rS^n \stackrel{\alpha}\lar \rS^n \to \mM \to 0,
$$
where $\alpha \in \Mat_n(\idn)$. Since $\mM$ is annihilated by $f$, we have the following diagram:
$$
\xymatrix
{
0 \ar[r] & \rS^n \ar[r]^\alpha  \ar[d]_{\cdot f} & \rS^n \ar[r]
\ar[ld]^\beta  \ar[d]^{\cdot f} & \mM \ar[r]
\ar[d]^{\cdot f}  & 0 \\
0 \ar[r] & \rS^n \ar[r]_\alpha & \rS^n \ar[r]  & \mM \ar[r]  & 0
}
$$
where $\beta$ is a chain homotopy $(\cdot f, \cdot f) \sim 0$. Hence, we have found a matrix $\beta \in \Mat_n(\rS)$
such that $\alpha \beta = \beta \alpha = f I_n$. In particular, it implies that $\beta \in \Mat_n(\idn)$.

Let $\bar\alpha$ and $\bar\beta$ be the images of $\alpha$ and $\beta$ in $\Mat_n(\idm)$,  where $\idm$ is the maximal ideal
of $\rA$, then  the sequence
$\rA^n \stackrel{\bar\alpha}\lar \rA^n \to M \to 0$ is exact.
Moreover, we claim that
$$
\dots \stackrel{\bar\alpha}\lar \rA^n \stackrel{\bar\beta}\lar \rA^n \stackrel{\bar\alpha}\lar \rA^n
\to \mM \to 0
$$
is a minimal free resolution of $\mM$. Indeed, let $\bar{x} \in \rA^n$ be such that $
\bar\alpha \bar{x} = 0$. This means that there exists an element $y \in \rS^n$ such that
$\alpha x = y$, where $x$ is some preimage of $\bar{x}$ in $\rS^n$. This implies:
$
f x = \beta \alpha x = f \beta y,
$
hence $x = \beta y$ and $\bar{x} = \bar{\beta}\bar{y}$. We have shown: $\ker(\alpha) = \im(\beta)$.
The remaining part is analogous.
\qed

\medskip
\begin{remark}
The pair $(\alpha, \beta)$ is a \emph{matrix factorization}, which corresponds
  to the Cohen-Macaulay module $\mM$.
In these  terms one  can write: $\mM = \mM(\alpha, \beta)$.
\end{remark}

\medskip
\begin{example}
Let $\kk$ be a field of characteristic zero and
$\rA = \kk\llbracket x,y\rrbracket /(y^2 - x^{2n})$  a simple curve singularity of type $A_{2n+1}$.
Then the module $\mM_{\pm} = \rA/(y \pm x^n)$ has the  minimal
free resolution
$$
\dots \lar \rA \stackrel{y \pm x^n}\lar  \rA  \stackrel{y \mp x^n}\lar \rA \lar  \mM_{\pm} \to 0.
$$
The rank one modules $I_i = \langle y, x^{i}\rangle,
\quad 2 \le i \le 2n-1$, have minimal free
resolutions of the following form:
$$
\dots  \lar \rA^2
\xrightarrow{
\left(
\begin{array}{cc}
y & x^{2n -i} \\
x^i & y
\end{array}
\right)}
\rA^2
\xrightarrow{
\left(
\begin{array}{cc}
y & x^{2n -i} \\
x^i & y
\end{array}
\right)}
\rA^2 \lar  I_i \lar  0.
$$
\qed
\end{example}

\medskip
\noindent
Recall the following standard result from the commutative algebra.

\begin{lemma}
Assume a Noetherian  ring $\rA$ is reduced. Then  the total ring of fractions $\rQ(\rA)$ is isomorphic to a direct product of
fields $\rQ_1 \times \rQ_2 \times \dots \times \rQ_t$, where
$\rQ_i = \rQ(\rA/\idp_i)$ is the field of fractions
of $\rA/\idp_i$ and $\idp_1, \idp_2,\dots, \idp_t$
are the minimal prime ideals of $\rA$.
\end{lemma}

\noindent
A proof of this lemma can be found
in \cite[Proposition 10, Section 2.5, Chapter IV]{Bourbaki}
 or \cite[Proposition 1.4.27 and Theorem 1.5.20]{deJongPfister}.

\medskip
\begin{definition}
Let $\rA$ be a reduced Noetherian ring,
$\rQ = \rQ_1 \times \rQ_2 \times \dots \times \rQ_t$ its total ring of fractions
 and $\mM$ a Noetherian $\rA$-module. Then the \emph{multi-rank}
 of $\mM$ is the tuple
of non-negative integers $(r_1, r_2, \dots, r_t)$ such that
$\rQ \otimes_\rA \mM \cong \rQ_1^{r_1} \times \rQ_2^{r_2} \times \dots \times \rQ_t^{r_t}$.
\end{definition}

\medskip
Let $(\rS, \idn)$ be a regular local ring of Krull dimension $d \ge 2$, $\rA = \rS/f$ a reduced hypersurface singularity.
 Since a regular ring is factorial (see \cite[Corollary IV.D.4]{Serre}), we
can write $f = f_1 f_2 \dots f_t$, where all elements $f_i \in \idn$
are irreducible. In these  terms
$$\rQ(\rA) \cong \rQ(\rS/f_1) \times \rQ(\rS/f_2) \times \dots \times \rQ(\rS/f_t).$$

\begin{lemma}\label{L:rankofmf} Let $(\rS, \idn)$ be a regular local ring, $\rA = \rS/f$
a reduced hypersurface singularity,  $\mM$  a Cohen-Macaulay module over $\rA$ and
$(\alpha, \beta) \in \Mat_n(\idn)$   the corresponding matrix factorization. Then
$\Det(\alpha) = u f_1^{r_1} f_2^{r_2} \dots f_t^{r_t}$,  where $u$ is an invertible element in  $\rS$ and
$(r_1, r_2,\dots, r_t)$ is the multi-rank of $\mM$.
\end{lemma}

\noindent
\emph{Proof}. Since $\alpha \beta = \beta \alpha = f^n I_n$, it follows that
$\Det(\alpha) = u f_1^{r_1} f_2^{r_2} \dots f_t^{r_t}$ 
for some unit $u \in \rS$ and 
non-negative integers $r_1, r_2, \dots, r_t$.

By the definition of matrix factorizations, an $\rA$--module
$\mM$ viewed as an $\rS$-module has a free resolution
$
0 \to \rS^n \stackrel{\alpha}\lar \rS^n \to \mM \to 0.
$
For any $i \in \{1, 2,\dots,t\}$ consider the prime ideal 
$\idp_i = (f_i) \subseteq \rS$. Then the localization
$\rS_{\idp_i}$ is a discrete valuation ring and we 
have a free resolution
$$
0 \to \rS_{\idp_i}^n 
\stackrel{\alpha_{\idp_i}}\lar \rS_{\idp_i}^n \to \mM_{\idp_i} \to 0.
$$
Since $\bigl(\Det(\alpha)\bigr)_{\idp_i} = 
\Det(\alpha_{\idp_i})$,
it is easy to see that $\mathsf{length}_{\rA_{\idp_i}}(\mM_{\idp_i})
= r_i$. Since the residue field 
$\kk(\idp_i)$
of the ring $\rS_{\idp_i}$ is isomorphic to $\rQ(\rS/f_i)$, 
it remains to note that
$$
\kk(\idp_i)^{r_i} \cong 
\mM_{\idp_i} \otimes_{\rS_{\idp_i}} \rQ(\rS/f_i) \cong 
\mM \otimes_{\rS} \rQ(\rS/f_i) \cong \mM \otimes_\rA \rQ(\rS/f_i).
$$
\qed

\medskip

\begin{definition}\label{D:isolsing}
Let $(\rA, \idm)$ be a Noetherian local ring. It is called an 
\emph{isolated singularity} if
for any $\mathfrak{p} \in \Spec(\rA)$, $\mathfrak{p} \ne \idm$  
the ring
$\rA_\mathfrak{p}$ is regular. In particular, a regular local ring
is an isolated singularity. 
\end{definition}

\medskip
\begin{lemma}\label{L:isolsing}
Let $\kk$ be a field and 
$\rA = \kk\llbracket x_1, \dots, x_n\rrbracket /f$, where $f \in \idm^2$. Then 
$\rA$ is an isolated singularity if and only if the 
\emph{Tyurina number}
$$
\tau(f) = \dim_{\kk} \bigl(\kk\llbracket x_1, \dots, x_n\rrbracket /(f,j(f))\bigr)
$$
is finite, where $j(f) = \bigl\langle \frac{\partial f}{\partial x_1}, \dots,
 \frac{\partial f}{\partial x_n}\bigr\rangle_\rA$ is the Jacobi ideal of $f$.
\end{lemma}

\noindent
\emph{Proof}. First of all note that $\rA$ is an isolated singularity
if and only if $\rA \otimes_\kk \bar{\kk}$ is, where $\bar\kk$ is the
algebraic closure of $\kk$. Hence, we may without loss of generality 
assume $\kk$ is algebraically closed. 

The singular locus of $\Spec(\rA)$ viewed as a subscheme of 
$\Spec\bigl(\kk\llbracket x_1, x_2, \dots, x_n\rrbracket\bigr)$ is 
$V\bigl(f, j(f)\bigr)$. The singularity $\rA$ is 
isolated if and only if $V\bigl(f, j(f)\bigr) = (0, 0, \dots,0)$.
By Hilbert-R\"uckert's Nullstellensatz this 
 is equivalent for  the ideal $\bigl(f, j(f)\bigr)$ to be 
$(x_1, x_2, \dots,x_n)$-primary in $\kk\llbracket x_1, \dots, x_n\rrbracket $.
The last condition is equivalent to the finiteness of 
$\tau(f)$.
\qed

\medskip

\begin{remark}\label{R:isolsing}
Let $\rA = \kk\llbracket x_1, \dots, x_n\rrbracket /f$, where $f \in \idm^2$ and $\kk$ be  
a  field of \emph{characteristic zero}. 
Then $\rA$ is an isolated singularity
if and only if the \emph{Milnor number}
$$
\mu(f) = \dim_{\kk} \bigl(\kk\llbracket x_1, \dots, x_n\rrbracket /j(f)\bigr)
$$
is finite,  see for example \cite[Lemma 23]{GLS}. 

Note that 
this is no longer true if the characteristic of $\kk$ is positive. 
For example, let $\mathsf{char}(\kk) = 3$ and 
$\rA = \kk\llbracket x,y\rrbracket /(x^2 + y^3)$. Since $\tau(f) = 3$, $\rA$ is 
an isolated singularity. However, $\mu(f) = \infty$.
\end{remark}

\medskip
\begin{example}
Consider the so-called
$T_{2,3,\infty}$ hypersurface singularity
 $\rA = \kk\llbracket x,y,z\rrbracket /(x^2 + y^3 + xyz)$.
Then the singular locus of $\Spec(\rA)$ is given by the ideal $I=(x, y)$.
Since $V(I) = \Spec(\kk\llbracket t\rrbracket )$,
this singularity is not isolated.
\end{example}

\medskip
\begin{definition}\label{D:locfr}
Let $(\rA, \idm)$ be a Cohen-Macaulay ring. A Noetherian module  $\mM$
is called \emph{locally free on the punctured spectrum} if the
$\rA_\idp$-module $\mM_{\mathfrak{p}}$
is free for any prime ideal $\mathfrak{p} \ne \idm$.
\end{definition}

\medskip

\noindent
The following lemma is straightforward.

\begin{lemma}\label{L:locfr}
Let $(\rA, \idm)$ be a Cohen-Macaulay  isolated singularity. Then any Cohen-Macaulay $\rA$-module
is locally free on the punctured spectrum.
\end{lemma}

\noindent
\emph{Proof}. Let $\mM$ be a Cohen-Macaulay $\rA$-module, then
for any  $\mathfrak{p} \in \Spec(\rA)\setminus \{\idm\}$ the ring
$\rA_\mathfrak{p}$ is regular and the module $\mM_\mathfrak{p}$
is Cohen-Macaulay, hence free by Corollary \ref{C:CMreg}.
\qed

\medskip
\begin{definition}
Let $\kk$ be a field. Its \emph{real valuation} is
a function $\| \,\|: \kk \to \mathbb{R}_{\ge 0}$ such that
\begin{itemize}
\item $\| ab \| = \| a\| \cdot \| b\| $
\item $\| a + b \| \le  \| a\| + \| b\|$.
\end{itemize}
For a formal power series
$f = \sum\limits_{\alpha \in \mathbb{N}^n}
c_\alpha x^\alpha \in \kk\llbracket x_1,x_2, \dots,x_n\rrbracket $ and
$\varepsilon \in \mathbb{R}_{> 0}$ we
denote
$$
\| f\|_\varepsilon :=
\sum\limits_{\alpha \in \mathbb{N}^n}
\|c_\alpha \|\varepsilon^{|\alpha|} \in \quad \mathbb{R} \cup
\{\infty\}.
$$
A power series $f$ is called \emph{convergent} with respect to a
valuation $\| \, \|$ if there exists
$\varepsilon \in \mathbb{R}_{> 0}$ such that
$\| f\|_\varepsilon < \infty$. Let $\kk\{x_1,x_2,\dots,x_n\}$
denote the ring of convergent  power series. A $\kk$-algebra 
$\rA$ is called \emph{analytic} if it is
isomorphic to an algebra of the form
$\kk\{x_1,x_2,\dots,x_n\}/I$.
\end{definition}

\medskip
\begin{remark}
If $\| a \| = 0$ for all $a \in \kk$ then the ring of convergent power
series coincides with the ring of formal power series.
\end{remark}

\medskip
\begin{definition}
Let $(\rA, \idm)$ be a Noetherian local ring and
$\kk = \rA/\idm$ be its residue field. The ring
$\rA$ is called \emph{Henselian} if for any polynomial
$p = p(t) \in \rA[t]$ such that $p(t)\equiv \bar p_1(t)\bar p_2(t)\ \mod\ \idm$,
 where $\bar p_1(t),\bar p_2(t)$ are coprime in $\kk[t]$, there exist
polynomials $p_1(t), p_2(t) \in \rA[t]$ such that
$p(t) = p_1(t) p_2(t)$ and $p_i \equiv \bar p_i \, \,  \mod \, \,\idm$
for $i=1,2$.
\end{definition}

\vspace{0.1cm}

\begin{theorem}
Let $\kk$ be a field and $\rA = \kk\{x_1, x_2, \dots, x_n\}/I$ 
be a local  analytic $\kk$-algebra. Then it is Henselian.
\end{theorem}

\noindent
A proof of this result can be found in \cite[Theorem 1.17]{GLS}.

\medskip

Let $(\rA, \idm)$ be a Henselian
Cohen-Macaulay local ring having a canonical module. First of all note that
the category $\rA-\mod$ of all Noetherian $\rA$-modules is a
\emph{Krull-Remak-Schmidt category}, see \cite[Proposition 30.6]{CurtisReiner} and
\cite[Theorem A.3]{Wahl}.
 Hence, the
category $\CM(\rA)$ of Cohen-Macaulay modules has
Krull-Remak-Schmidt property, too.

The main goal  of this
article is to describe the  surface singularities having  finite, discrete and
tame Cohen-Macaulay representation type. In particular,
one can  pose a question about the existence of almost
split sequences.

\begin{theorem}[Auslander]
Let $(\rA, \idm)$ be a Henselian local ring having a canonical module, $\mM$  a non-free
indecomposable  Cohen-Macaulay
$\rA$-module, which is locally free on the punctured spectrum. Then there exists an almost split sequence
$$
0 \to \tau(\mM) \to \mN \to \mM \to 0
$$
ending at $\mM$. Moreover, the category $\CM(\rA)$ admits almost split sequences if and
only if $\rA$  is an isolated singularity.
\end{theorem}

\noindent
For a proof of this Theorem, see 
\cite{PhilNotes} and \cite[Chapter 3]{Yoshino}.

\medskip

This theorem naturally raises the question about an explicit description of the Auslander-Reiten translation
$\tau$.

\medskip

\begin{definition}\label{D:AuslTr} Let $(\rA, \idm)$ be  a local Cohen-Macaulay ring, $\mM$ a Noetherian $\rA$-module
and $\mG \stackrel{\varphi}\lar \mF \lar \mM \to  0$ a free presentation of $M$.
The Auslander-transpose of $\mM$ is defined via the following exact sequence:
$$
0 \to \mM^* \lar \mF^* \stackrel{\varphi^*}\lar  \mG^* \lar \Tr(\mM) \to 0.
$$
\end{definition}

\noindent
In these  terms,  we have the following theorem.

\medskip
\begin{theorem}[Auslander]
Let $(\rA, \idm)$ be a Henselian local ring of Krull dimension $d$, $\mK$ the canonical
$\rA$-module and $\mM$ a non-free indecomposable
Cohen-Macaulay $\rA$-module, which is locally free on the punctured spectrum. Then
$$
\tau(\mM) \cong  \syz^d\bigl(\Tr(\mM)\bigr)^\vee,
$$
where $\mN^\vee = \Hom_\rA(\mN, \mK)$. If $\rA$ is moreover Gorenstein, then
$\tau(\mM) \cong \syz^{2-d}(\mM)$.
\end{theorem}

\noindent
For a proof we refer to \cite[Proposition 8.7, page 105]{PhilNotes}
and \cite[Proposition 1.3, page 205]{PhilNotes}, see also 
\cite[Proposition 3.11]{Yoshino}.

\section{Cohen-Macaulay modules over surface singularities}

Throughout this section,  let $(\rA,\idm)$ be a \emph{reduced} Cohen-Macaulay singularity
of Krull dimension \emph{two}, having a canonical module $\mK_\rA$. Let $\sP$ denote the
set of prime ideals in $\rA$ of height $1$.

\begin{lemma}\label{L:CMdim2}
Let $\mN$ be a Cohen-Macaulay $\rA$-module and $\mM$ a Noetherian $\rA$-module. Then
the $\rA$-module $\Hom_\rA(\mM, \mN)$ is Cohen-Macaulay.
\end{lemma}

\noindent
\emph{Proof}. From  a free presentation $\rA^n \stackrel{\varphi}\to \rA^m \to M \to 0$
of $M$ we obtain an exact sequence:
$$
0 \to \Hom_\rA(\mM, \mN) \to \mN^m \stackrel{\varphi^*}\lar \mN^n \to \coker(\varphi^*) \to 0.
$$
Since $\depth_\rA(\mN) = 2$, applying the Depth Lemma twice we obtain:
$$\depth_\rA\bigl(\Hom_\rA(\mM, \mN)\bigr) \ge  2.$$ Hence, $\Hom_\rA(\mM, \mN)$ is Cohen-Macaulay.
\qed

\medskip
\begin{proposition}\label{P:Macaulafic}
In the notations of this  section, the canonical embedding functor
$\CM(\rA) \to \rA-\mod$ has a left adjoint $$\mM \mapsto \mM^\dagger := \mM^{\vee\vee} =
\Hom_\rA\bigl(\Hom_\rA(\mM, \mK_\rA), \mK_\rA\bigr).$$
\end{proposition}

\noindent
\emph{Proof}. Note that for any Noetherian module $\mM$
the   $\rA$-module $\mM^\dagger$ is Cohen-Macaulay by Lemma \ref{L:CMdim2}. Next, for any Noetherian $\rA$-module $\mM$
there exists an  exact sequence
$$
0 \to \tor(\mM) \to \mM \stackrel{i_\mM}\lar  \mM^\dagger \to \mT \to 0,
$$
where $\tor(\mM)$ is the torsion part of $\mM$ and $i_\mM$ is the canonical morphism.

\vspace{0.3cm}

\noindent $\bullet$ Let us first assume $\mM$ to be 
torsion free, so we have a short exact
sequence $$0  \to \mM \stackrel{i_\mM}\lar  \mM^\dagger \to T \to 0.$$
 Since  for any $\idp \in \sP$  the ring $\rA_\idp$ is reduced and Cohen-Macaulay of Krull dimension one and
the module
$\mM_\idp$ is torsion free,  it is Cohen-Macaulay.
By Theorem  \ref{T:CMdual}  we have: $(\mK_\rA)_\idp \cong \mK_{\rA_\idp}$, hence
the morphism  $(i_\mM)_\idp: \mM_\idp \to (\mM^\dagger)_\idp$ is an isomorphism.
This means that for any $\idp \in \sP$ we have: $\mT_\idp = 0$, hence $\mT$ is a finite length module.
If $\mN$ is a Cohen-Macaulay module, it follows from the  exact sequence
$$
\Hom_\rA(\mT, \mN) \to \Hom_\rA(\mM^\dagger, \mN) \to \Hom_\rA(\mM, \mN) \to
\Ext^1_\rA(\mT, \mN)$$
and equalities  $\Hom_\rA(\mT, \mN) = 0 = \Ext^1_\rA(\mT, \mN)$, that the  canonical morphism
$(i_\mM)_*: \Hom_\rA(\mM^\dagger, \mN) \to \Hom_\rA(\mM, \mN)$ is an isomorphism.

\vspace{0.3cm}

\noindent $\bullet$  Now let $\mM$ be an arbitrary $\rA$-module. Since a Cohen-Macaulay module
$\mN$ is always torsion-free, we have a canonical
 isomorphism $$\Hom_\rA\bigl(\mM/\tor(\mM), \mN\bigr) \cong \Hom_\rA\bigl(\mM, \mN).$$
In the commutative diagram
$$
\xymatrix
{
\Hom_\rA\bigl((\mM/\tor(\mM))^\dagger, \mN\bigr) \ar[rr] \ar[d] & &
\Hom_\rA\bigl(\mM/\tor(\mM), \mN\bigr) \ar[d] \\
\Hom_\rA(\mM^\dagger, \mN) \ar[rr] & & \Hom_\rA(\mM, \mN)
}
$$
both vertical arrows and the first horizontal arrow are isomorphisms. Hence, the canonical
morphism  $(i_\mM)_*: \Hom_\rA(\mM^\dagger, \mN) \to  \Hom_\rA(\mM, \mN)$ is an isomorphism for any $\rA$-module
$\mM$ and a Cohen-Macaulay module $\mN$.
\qed

\medskip

\begin{definition}\label{D:biratisom} Let $\mM$ and $\mN$ be two Noetherian $\rA$-modules.
A morphism $f: \mM \to \mN$ is called \emph{birational isomorphism} if
the induced map $1 \otimes f: Q(\rA) \otimes \mM \to Q(\rA) \otimes \mN$ is an isomorphism of $Q(\rA)$-modules.
\end{definition}

\medskip
\begin{lemma}\label{L:biratisom}
For a Noetherian module $\mM$ the canonical morphism $\mM \to \mM^\dagger$ is a birational isomorphism.
\end{lemma}

\noindent
\emph{Proof}. Indeed, since in the exact sequence
$$
0 \to \tor(\mM) \to \mM \stackrel{i_\mM}\lar  \mM^\dagger \to \mT \to 0
$$
the modules $\tor(\mM)$ and $\mT$ are torsion Noetherian modules and $Q(\rA)$ is a flat
$\rA$-module, the claim follows.
\qed

\medskip
\begin{lemma}
A Noetherian module $\mM$ is Cohen-Macaulay if and only if the canonical morphism
$\mM \to \mM^\dagger$ is an isomorphism.
\end{lemma}

\noindent
\emph{Proof}. Since $\mM^\dagger$ is always Cohen-Macaulay, one direction is clear.
Assume now  that $\mM$ is Cohen-Macaulay. From the canonical isomorphism
$\End_\rA(\mM) \cong  \Hom_\rA(\mM^\dagger, \mM)$ 
we obtain  a morphism
$\pi_\mM: \mM^\dagger \to \mM$ such that $\pi_\mM \circ i_\mM = 1_\mM$.
It is clear that $\pi_\mM$ is an epimorphism. Since it is a birational isomorphism, the multi-rank of $\ker(\pi_\mM)$ is zero, hence it it a torsion module. But
$\mM$ is torsion free, hence $\pi_\mM$ is a monomorphism.
\qed

\medskip
\begin{lemma}\label{L:isocodone}
Let $\mM$ and $\mN$ be two Cohen-Macaulay $\rA$-modules. Then a morphism
$f: \mM \to \mN$ is an isomorphism if and only if for all 
$\idp \in \sP$ the morphism 
$f_\idp: \mM_\idp \to \mN_\idp$ is an isomorphism.
\end{lemma}

\noindent
\emph{Proof}. If $f_\idp$ is an isomorphism
for all $\idp \in \sP$, then
$f$ is a birational isomorphism, hence $\ker(f)$ has zero multi-rank. Since
$\mM$ is torsion free, $\ker(f) = 0$ and $f$ is injective.

Next,  consider  the short exact sequence
$
0 \to \mM \stackrel{f}\to \mN \to \mT \to 0.
$
Since $\mT_\idp = 0$ for all $\idp \in \kP$, 
the module $\mT$ is of  finite length. Moreover, since
$\Hom_\rA(\mT, \mK_\rA) = 0 = \Ext^1_\rA(\mT, \mK_\rA)$, the induced
morphism  $\Hom_\rA(\mN, \mK_\rA) \to \Hom_\rA(\mM, \mK_\rA)$ is an isomorphism, hence $f^\dagger: \mM^\dagger \to
\mN^\dagger$ is an isomorphism, too.
The claim follows now from the commutative  diagram
$$
\xymatrix
{
\mM \ar[rr]^f \ar[d] & & \mN \ar[d]  \\
\mM^\dagger \ar[rr]^{f^\dagger} & & \mN^\dagger
}
$$
and the fact that the vertical arrows and $f^\dagger$ are isomorphisms.
\qed

\medskip

\begin{proposition}\label{P:CMrefl} Let $(\rA, \idm)$ be a Noetherian local ring
of Krull dimension two, which  is Gorenstein in codimension one. Then for any Noetherian
 $\rA$-module $\mM$ we have a functorial
isomorphism $\mM^* \to \mM^\vee$, where $\mX^* = \Hom_\rA(\mX, \rA)$.
In particular, $\mM^\dag\simeq\mM^{**}$, so the Cohen-Macaulay modules over $\rA$ are precisely the
reflexive modules.
\end{proposition}

\noindent
\emph{Proof}. The ring $\rA$ is Gorenstein in codimension one if any only if
for all $\idp \in \sP$ we have: $\mK_{\rA_\idp} \cong \rA_\idp$. 
Let $\rQ = \rQ(\rA)$ be the total ring of fractions of $\rA$. 
Since $\rQ = \rQ(\rA_\idp)$ for $\idp \in \sP$, by Theorem 
\ref{T:CMdual} we get: $\rQ$ is a  Gorenstein 
ring of Krull dimension zero and 
$$\mK_\rA \otimes \rQ \cong \mK_{\rQ} \cong 
\rQ \cong \rA \otimes_\rA \rQ.$$
In particular,
the modules $\rA$ and $\mK_\rA$ are birationally isomorphic. Moreover, for any
birational isomorphism $j: \rA \to \mK_\rA$ and any Noetherian module $\mM$ the induced morphism
of Cohen-Macaulay modules
 $j_*: \Hom_\rA(\mM, \rA) \to \Hom_\rA(\mM, \mK_\rA)$ is an isomorphism in codimension one.
By Lemma \ref{L:isocodone} it is an isomorphism.
\qed

\medskip
\noindent
Recall the following well-known result of Serre.

\begin{theorem}[Serre]
A two-dimensional singularity $(\rA, \idm)$ is normal if and only if it is
Cohen-Macaulay and isolated.
\end{theorem}

\noindent
For a proof, see  \cite[Theorem IV.D.11]{Serre}. 

\medskip
\noindent
From Proposition \ref{P:CMrefl} and Serre's theorem we immediately obtain the following
corollary.

\medskip
\begin{corollary}\label{C:CMrefl}
Let $(\rA, \idm)$ be a normal two-dimensional singularity. Then a Noetherian
$\rA$-module is Cohen-Macaulay if and only if it is reflexive.
\end{corollary}

\medskip
Let $(\rA, \idm)$ be a reduced Cohen-Macaulay ring of Krull dimension two,
 $\bX = \Spec(\rA)$, $\bU = \bX\setminus\{\idm\}$ and $i: \bU \to \bX$
the inclusion map. Let $\mM \mapsto \widetilde\mM$ and $\kF \mapsto \Gamma(\kF)$
denote  the quasi-inverse equivalences $\rA-\mod \to \Coh(\bX)$ and
$\Coh(\bX) \to \rA-\mod$.

\begin{proposition}\label{P:geominterpr}
Let $\mM$ be a torsion free Noetherian $\rA$-module. Then there is a functorial isomorphism
$\mM^\dagger \to \Gamma\bigl(i_* i^*\widetilde\mM\bigr)$.
\end{proposition}

\noindent
\emph{Proof}. First of all, if $\mM$ is a Cohen-Macaulay $\rA$-module then
the canonical morphism $M \to \Gamma\bigl(i_* i^*\widetilde\mM\bigr)$ is an isomorphism.
Indeed, by \cite[Corollaire 2.9]{Grothendieck} we have an exact sequence
$$
0 \to H^0_{\idm}(\mM) \to \mM \to \Gamma\bigl(i_* i^*\widetilde\mM\bigr)
\to  H^1_{\idm}(\mM) \to 0.
$$
It remains to note that 
$H^0_{\idm}(\mM) = 0 =  H^1_{\idm}(\mM)$ for a Cohen-Macaulay module
$\mM$.

Next, if $\mM$ is torsion free then we have an exact sequence
$0 \to \mM \to \mM^\dagger \to \mT \to 0$ where $\mT$ is a module of finite length.
Since $i$ is an open embedding, the functor $i^*$ is exact and $i^*\widetilde\mT = 0$, so
we obtain an isomorphism
$i^*\widetilde\mM \to i^*\widetilde\mM^\dagger$. The claim follows from the
commutative diagram
$$
\xymatrix
{ \mM \ar[rr] \ar[d] & & \mM^\dagger \ar[d]^\cong \\
\Gamma(i_*i^*\widetilde\mM) \ar[rr]^\cong  & & \Gamma(i_*i^*\widetilde\mM^\dagger).
}
$$
\qed

\medskip
\begin{corollary}
Let $(\rA,\idm)$ be a reduced Cohen-Macaulay ring of Krull dimension two,
$\mM$ and $\mN$ two Cohen-Macaulay modules. Then a morphism $f: \mM \to \mN$ is an epimorphism if and only if
for all $\idp \in \kP$ the morphism $f_\idp: \mM_\idp \to \mN_\idp$ is an epimorphism.
\end{corollary}

\noindent
\emph{Proof}.  The second condition is equivalent to the fact that $\mT := \coker(f)$ is a module of finite length.
Since the functor $i^*$ is exact and $i^*(\tilde\mT) = 0$, the result follows from the functorial isomorphisms
$\mM \to \Gamma(i_* i^*\widetilde\mM)$ and $\mN  \to \Gamma(i_* i^*\widetilde\mN)$.
\qed

\medskip
\noindent
The following result is also a consequence of Proposition \ref{P:geominterpr}.

\medskip
\begin{corollary}
Let $(\rA, \idm)$ be a normal surface singularity, $\VB(U)$ the category of locally free $\kO_U$--modules.
 Then  the functor $i^*: \CM(\rA) \to \VB(\bU)$, mapping a Cohen-Macaulay module
$\mM$ to the locally free sheaf $i^*\widetilde\mM$,  is an equivalence of categories.
\end{corollary}

\noindent
\emph{Proof}. Let $\mM$ be a Cohen-Macaulay $\rA$-module. 
Since for any $\idp \in \sP$ the ring $\rA_\idp$ is regular,
the module $\mM_\idp$ is free. Hence, the coherent sheaf $i^*\widetilde\mM$
is indeed locally free.

Let $\kF$ be a locally free sheaf on $\bU$, then the direct image sheaf
$\kG := i_*\kF$  is quasi-coherent. However, any quasi-coherent sheaf on a Noetherian scheme
can be written  as  the direct limit of an increasing  sequence of coherent subsheaves
$\kG_1 \subseteq \kG_2 \subseteq \dots \subseteq \kG$. Since the functor
$i^*$ is exact, we obtain an increasing filtration
$i^*\kG_1  \subseteq i^*\kG_1  \subseteq \dots \subseteq
i^*\kG$. But $i^*\kG = i^*i_*\kF \cong \kF$. Since the scheme $\bU$ is Noetherian and $\kF$ is coherent,
it implies that $\kF \cong i^*\kG_t$ for some  $t \ge 1$. Moreover, since
the module $\mG_t = \Gamma(\kG_t)$ is torsion-free on the punctured spectrum,
the morphism  $\mG_t \to \mG_t^\dagger =: \mM$ 
induces an isomorphism
$i^*\widetilde\mG_t \cong i^*\widetilde\mG_t^\dagger \cong \kF$. Hence, the functor
$i^*$ is dense. Moreover, by Proposition \ref{P:geominterpr} 
we have: $i_*\kF \cong i_* i^*\widetilde\mM \cong \widetilde\mM$,
hence $i_*\kF$ is always coherent. Since $i^* i_* \kF \cong \kF$ for any
$\kO_\bU$-module  $\kF$, it is easy to see that the functor $\VB(\bU) \to
\CM(\rA)$ given by
$\kF \mapsto \Gamma(i_*\kF)$  is  quasi-inverse to $i^*$.
\qed

\medskip

\begin{remark}
It can be shown that for an isolated surface singularity $(\rA, \idm)$ the abelian category
$\Coh(U)$ is hereditary. Hence, the category of Cohen-Macaulay modules on a normal surface singularity
 can be interpreted as the category of vector bundles on a certain ``non-compact''
 smooth curve.
\end{remark}

\medskip

\medskip
\begin{definition}
Let $(\rA, \idm)$ be a reduced Cohen-Macaulay ring of Krull dimension two
and $\mM$ be a Cohen-Macaulay module over $\rA$. The interior
tensor product functor  $\mM \boxtimes_\rA \,-\,\,: \CM(\rA) \to \CM(\rA)$ is defined as
$\mM \boxtimes \mN := (\mM \otimes_\rA \mN)^\dagger$.
\end{definition}

\medskip
\begin{proposition}\label{P:tensandhom}
Let $\mN$ be a Cohen-Macaulay module. Then the interior tensor product functor
 $\mN \boxtimes_\rA \,-\,\,$ is left adjoint to the interior $\Hom$-functor
$\Hom_\rA(\mN, \,-\,)$.
Moreover, one also has the following canonical isomorphisms in the category $\CM(\rA)$:
$$
\mM_1 \boxtimes \mM_2 \cong \mM_2 \boxtimes \mM_1, \quad
(\mM_1 \boxtimes \mM_2) \boxtimes \mM_3 \cong
\mM_1 \boxtimes (\mM_2 \boxtimes \mM_3).
$$
\end{proposition}

\noindent
\emph{Proof}. For any $\rA$-modules $\mM, \mN$ and $\mK$ there is a canonical isomorphism
$$
\Hom_\rA(\mM \otimes \mN, \mK) \cong \Hom_\rA\bigl(\mM, \Hom_\rA(\mN, \mK)\bigr).
$$
If $\mK$ is Cohen-Macaulay then $\Hom_\rA(\mN, \mK)$ is Cohen-Macaulay, too. Moreover,
we have a natural  isomorphism
$
\Hom_\rA\bigl((\mM \otimes \mN)^\dagger, \mK\bigr)
\cong \Hom_\rA(\mM \otimes \mN, \mK)
$
implying that
$$
\Hom_\rA(\mM \boxtimes_\rA \mN, \mK) \cong \Hom_\rA\bigl(\mM,
\Hom_\rA(\mN, \mK)\bigr)
$$
for any Cohen-Macaulay modules $\mM$, $\mN$ and $\mK$.

\vspace{0.3cm}

In order to prove the second part of the proposition, consider  the composition map
$$
\imath: (\mM_1 \otimes \mM_2) \otimes \mM_3 \to
(\mM_1 \otimes \mM_2)^\dagger  \otimes \mM_3 \to
\bigl((\mM_1 \otimes \mM_2)^\dagger \otimes \mM_3)\bigr)^\dagger
$$
which
is a birational isomorphism. Moreover, the cokernel of $\imath$  is a module of finite length. By the universal property of the
Macaulayfication functor  we obtain a commutative diagram
$$
\xymatrix
{
\mM_1 \otimes \mM_2 \otimes \mM_3 \ar[rr]^\imath \ar[d]_\jmath & &
\bigl((\mM_1 \otimes \mM_2)^\dagger \otimes \mM_3)\bigr)^\dagger \\
(\mM_1 \otimes \mM_2 \otimes \mM_3)^\dagger \ar[urr]_{\varphi} & &
}
$$
where all morphisms $\imath$, $\jmath$ and $\varphi$ are birational isomorphisms.
Moreover, the cokernel of $\varphi$
is a quotient of the cokernel of $\imath$, hence it has finite length.
By Lemma \ref{L:isocodone} this implies that  $\varphi$  is an isomorphism. Hence,
$$
(\mM_1 \boxtimes \mM_2) \boxtimes \mM_3 \cong
(\mM_1 \otimes \mM_2 \otimes \mM_3)^\dagger \cong
\mM_1 \boxtimes (\mM_2 \boxtimes \mM_3),
$$
implying the claim. The proof of the remaining statement is similar.
\qed

\medskip
\begin{remark}
Let $(\rA, \idm)$ be a normal surface singularity. Then the equivalence of categories
$i^*: \CM(\rA) \to \VB(\bU)$ additionally satisfies the following property:
$$
i^*(\mM_1 \boxtimes \mM_2) \cong i^*\mM_1 \otimes i^*\mM_2.
$$
\end{remark}

\medskip
\medskip

Let $(\rA, \idm)$ be a Henselian Cohen-Macaulay local ring and
$\rA \subseteq \rB$ a finite ring extension. Then the ring $\rB$ is semi-local.
Moreover,
$\rB \cong  (\rB_1, \idn_1) \times (\rB_2, \idn_2) \times \dots \times (\rB_t, \idn_t)$, where all
$(\rB_i, \idn_i)$ are local. Assume all the rings  $\rB_i$ are Cohen-Macaulay.

\begin{proposition}
The functor $\rB \boxtimes_\rA - : \CM(\rA) \to \CM(\rB)$ mapping
 a Cohen-Macaulay module
$\mM$ to $(\rB \otimes_\rA \mM)^\dagger$ is left adjoint to the forgetful functor
$\CM(\rB) \to \CM(\rA)$.
\end{proposition}

\noindent
\emph{Proof}. Let $\mN$ be a Cohen-Macaulay $\rB$--module and $\mM$ a Cohen-Macaulay $\rA$--module.
Then we have functorial isomorphisms:
$$
\Hom_\rA(\mM, \mN) \cong \Hom_\rB(\rB \otimes_\rA \mM, \mN) \cong \Hom_\rB\bigl((\rB \otimes_\rA \mM)^\dagger, \mN\bigr),
$$
implying  the claim.
\qed

\begin{proposition}\label{P:Macdiffbase}
Let $(\rA, \idm)$ be a reduced  Henselian Cohen-Macaulay  local ring of Krull dimension two,
$\rA \subseteq \rB$ a finite ring extension such that $\rB$ is reduced 
and Cohen-Macaulay. Let $\mM$ be a Noetherian
$\rB$--module, then $\mM^{\dagger_\rA} \cong \mM^{\dagger_\rB}$ as $\rA$--modules.
\end{proposition}

\noindent
\emph{Proof}.  By Corollary \ref{C:CMringext}, the $\rB$-module
$\mM^{\dagger_\rB}$ is Cohen-Macaulay as an $\rA$--module.  
Let  \linebreak  $\bigl\{\idp_1, \idp_2, \dots,\idp_n\bigr\}$  be the set 
of  minimal prime ideals of $\rA$, then 
$\rQ(\rA) \cong \rA_{\idp_1} \times \rA_{\idp_2} 
\times \dots  \times \rA_{\idp_n}$. Since the ring extension
$\rA \subseteq \rB$ is finite, for a given minimal prime ideal 
$\idp \in \Spec(\rA)$ 
the set $I(\idp) := \bigl\{\idq \in \Spec(\rB) \, \big| \, 
\idp \subseteq \idq, \mathsf{ht}(\idq) = 0\bigr\}$ 
is non-empty and finite.  Moreover, we have:
$\rB_\idp \cong  \prod_{\idq \in I(\idp)} \rB_\idq$. This implies that
 $\rQ(\rB)\cong \rQ(\rA)\otimes_{\rA}\rB$, hence
 $\rQ(\rA) \otimes_\rA \mM \cong \rQ(\rB) \otimes_\rB \mM$
and 
the torsion part $\tor_\rB(\mM)$ of the $\rB$-module $\mM$ coincide with
 the torsion part $\tor_\rA(\mM)$   of $\mM$ viewed as an $\rA$--module.
 Hence, we may without loss of generality assume $\mM$ is torsion free, both as
$\rA$-- and $\rB$--module. Moreover, by the universal property of the Macaulayfication functor we obtain a morphism
$\varphi: \mM^{\dagger_\rA} \to \mM^{\dagger_\rB}$ making the following diagram
$$
\xymatrix
{ & \mM \ar[ld]_{\imath} \ar[rd]^\jmath & \\
\mM^{\dagger_\rA} \ar[rr]^\varphi & & \mM^{\dagger_\rB}
}
$$
commutative.  However, the cokernels of $\imath$ and $\jmath$ have finite length over $\rA$, hence
the cokernel of $\varphi$ has finite length, too.  Moreover,
$\varphi$ is a birational isomorphism, hence it is a monomorphism. By Lemma \ref{L:isocodone} it is
an isomorphism.
\qed

\medskip
\begin{proposition}
Let $\kk$ be an algebraically closed field of characteristic zero,
$ (\rA, \idm)$ be a normal $\kk$--algebra of Krull dimension two 
such that $\kk = \rA/\idm$ and 
$\Omega^1_{\rA}$ be the module
of K\"ahler differentials and $\Omega^2_{\rA} = \Omega^1_{\rA} \wedge \Omega^1_{\rA}$.
Then the canonical module $\mK = \mK_A$ is isomorphic to $(\Omega^2_{\rA})^\dagger \cong  (\Omega^2_{\rA})^{**}$.
In the geometrical terms, if
$X = \Spec(\rA)$ and
$U = X \setminus \{\idm\}$ then $\mK  \cong \Gamma(i_*\Omega^2_U)$, where
$\Omega^2_U$ is the sheaf of regular differential two-forms  on the  smooth two-dimensional scheme $U$.
\end{proposition}

\noindent
\emph{Proof}.  Since the scheme $U$ is smooth and two-dimensional,
we have: $\omega_U \cong \Omega_U^2$, where $\omega_U$ is the canonical module of $U$.
Moreover, $\widetilde{\mK}|_U \cong \omega_U$, see \cite{RD}, hence
$\mK = \Gamma(i_* i^* \widetilde{\mK}) \cong  \Gamma(i_*\Omega^2_U)$.
\qed

\medskip
Having this description of $\mK_\rA$ in mind, one may ask about a possible interpretation
of the module $(\Omega^1_{\rA})^{**}$. In turns out that this question is closely related with the
theory of almost split sequences of Cohen-Macaulay modules over
 normal surface singularities.

\begin{proposition}\label{P:fundseq}
Let $(\rA, \idm)$ be a normal two-dimensional Noetherian ring   having a canonical module
$\mK$. Consider the  exact sequence
\begin{equation}\label{E:fund}
0 \to \mK \to \mD \to \rA \to \kk \to 0
\end{equation}
corresponding to a  generator of $\Ext^2_\rA(\kk, \mK) \cong \kk$. Then the module $\mD$ is Cohen-Macaulay.
\end{proposition}

\noindent
\emph{Proof}. From the exact sequence $0 \to \mK \to \mD \to \idm \to 0$ it follows that
$H^0_{\idm}(\mD) = 0$.
Moreover, since $\Ext^1_\rA(\rA, \mK) = 0 = \Ext^2_\rA(\rA, \mK)$, we get:
$\Ext^1_\rA(\idm, \mK) \cong \Ext^2_\rA(\kk, \mK) \cong \kk$.

Let $\DD = \Hom_\rA\bigl(-,\mE(\kk)\bigr)$ be the Matlis functor. If $\rA$ in not complete,
it need not be a duality, but anyway it is exact and
by the Local Duality theorem  we have a commutative diagram
$$
\xymatrix
{
0 \ar[r] & H^1_{\idm}(\mD) \ar[r] \ar[d] & H^1_{\idm}(\idm) \ar[r] \ar[d] & H^2_{\idm}(\mK)  \ar[d] \\
0 \ar[r] & \DD\bigl(\Ext^1_\rA(\mD, \mK)\bigr) \ar[r] &   \DD\bigl(\Ext^1_\rA(\idm, \mK)\bigr) \ar[r] &
\DD\bigl(\Hom_\rA(\mK, \mK)\bigr),
 }
$$
where all horizontal maps are induced by (\ref{E:fund}) and all vertical arrows are isomorphisms.
But the morphism $\Hom_\rA(\mK, \mK) \to \Ext^1_\rA(\idm, \mK) \cong \kk$
is obviously non-zero, hence it is surjective and
the corresponding dual morphism is injective. Thus, $\kk \cong
H^1_{\idm}(\idm) \to H^2_{\idm}(\mK)$
is a monomorphism and  $H^1_{\idm}(\mD) = 0$.
\qed

\begin{remark}
The constructed  exact sequence
(\ref{E:fund})  is called \emph{fundamental}, the corresponding
Cohen-Macaulay module $\mD$ is called \emph{fundamental module} or
\emph{Auslander module}. The reason for this terminology will
be explained below in Remark \ref{R:fundmod}.
\end{remark}

\medskip

\begin{lemma}
Let
$(\rA, \idm)$ be a local normal domain of Krull dimension two having a canonical module
$\mK$.  Then we have:
$(\mD \wedge \mD)^{**} \cong \mK.$
\end{lemma}

\noindent
A proof of this lemma can be found in \cite[Lemma 1.2]{YoshinoKawamoto}.

\medskip

\begin{proposition}[Martsinkovsky]
Let $\kk$ be an algebraically closed field of characteristic zero,
$\rA = \kk\{x,y,z\}/f$ an isolated analytic hypersurface singularity.
Then
$\mD \cong (\Omega_\rA^1)^{**}$ if any only if
$\rA$ is \emph{quasi-homogeneous}.
\end{proposition}

\noindent
For a proof of this result, see \cite[Theorem 1]{Martsinkovsky}.

\medskip

\begin{remark}
If $(\rA, \idm)$ is a quasi-homogeneous analytic normal surface singularity
over an algebraically closed field $\kk$ of characteristic zero
 then $\mD \cong (\Omega_\rA^1)^{**}$, see
\cite[Proposition 2.1]{Behnke} and \cite[Proposition 2.35]{Kahn2}.
It was shown by Herzog \cite{Herzog}
that if the ring $\rA$ is Gorenstein, $\mD \cong (\Omega_\rA^1)^{**}$ and the canonical
morphism $\Omega_\rA^1 \otimes_\rA \kk \to (\Omega_\rA^1)^{**} \otimes_\rA \kk$ is injective, then
$\rA$ is quasi-homogeneous.
\end{remark}

\medskip
\noindent
Moreover, the following conjecture was posed by Martsinkovsky.

\begin{conjecture}
Let $\rA$ be   a
 normal analytic  algebra  of Krull dimension two over an algebraically
closed field of characteristic zero.
Then the isomorphism
$\mD \cong (\Omega_\rA^1)^{**}$ is equivalent to the quasi-homogenity of
  $\rA$.
\end{conjecture}

\medskip

\begin{remark}
Let $(\rS, \idn)$ be a regular Noetherian ring of Krull dimension three,
$(\rA, \idm) = \rS/f$ a normal hypersurface singularity and   $\kk = \rA/\idm$ the residue field.
Then we have:  $\mD \cong  \syz^3(\kk)$, see 
\cite[Lemma 1.5]{YoshinoKawamoto}.
\end{remark}

\medskip
\noindent
The interest to the fundamental module is explained by the following result of Auslander.

\begin{theorem}[Auslander]
Let $\kk$ be an algebraically closed field of characteristic zero,
$(\rA, \idm)$ be an analytic local normal $\kk$-algebra and
 $\mM$ an indecomposable non-free Cohen-Macaulay
$\rA$--module. Then the almost split sequence ending in $\mM$ has the following form:
$$
0 \to \mK \boxtimes_\rA \mM \to \mD \boxtimes_\rA \mM \to \mM \to 0,
$$
in particular, $\tau(\mM) \cong \mK \boxtimes_\rA \mM$.
\end{theorem}

\noindent
For a proof of this Theorem, see \cite[Theorem 6.6]{Auslander} and \cite[Chapter 11]{Yoshino}.

\section{Cohen-Macaulay modules
over two-dimensional quotient singularities}
In this section we deal with  Cohen-Macaulay modules over quotient
two-dimensional singularities. First let us recall  some basic
properties of rings of invariants with respect to an action
of  a finite group.

\begin{theorem}\label{T:quotbasicprop}
Let $(\rA, \idm)$ be a Noetherian local normal domain, $G \subset
\Aut(\rA)$ a finite group of invariants of $\rA$ such that the order of the
group $t = |G|$ is invertible in $\rA$.
 Then the ring of invariants $\rA^G$ is again
\begin{enumerate}
\item a Noetherian local normal domain;
\item
the ring extension $\rA^G \subset \rA$ is finite;
\item  moreover, if the ring
$\rA$ is complete then $\rA^G$ is complete, too;
\item  if $\rA$ is Cohen-Macaulay then
$\rA^G$ is Cohen-Macaulay as well.
\end{enumerate}
\end{theorem}

\noindent
\emph{Proof}. We prove this theorem step by step. Without
loss of generality we assume the action of $G$ is effective, i.e.
for any $g \in G$ there exists $a \in \rA$ such that $g(a) \ne a$.

\hspace{0.2cm}

\noindent  1. First of all we show the ring $\rA^G$ is local and
$\idn = \idm \cap \rA^G$ is its unique maximal ideal. Indeed,
let $x \notin \idn$ then $x$ is invertible in $\rA$. Since its inverse
is again $G$-invariant,  $x$ is invertible in $\rA^G$.
Hence, $\idn$ is the unique maximal ideal of $\rA^G$.

Now we  prove the ring $\rA^G$ is Noetherian. For this it  is sufficient to show
that for any ideal  $I$ in $\rA^G$ we have: $(I\rA)\cap \rA^G = I$.
Indeed, for any increasing chain of ideals
 $$I_1 \subseteq I_2 \subseteq \dots \subseteq I_n \subseteq \dots \subseteq \rA^G$$
 we consider the induced chain
$
I_1\rA  \subseteq I_2 \rA \subseteq \dots \subseteq I_n \rA
 \subseteq \dots \subseteq \rA.
$
Since $\rA$ is Noetherian, $I_m \rA = I_n \rA$ for some big $n$ and all $m \ge n$.
Hence, $I_m = (I_m \rA)\cap \rA^G = (I_n \rA)\cap \rA^G = I_n$.

Let $f_1, f_2, \dots, f_m \in I$ and $r_1, r_2, \dots, r_m \in \rA$ be such that
$f = \sum_{i = 1}^m f_i r_i \in \rA^G$. Then we have:
$$
tf = \sum\limits_{\sigma \in G} \sigma(f) =
\sum\limits_{i = 1}^m \sum\limits_{\sigma \in G} \sigma(f_i r_i) =
t \sum\limits_{i = 1}^m f_i \frac{1}{t} \sum\limits_{\sigma \in G} \sigma(r_i) =
t  \sum\limits_{i = 1}^m f_i \tilde{r}_i,
$$
where $\tilde{r}_i = \frac{\displaystyle 1}{\displaystyle t} \sum_{i=1}^t \sigma(r_i) \in \rA^G$. Hence,
$$
f = \sum_{i = 1}^m f_i \tilde{r}_i
$$
and $(I\rA)\cap \rA^G = I$.

\vspace{0.2cm}
 \noindent
 2.
In order to prove that the ring of invariants
$\rA^G$ is normal observe first that its
ring of fractions $\rK:= Q(\rA^G)$ coincides with $\rL^G$, where
$\rL = \rQ(\rA)$. 
Indeed, one inclusion $\rK = Q(\rA^G) \subset Q(\rA)^G =: \rL$ is clear.
To prove another one,
take any fraction $\frac{\displaystyle a}{\displaystyle b} \in Q(\rA)^G$. Let
$G = \{g_1 = e, g_2, \dots, g_t\}$ then
$$
\frac{a}{b} =
\frac{a g_2(b) \dots  g_t(b)}{b g_2(b)
 \dots  g_t(b)} = \frac{\tilde a}{\tilde b}
$$
where $\tilde b \in \rA^G$. Since $\frac{\displaystyle a}{\displaystyle b}$ is invariant under the
action
of $G$, $\tilde a \in \rA^G$ and $\frac{\displaystyle \tilde a}{\displaystyle \tilde b} \in
Q(\rA^G)$.
Now, assume  $x \in K$ is such that there exist
elements $c_1, c_2,\dots, c_l \in \rA^G$ such that
$x^l + c_1 x^{l-1} + \dots + c_l = 0.$
 Then
$x \in \rA \cap Q(\rA^G) = \rA^G$, hence $\rA^G$ is normal.

To prove the ring extension $\rA^G \subset \rA$ is finite, first note
that  by Artin's Lemma the field extension
$\rK \subset \rL$ is Galois, hence separable, see 
\cite[Theorem VI.1.8]{Lang}.
Let $a \in \rL$ and $\varphi_a(y) = y^t + c_1 y^{t-1} + \dots
+ c_t \in \rK[y]$ be its characteristic polynomial.
Recall that
$\mathrm{tr}(a) = \mathrm{tr}_{\rL/\kK}(a) = - c_1 \in \rK$ is
$\rK$--linear. Moreover, since the  field extension
$\rK \subset \rL$ is separable,
the $\rK$--bilinear form
$$\rL \times \rL \to \rK, \quad  (a,b) \mapsto \mathrm{tr}(ab)$$
is non-degenerate, see  \cite[Theorem VI.5.2]{Lang}.
Since the field extension $\rK \subseteq \rL$ is Galois, 
 the characteristic  polynomial of
$a \in \rA$ is
$$\varphi_a(y) = \bigl(y - h_1(a)\bigr)\bigl(t- h_2(a)\bigr)\dots 
\bigl(y - h_t(a)\bigr),$$ where
$G= \{h_1, h_2,\dots, h_t\}$. In particular,
$\mathrm{tr}(a) \in \rA^G$.

Let $\alpha_1, \alpha_2,\dots,\alpha_n \in \rL$ be a basis
of $\rL$ over $\rK$. Note that without loss of generality we may assume
all elements $\alpha_i$ actually belong to  $\rA$.
Denote  $$\rA^{!} 
:= \bigl\{x \in \rL \, \big| \, \mathrm{tr}(xa) \in \rA^G \,\, \mbox{for all} \, \, 
\, \, a \in \rA\bigr\}.$$
Since $\mathrm{tr}$ is $\rA^G$--linear, $\rA^{!}$ is an $\rA^G$--module;
moreover, $\rA \subseteq \rA^{!}$. Let $$\mM :=
\langle \alpha_1, \alpha_2, \dots, \alpha_n\rangle_{\rA^G} \subset
\rA.$$ It is easy to see that  $\rA^{!} \subseteq \mM^{!} :=
\langle \alpha_1^*, \alpha_2^*, \dots, \alpha_n^*\rangle_{\rA^G},$
where $\alpha_i^*(\alpha_j) = \delta_{i,j}$.
Summing everything up, we see that $\rA$ is a submodule of
a finitely generated $\rA^G$--module $\mM^{!}$. Since $\rA^G$ is
Noetherian, $\rA$ is a Noetherian $\rA^G$--module, hence finite
over $\rA^G$.

 \vspace{0.2cm}
 \noindent
 3. Assume the ring $\rA$ is complete. We already know the ring
 $\rA^G$ is local and $\idn = \idm \cap \rA^G$ is its unique
 maximal ideal.
In order to show the local ring
$(\rA^G,\idn)$ is complete, take any sequence
$(a_n)_{n \ge 1}$ of elements of $\rA^G$  such that
$a_n - a_m \in \idn^l$ for all $n,m \ge 1$, where
$l = \min(m,n)$. By completeness of  $\rA$  there exists
a unique element $a \in \rA$ such that $a \equiv a_l \, \,
 \mod \,\, \idm^l$ for all $l \ge 1$. Since for all
 $g \in G$ we have $g(a_l) = a_l$, by Krull's intersection
 theorem
 $$
 g(a) - a \in
 \bigcap\limits_{l \ge 1} \idm^l = 0.
 $$

 \vspace{0.2cm}
 \noindent
 4. Since  $t = |G|$ is invertible in $\rA$, we can consider the
 Reynold's operator $p: \rA \to \rA^G$, given by the rule
 $a \mapsto p(a) := \frac{1}{|G|}
 \sum_{g \in G} g(a)$. It is clear that $p$ is $\rA^G$--linear
 and $p(a) = a$ for $a \in \rA^G$,
 hence $\rA \cong  \rA^G \oplus \rA'$ viewed as  an
 $\rA^G$--module. Moreover,
 since the ring extension $\rA^G \subset \rA$ is known to be finite,
 $H^t_{\idn}(\rA)= H^t_{\idm}(\rA)$ for all $t \ge 0$.
 If $\rA$ is Cohen-Macaulay, then by Corollary \ref{C:CMandloccoh} we
 have
 $H^t_{\idn}(\rA) = 0$ for $t \ne d$. Thus,  $H^t_{\idn}(\rA^G) = 0$
 and $\rA^G$ is Cohen-Macaulay, too.
 \qed

\medskip

\begin{remark}\label{R:funnyremark}
If the order of the group $G$ is not invertible in $\rA$, then the ring of invariants
 $\rA^G$ can be \emph{not} Noetherian! Such an example was constructed for the
first time by Nagata in  \cite{Nagata}. Moreover, Fogarty gave an example
of a finite group $G$ acting on  a Cohen-Macaulay ring $\rA$ such that
the ring of invariants $\rA^G$ is Noetherian but not Cohen-Macaulay, see \cite{Fogarty}.
\end{remark}

\medskip

\begin{remark}
If $(\rA, \idm)$ is a Noetherian local ring and $G$ a finite group of automorphisms of $\rA$
such that $|G|$ is invertible in $\rA$ and $\rA^G$ is Cohen-Macaulay, then
$\rA$ is \emph{not} necessary Cohen-Macaulay. Indeed, let $\kk$ be a field of
characteristic different from two, $\rA = \kk\llbracket x,y\rrbracket /(xy, y^2)$ and
$G = \langle \sigma\,|\, \sigma^2 = e\rangle$ acts on $\rA$  $\kk$-linearly
by the rule $\sigma(x) = x$ and $\sigma(y) = -y$. Then
$\rA^G = \kk\llbracket x\rrbracket $ is Cohen-Macaulay, but $\rA$ is not.
\end{remark}

\medskip

Let $\kk$ be an algebraically closed field,
$\rR = \kk\llbracket x_1,x_2,\dots,x_n\rrbracket $ the ring of formal power series, $\idn$ its maximal ideal
 and $G \subset \Aut(\rR)$  a finite group of ring automorphisms of $\rR$ such that $|G|$
 is invertible in $\kk$. Then  $G$ acts on the cotangent  space  $V:= \idn/\idn^2$ and
we obtain a  group homomorphism $\rho: G \to
\GL(V)$.

\begin{proposition}\label{P:Cartan}
In the notations as above denote  $G' = \im(\rho)$.
Then we have: $R^G \cong R^{G'}$.
\end{proposition}

\noindent
\emph{Proof}. The following argument is due to Cartan \cite{Cartan}.
For an element $g \in G$ let $\bar{g}$ be the corresponding element of
$G'$. Let
$$
\bar{g}x_i = \sum\limits_{j=1}^n \bar{g}_{ij} x_j
$$
be the linearized action of $G$ on $\rR$.
Consider the ring endomorphism
$\tau: \rR \to \rR$
given by the rule
$$
\tau(x_i) = \frac{1}{|G|}
\sum\limits_{g \in G} \bar{g}^{-1} g (x_i), \quad  1 \le i \le n.
$$
Since the induced  action of $\tau$
 on the cotangent space $\idn/\idn^2$ is identity,
it is an automorphism. Moreover, it is easy to see that for
any element $h \in G$ we have: $\bar{h} \circ \tau = \tau \circ h$,
so these two actions are conjugate.  Therefore, they have isomorphic rings
 of invariants, namely, $a\mapsto \tau(a)$ gives an isomorphism
 $\rR^{G}\to\rR^{G'}$.
\qed

\medskip

\begin{remark}\label{R:Cartan}
Proposition \ref{P:Cartan} remains true if we replace $\rR$ by the ring
of convergent power series $\mathbb{C}\{x_1, x_2,\dots,x_n\}$. Moreover,
it was shown by Cartan that if $\rA = \mathbb{C}\{x_1, x_2,\dots,x_n\}/I$ is an
analytic algebra over $\mathbb{C}$ and $G$ a finite group of invariants of $\rA$ then
$\rA^G$ is again analytic, \linebreak see \cite[Th\'eor\`eme 4]{Cartan}.
\end{remark}

\medskip

In other words, Proposition \ref{P:Cartan} and Remark \ref{R:Cartan} mean
that dealing with  quotient singularities
we  may  without loss of generality
assume that group  $G$ acts
 linearly on $\rR$. Moreover, it is sufficient
to consider the so-called
small subgroups.

\medskip

\begin{definition}
Let $G \subset \GL_n(\kk)$ be a finite subgroup. An element $g \in G$ is a
 \emph{pseudo-reflection} if
$\rank(1 - g) = 1$.
A subgroup $G$ is called \emph{small} if it contains no pseudo-reflections.
\end{definition}

\medskip

\begin{example}
Let $\rR = \kk\llbracket x_1, x_2\rrbracket $ and
$G = \mathbb{Z}_2 = \langle \sigma\rangle$, where $\sigma(x_1) = x_2$ and
$\sigma(x_2) = x_1$. Then $\sigma$ is a pseudo-reflection and
$\rR^G = \kk\llbracket x_1 + x_2, x_1 x_2\rrbracket  \cong \kk\llbracket u,v\rrbracket $. Hence, if a finite group
is not small, one can obtain   a smooth ring as a quotient.
\end{example}

\medskip
\noindent
The following theorem lists  some basic results on  quotient singularities.

\begin{theorem}
Let $\kk$ be an algebraically closed field of characteristic zero,
$\rR = \kk\llbracket x_1, x_2, \dots,x_n\rrbracket $ and $G \subset \GL_n(\kk)$  be
a finite subgroup. Then the following properties
hold:
\begin{enumerate}
\item  The invariant ring $\rR^G$ is always normal and Cohen-Macaulay.
\item  There exists a finite small group $G' \subset \GL_n(\kk)$ such that
$\rR^G \cong  \rR^{G'}$.
\item  Let $G_1$ and $G_2$ be two finite small subgroups of $\GL_n(\kk)$. Then
      $\rR^{G_1} \cong  \rR^{G_2}$ if and only if there exists an element $g \in \GL_n(\kk)$ such that
      $g^{-1} G_1 g = G_2$.
\item Let $G \subset \GL_n(\kk)$ be a finite small subgroup. Then $\rR^G$ is Gorenstein if and only if
      $G \subset \SL_n(\kk)$.
\end{enumerate}
\end{theorem}

\noindent
The property (1) follows is proven in Theorem \ref{T:quotbasicprop},
(2) and (3) are due to Prill, see
\cite[Proposition 6]{Prill} and
\cite[Theorem 2]{Prill}, respectively. Finally, the part (4)
is a result of Watanabe \cite{WatanabeI} and \cite{WatanabeII}, see also \cite{Hinich}.

\medskip

Let $\rR = \mathbb{C}\llbracket x_1, x_2\rrbracket $. All finite subgroups of
$\SL_2(\mathbb{C})$ modulo conjugation are known
and the corresponding quotient singularities are precisely the \emph{simple hypersurface singularities}.

\begin{enumerate}
\item If $G = \langle g| g^n = e\rangle \cong  \mathbb{Z}_n$ is a cyclic subgroup of order $n \ge 2$, where
$$
g =
\left(
\begin{array}{cc}
\xi & 0 \\
0 & \xi^{-1}
\end{array}
\right), \quad \xi = \exp\Bigl(\frac{\pi i}{n}\Bigr)
$$
then the corresponding ring of invariants is
$$\mathbb{C}\llbracket x_1^n, x_2^n, x_1 x_2\rrbracket  \cong \mathbb{C}\llbracket u,v,w\rrbracket /(uv - w^n).$$
 This is the so-called
simple singularity of type  $A_{n-1}$.

\item Binary dihedral group
$\DD_{n}$  is generated by two elements $a$ and $b$, where
$$
a =
\left(
\begin{array}{cc}
\xi & 0 \\
0 & \xi^{-1} \\
\end{array}
\right), \quad
\xi = \exp\Bigl(\frac{\pi i}{n}\Bigr),\quad
b =
\left(
\begin{array}{cc}
0 & 1 \\
-1  & 0 \\
\end{array}
\right).
$$
It corresponds to the $D_n$--singularity ($n \ge 4$), given by the equation
$$
u^{n-1} + uv^2 + w^2.
$$
\item Binary tetrahedral group $\TT$ is generated by three elements
$\sigma, \tau$ and
$\mu$ where
$$
\sigma =
\left(
\begin{array}{cc}
i & 0 \\
0 & -i \\
\end{array}
\right), \quad
\tau =
\left(
\begin{array}{cc}
0 & 1 \\
-1  & 0 \\
\end{array}
\right), \quad
\mu =
\frac{1}{\sqrt{2}}
\left(
\begin{array}{cc}
\xi^7 & \xi^7 \\
\xi^5 & \xi \\
\end{array}
\right), \quad \xi = \exp\Bigl(\frac{2\pi i}{8}\Bigr).
$$
It corresponds to $E_6$-singularity
$$u^3 + v^4 + w^2 = 0.$$
\item Binary octahedral group
$\OO$ is generated by the matrices  $\sigma, \tau,
\mu$  occurring in the description  of $\TT$ and by
$$
\kappa =
\left(
\begin{array}{cc}
\xi & 0 \\
0 & \xi^7 \\
\end{array}
\right).
$$
The corresponding singularity is
$$E_7: \quad u^3 v + v^3 + w^2.$$
\item Finally,  we have the binary icosahedral subgroup
$\II = \langle \sigma, \tau\rangle$, where
$$
\sigma = -
\left(
\begin{array}{cc}
\xi^3 & 0 \\
0 & \xi^2 \\
\end{array}
\right), \quad
\tau = \frac{1}{\sqrt{5}}
\left(
\begin{array}{cc}
-\xi + \xi^4 & \xi^2 - \xi^3 \\
\xi^2 - \xi^3 & \xi - \xi^4\\
\end{array}
\right), \quad \xi = \exp\Bigl(\frac{2\pi i}{5}\Bigr).
$$
The corresponding singularity is
$$E_8: \quad u^3 + v^5 + w^2.$$
\end{enumerate}

\medskip
\begin{remark}
In fact, all finite subgroups of $\SL_2(\mathbb{C})$ can be
described in the following elegant way: they are parameterized
by triples of positive integers
$$\Bigl\{(p, q, r) \, \Big| \,  p \le q \le r, \quad  \frac{1}{p} +
\frac{1}{q} + \frac{1}{r} > 1\Bigr\}$$ and are given by the
presentation
$$
G_{p,q,r} = \Bigl\langle x,y,z \, \Big| \,  x^p = y^q = z^r = xyz\Bigr\rangle.
$$
Moreover,
\begin{itemize}
\item If $p = 1$ then $G_{1,q,r} \cong \mathbb{Z}_{q+r}$;
\item $G_{2,2,n}$ is the binary dihedral group $\mathbb{D}_n$;
\item $G_{2,3,3}$ is the binary tetrahedral group $\mathbb{T}$;
\item $G_{2,3,4}$ is the binary octahedral  group $\mathbb{O}$;
\item $G_{2,3,5}$ is the binary icosahedral  group $\mathbb{I}$.
\end{itemize}
\end{remark}

\medskip
\begin{remark}
A classification of all small subgroups of $\GL_2(\mathbb{C})$ modulo conjugation is
also known, see for example \cite[Satz 2.3]{Brieskorn}.
\end{remark}

\medskip
\begin{theorem}[Herzog \cite{Herzog1}]\label{T:Herzog}
Let $\kk$ be an algebraically closed field
of characteristic zero, $\rR = \kk\llbracket x_1, x_2\rrbracket $ or
$\kk\{x_1, x_2\}$,  $G \subseteq \GL_2(\kk)$ be a finite
 small  subgroup and
$\rA = \rR^G$. Then $\CM(\rA) = \add_\rA(\rR)$, i.e.
any  indecomposable
Cohen-Macaulay $\rA$--module is a direct summand of $\rR$ viewed as as $\rA$--module.
\end{theorem}

\noindent
\emph{Proof}.  The inclusion map $\imath: \rA \to \rR$ has an inverse:
$p: \rR \to \rA$ given by the formula
$$
p(r) =
\frac{1}{|G|} \sum\limits_{g \in G} g(r).
$$
Moreover, $p$ is a morphism of $\rA$--modules, hence $\rR \cong \rA \oplus \rR'$, where
$\rR'$ is a certain $\rA$--module.
Hence, for any  $\rA$--module $\mM$ we have:
$$
\rR \otimes_\rA \mM \cong \mM \oplus \mN,
$$
where $\mN = \rR \otimes_\rA \rR'$. This implies:
$$
(\rR \otimes_\rA \mM)^\dagger \cong \mM \oplus \mN^\dagger
$$
in the category of $\rA$--modules.
Since $\rR \otimes_\rA \mM$ is also an $\rR$--module, by Proposition \ref{P:Macdiffbase}
$(\rR \otimes_\rA \mM)^{\dagger_\rR} \cong (\rR \otimes_\rA \mM)^{\dagger_\rA}$ as $\rA$--modules.
Moreover,  $\rR$ is regular, hence by Corollary \ref{C:CMreg} we have:
 $(\rR \otimes_\rA \mM)^{\dagger_\rR} \cong \rR^n$ for
some positive integer $n$. Since the category $\rA-\mod$ is Krull-Schmidt,
$\mM$ is a direct summand of $\rR$.
\qed

\medskip
Theorem \ref{T:Herzog} shows that a quotient surface singularity always has
finite Cohen-Macaulay representation type. However, one would wish a more explicit description
of indecomposable  Cohen-Macaulay modules.

\medskip
\noindent
Recall the following easy fact from category theory.

\begin{lemma}\label{L:addclosure}
Let $\cA$ be an additive category, $\mX$ an object of $\cA$ and $\Gamma = \End_\cA(\mX)$ its endomorphism ring.
Then the functor $\Hom_\cA(\mX, \,-): \quad \cA \lar  \mod-\Gamma$ induces an equivalence of categories
$$
\Hom_\cA(\mX, \,-): \quad \add(\mX)  \lar \Pro(\Gamma),
$$
where $\Pro(\Gamma)$ is the category of finitely generated
projective right $\Gamma$--modules.
\end{lemma}

\medskip
Since by Herzog's result we know that $\CM(\rA) = \add_\rA(\rR)$ for a quotient singularity $\rA = \rR^G$, it raises a question
about a description of the ring $\End_\rA(\rR)$.

\begin{definition}
Let $G$ be a finite group acting on a ring $\rS$. Then the skew group ring
$\rS*G$ is a free left $\rS$-module
$$
\rS*G = \Bigl\{\sum\limits_{g \in G} a_g [g] \, \big| \, a_g \in \rS\Bigr\}
$$
and the multiplication  is given by the rule : $a_g[g]\cdot a_h[h] = a_g g(a_h) [gh]$.
\end{definition}

\medskip
\begin{theorem}[Auslander]\label{T:endquot}
Let $\kk$ be an algebraically closed field of characteristic zero,
$G \subset \GL_2(\kk)$ a small subgroup, $\rR =
\kk\llbracket x_1, x_2\rrbracket $ and $\rA  = \rR^G$.
 Then the $\kk$-linear map
$$
\theta:  \rR*G \lar \End_\rA(\rR)
$$
mapping an element $s[g]$ to the morphism $r \mapsto s g(r)$,  is an isomorphism of algebras.
\end{theorem}

\noindent
We refer to \cite{Auslander}  (see also \cite[Proposition 10.9]{Yoshino}) for a proof of this
theorem.

\medskip
\noindent
Since $\chr(\kk) = 0$, the Jacobson's radical of $\kk\llbracket x_1,x_2\rrbracket *G$ is
$\idm*G$ and the semi-simple algebra
$\rR*G/\rad(\rR*G)$ is isomorphic to the group algebra of $G$.
Moreover, the functor
$$
\rR*G-\mod \lar \bigl(\rR*G/\rad(\rR*G)\bigr)-\mod
$$
mapping an $\rR*G$--module $\mM$ to $\mM/\rad(\mM)$ induces  a bijection between indecomposable
projective $\rR*G$--modules and irreducible representations of the group $G$.

\medskip

\begin{theorem}[Auslander]\label{T:CMquot}
Let $G \subset \GL_2(\kk)$ be a small subgroup,
$\rR = \kk\llbracket x_1, x_2\rrbracket $ and $\rA  = \rR^G$.
Then the exact functor
$
\mathbb{F}: \quad \rR*G-\mod \lar \rA-\mod
$
mapping an $\rR*G$--module $\mM$ to
$\mM^G = \bigl\{m \in \mM\,\big| \, gm = m \quad 
\mbox{for all} \quad g \in G\bigr\}$
and a morphism of  $\rR*G$--modules $f: \mM \to \mN$ to $f|_{\mM^G}$, induces an equivalence of categories
$
\Pro(\rR*G) \lar \CM(\rA).
$
\end{theorem}

\noindent
\emph{Proof}. Let us check the functor $\mathbb{F}$ is exact.
Let
$$
0 \to \mM \stackrel{\varphi}\lar \mN \stackrel{\psi}\lar \mK \to 0
$$
be an exact sequence of $\rR*G$--modules. It is clear that
 $\varphi^G$ is injective and $\ker(\psi^G) = \im(\varphi^G)$.  Let us check $\psi^G$ is an epimorphism.
For  $z \in \mK^G$ take $y \in \mN: \psi(y) = z$ and put
$\tilde{y} = \frac{1}{|G|} \Bigl(\sum_{g \in G} g\cdot y\Bigr) \in \mN^G$.
Then $\psi(\tilde{y}) = z$, hence $\mathbb{F}$ is exact.

\hspace{0.1cm}

\noindent
Next, it is easy to see that
$
\mathbb{F}(\rR*G) = (\rR*G)^G =  \rS_1 :=
\Bigl\{\tilde{r} = \sum\limits_{g \in G} g(r)[g]\, \big| \,  r \in \rR\Bigr\}
$
is isomorphic to $\rR$ as an $\rA$--module.  Moreover,
in the category of left $\rR*G$--modules we have:
$$\End_{\rR*G}(\rR*G) \cong ({\rR*G})^{\mathsf{op}}.$$
Consider the following chain of morphisms of algebras:
$$
\rR*G \stackrel{\alpha}\lar (\rR*G)^{\mathsf{op}} \stackrel{\beta}\lar \End_{\rR*G}(\rR*G)
\stackrel{\gamma}\lar \End_{\rA}(\rS_1)
$$
where $\alpha(p[g]) = g^{-1}(p)[g^{-1}]$, $\beta(\eta)(\zeta) = \zeta \eta$ and $\gamma$ is the morphism induced
by $\mathbb{F}$.
Let $r \in \rR$ and $\tilde{r} = \sum\limits_{h \in G} h(r)[h]$. Then
$$\bigl((\gamma \circ \beta \circ \alpha)(p[g])\bigr)(\tilde{r}) =
\bigl(\sum\limits_{h \in G} h(r)[h])(g^{-1}(p)[g^{-1}]\bigr) =
\sum\limits_{h \in G} h(r g^{-1}(p)) h g^{-1} =
$$
$$
= \sum\limits_{h \in G} hg^{-1}(g(r)p) hg^{-1} = \widetilde{p g(r)}.
$$
Hence, the morphism of $\kk$--algebras $\gamma \circ \beta \circ \alpha$ coincides with the morphism
$\theta$ from Theorem \ref{T:endquot}, which is known to be an isomorphism.
Since  $\mathbb{F}(\rR*G) \cong \rR$ and
$\mathbb{F}: \End_{\rR*G}(\rR*G) \to \End_\rA(\rR)$ is an isomorphism,
by Theorem \ref{T:Herzog} we get  a chain of equivalences
of categories
$$
\Pro(\rR*G) = \add_{\rR*G}(\rR*G) \cong \add_{\rA}(\rR) \cong \CM(\rA).
$$
\qed

\medskip
Let $\mV$ be an irreducible representation of the group $G$, then
$\mP= \mP_\mV := \kk\llbracket x_1,x_2\rrbracket  \otimes_\kk \mV$
is an  indecomposable projective left $\rR*G$--module such that
$\mP/\rad(\mP) \cong \mV$, where $p[g] \in \rR*G$ acts on a simple tensor $q \otimes v$ by the rule
$p[g](q \otimes v) = p g(q) \otimes g(v)$.

\begin{corollary}[McKay Correspondence \`a la Auslander]\label{C:McKay}
Let $G \subset \GL_2(\kk)$ be a finite small subgroup,
$\rR = \kk\llbracket x_1,x_2\rrbracket $ and $\rA = \rR^G$.
Then there exists a bijection
between irreducible representations of $G$ and indecomposable
Cohen-Macaulay $\rA$--modules given by the
functor
$$
\Rep(G) \ni \mV  \mapsto
\bigl(\kk\llbracket x_1,x_2\rrbracket  \otimes_\kk \mV\bigr)^G \in \rA-\mod.
$$
\end{corollary}

\medskip

\begin{remark}\label{R:fundmod}
In the notations of Corollary \ref{C:McKay},
let $V  = \langle e_1, e_2\rangle_\kk \cong  \kk^2$ be a two-dimensional
vector space. The given embedding
$\rho: G \subset \GL(V) = \GL_2(\kk)$
defines a two-dimensional representation
$g \in G \mapsto \rho(g)$ called the \emph{fundamental representation}.
By a result of Watanabe \cite{WatanabeI} and \cite{WatanabeII}, see also
\cite{Hinich}, the Cohen-Macaulay $\rA$--module
$\mK = (\kk\llbracket x_1, x_2\rrbracket  \otimes_{\kk} V_{\Det})^G$ is canonical, where
$V_{\Det} = \kk$ and $\Det(g) := \Det(\rho(g))$. Note that we have
$V_{\Det} = \wedge^2(V)$.
For $f = \alpha_1 e_1 + \alpha_2 e_2 \in V$  denote
$\tilde f = \alpha_1 x_1 + \alpha_2 x_2$.
Then  the Koszul resolution
$$
0 \to \kk\llbracket x_1, x_2\rrbracket  \otimes \wedge^2(V)
\stackrel{\alpha}\lar \kk\llbracket x_1, x_2\rrbracket  \otimes V
\stackrel{\beta}\lar  \kk\llbracket x_1, x_2\rrbracket  \stackrel{\phi}\lar \kk \to 0
$$
where $\alpha\bigl(p \otimes (f_1\otimes f_2 -
f_2 \otimes f_1)\bigr) = p \tilde{f}_1 \otimes f_2 -
p \tilde{f}_2 \otimes f_1, \quad \beta(q \otimes f) = q \tilde{f}$ and
$\phi(t) = t(0,0)$ is also a minimal free resolution of
$\kk$ in the category of $\kk\llbracket x_1, x_2\rrbracket *G$--modules.
Hence, taking $G$--invariants we obtain an exact sequence
$$
\omega: \quad
0 \to \mK \to \mD \to \rA \to \kk \to 0,
$$
where $\omega$ denotes the corresponding element in
$\Ext^2_\rA(\kk, \mK) \cong \kk$. Moreover, since the morphism
$\alpha: \kk\llbracket x_1, x_2\rrbracket  \otimes \wedge^2(V) \stackrel{\alpha}\lar
\kk\llbracket x_1, x_2\rrbracket \otimes V$ is non-split, the sequence $\omega$
is non-split, too, and $\omega \ne 0$. Hence, the 
module
$\mD  = \bigl(\kk\llbracket x_1,x_2\rrbracket  \otimes V\bigr)^G$ is the fundamental module
of the quotient singularity $\rA = \kk\llbracket x_1, x_2\rrbracket ^G$.
\end{remark}

\medskip
\begin{remark}\label{R:ARs}
Let $W$ be a non-trivial  irreducible $\kk[G]$--module.
Then its minimal free
projective resolution as $\kk\llbracket x_1, x_2\rrbracket *G$--module is
$$
0 \to \kk\llbracket x_1, x_2\rrbracket  \otimes (\wedge^2(V)\otimes W)
\stackrel{\alpha \otimes 1}\lar \kk\llbracket x_1, x_2\rrbracket  \otimes (V \otimes W)
\stackrel{\beta \otimes 1}\lar  \kk\llbracket x_1, x_2\rrbracket  \otimes W \lar W  \to 0.
$$
Since morphisms $\alpha \otimes 1$ and $\beta \otimes 1$ are almost split
in the category of projective $\kk\llbracket x_1, x_2\rrbracket *G$--modules,  $W^G = 0$ and
the functor $\Pro(\kk\llbracket x_1, x_2\rrbracket *G) \lar \CM(\rA)$ is
an equivalence of categories, we obtain an exact
sequence of Cohen-Macaulay $\rA$--modules
$$
0 \to \bigl(\kk\llbracket x_1, x_2\rrbracket  \otimes (\wedge^2(V)\otimes W)\bigr)^G
\stackrel{\alpha \otimes 1}\lar
\bigl(\kk\llbracket x_1, x_2\rrbracket  \otimes (V \otimes W)\bigr)^G
\stackrel{\beta \otimes 1}\lar
\bigl(\kk\llbracket x_1, x_2\rrbracket  \otimes W\bigr)^G   \to 0,
$$
which is precisely the almost split sequence ending
at the indecomposable  Cohen-Macaulay module
$\bigl(\kk\llbracket x_1, x_2\rrbracket \otimes W\bigr)^G$.
\end{remark}

\medskip

\begin{corollary}
Let $W$ be an irreducible representation of $G$ over $\kk$
and $\mM_W = \bigl(\kk\llbracket x_1, x_2\rrbracket \otimes  W\bigr)^G$  the corresponding
indecomposable Cohen-Macaulay $\rA$--module. Then
$\rank(\mM_W) = \dim_\kk(W)$.
\end{corollary}

\medskip

\begin{example}\label{E:quot}
Let $\kk$ be an algebraically closed field of characteristic zero,
$G = \langle g\rangle = \mathbb{Z}_n$ be a cyclic group of order $n \ge 2$,
$\varepsilon \in \kk$ a primitive $n$-th root of unity  and
$0 <  m < n$ an integer such that $\gcd(n,m) = 1$. Define
an embedding $\rho: G \to \GL_2(\kk)$ by the rule
$$
\rho(g) =
\left(
\begin{array}{cc}
\varepsilon & 0 \\
0 & \varepsilon^m
\end{array}
\right).
$$
Let $\rA = \kk\llbracket x_1, x_2\rrbracket ^G$ be the corresponding quotient
singularity. Since the group $G$ is cyclic, all irreducible
representations of $G$ are one-dimensional. Moreover, there
are precisely  $n$ non-isomorphic
irreducible representations $V_1, V_2, \dots, V_n$ of  $G$ defined  by the action
$g \cdot 1 = \varepsilon^{-l}, \, \, 1 \le l \le n$. By Theorem \ref{T:CMquot},
the corresponding indecomposable Cohen-Macaulay $\rA$--modules are
$$
\kk\llbracket x_1, x_2\rrbracket  \supseteq \mM_l := 
\Bigl\{\sum\limits_{i,j = 0}^\infty a_{ij} x_1^i x_2^j \, \Big| \, 
\,\, a_{ij} \in \kk, \, \, i+ mj \equiv
l \, \,  \mod  \, \, n\Bigr\}, \quad 1 \le l \le n.
$$
\end{example}

\medskip
\begin{remark}
In the notations of Example  \ref{E:quot},  put  $m = n - 1$.
Then the corresponding invariant ring $\rA =
\kk\llbracket x_1, x_2\rrbracket ^G$ is
the simple singularity $\kk\llbracket u,v,w\rrbracket /(uv - w^n)$ of type
$A_{n-1}, \, \, n \ge 2$.
In this case all non-free indecomposable Cohen-Macaulay
$\rA$--modules are given by the matrix factorizations
$\mM = \mM(\varphi_l, \psi_l)$, where
$$
\varphi_l =
\left(
\begin{array}{cc}
u & - w^{n-l} \\
- w^l & v
\end{array}
\right),
\quad
\psi_l =
\left(
\begin{array}{cc}
v & w^{n-l} \\
w^l & u
\end{array}
\right),
\quad 1 \le l < n.
$$
Matrix factorizations describing indecomposable Cohen-Macaulay modules
over the other two-dimensional simple hypersurface singularities
can be found for example in
\cite{Kajiuraetal}.
\end{remark}

By Theorem \ref{T:Herzog} all quotient singularities have
finite Cohen-Macaulay representation type. Moreover, Auslander \cite{Auslander}
and Esnault \cite{Esnault} have  shown that
the converse is true  in the case $\kk = \mathbb{C}$.

\begin{theorem}\label{T:finisquot}
Let $(X,x)$ be a normal complex two-dimensional singularity
of finite Cohen-Macaulay representation type. Then $(X,x)$ is
a quotient singularity.
\end{theorem}

\noindent
A proof of this theorem can be found in \cite[Chapter 11]{Yoshino}.

\medskip
It is interesting to note that the normal complex two-dimensional
singularities of finite Cohen-Macaulay representation type can be described
in purely topological terms.

\begin{theorem}[Prill]\label{T:Prill}
Let $(X,x)$ be a complex normal  surface singularity. Then the following
statements are equivalent:
\begin{itemize}
\item $(X, x)$ is a quotient singularity;
\item The local fundamental group of $(X, x)$ is finite.
\end{itemize}
\end{theorem}

\noindent
For a proof of this theorem we refer to
\cite[Theorem 3]{Prill}, see also \cite[Satz 2.8]{Brieskorn}.

\medskip
\noindent
Combining theorems of Prill and Auslander, we obtain the
following interesting corollary.

\begin{corollary}\label{C:fintypeisfingroup}
 Let $(X,x)$ be a complex normal surface singularity. Then the following
statements are equivalent:
\begin{itemize}
\item $(X, x)$ is a quotient singularity;
\item The local fundamental group of $(X, x)$ is finite.
\item $(X, x)$ has finite Cohen-Macaulay representation type.
\end{itemize}
\end{corollary}

\section{Cohen-Macaulay modules over certain
non-isolated singularities}

Let $\kk$ be an algebraically closed field of arbitrary characteristic. 
The aim of this section is to classify indecomposable Cohen-Macaulay
modules over non-isolated singularities
$\kk\llbracket x,y,z\rrbracket /xy$ (type $A_\infty$)
and $\kk\llbracket x,y,z\rrbracket /(x^2 y - z^2)$ (type $D_\infty$).
By a result of Buchweitz, Greuel and Schreyer these singularities
have \emph{discrete} (also called  \emph{countable}) Cohen-Macaulay
representation type. Note that  in their paper \cite{BGS}
the authors obtain a complete classification of indecomposable
Cohen-Macaulay modules over the corresponding  \emph{curve} singularities
$\kk\llbracket x,y\rrbracket /x^2$ and $\kk\llbracket x,y\rrbracket /x^2y$. Although, by Kn\"orrer's
correspondence \cite{Knoerrer2} the surface
singularities $A_\infty$ and $D_\infty$ have the same Cohen-Macaulay
representation type as the corresponding curve singularities, the
problem to describe  all indecomposable  matrix factorizations
for $\kk\llbracket x,y,z\rrbracket /xy$
and $\kk\llbracket x,y,z\rrbracket /(x^2 y - z^2)$
remained to be done.

Let $\kk$ be an algebraically closed field and
$(\rA, \idm)$ be a reduced analytic Cohen-Macaulay $\kk$--algebra
of Krull dimension two, which is not an isolated singularity,
 and  let  $\rR$  be its normalization.
 It is well-known that $\rR$ is again analytic and the
ring extension $\rA \subset \rR$ is finite, see \cite{GLS} or
\cite{deJongPfister}.  Moreover, the
ring $\rR$ is isomorphic to the  product of a finite number of normal
local rings:
$$
\rR \cong  (\rR_1, \idn_1) \times (\rR_1, \idn_1) \times \dots \times
(\rR_t, \idn_t).
$$
 Note that all rings $\rR_i$ are automatically Cohen-Macaulay, see
\cite[Theorem IV.D.11]{Serre}.

Let $I = \Ann(\rR/\rA) \cong \Hom_\rA(\rR, \rA)$
be the conductor ideal. Note that $I$ is also an ideal in $\rR$, denote
$\bar\rA = \rA/I$ and $\bar\rR = \rR/I$.
By Lemma \ref{L:CMdim2} the ideal $I$ is Cohen-Macaulay, both  as $\rA$- and
$\rR$-module. Moreover, $V(I) \subset \Spec(\rA)$ is exactly the locus
where the ring $\rA$ is not normal. It is not difficult to show
that both rings $\bar\rA$ and $\bar\rR$ have Krull dimension one and
are Cohen-Macaulay (but not necessary reduced). Let
$\rQ(\bar\rA)$ and $\rQ(\bar\rR)$ be the corresponding total rings
of fractions, then the inclusion $\bar\rA \to \bar\rR$ induces
an inclusion $\rQ(\bar\rA) \to \rQ(\bar\rR)$.

Let $\mM$ be a Cohen-Macaulay $\rA$--module. Recall
that $ \rR \boxtimes_\rA \mM = (\rR \otimes_\rA \mM)^\dagger$
and for any Noetherian $\rR$--module $\mN$ we have an exact sequence
$$
0 \to \tor(\mN) \to \mN \stackrel{i_\mN}\lar  \mN^\dagger \to \mT \to 0,
$$
where $\mT$ is an $\rR$--module of finite length.
Hence,  the canonical morphism
$$
\theta_\mM: \rQ(\bar\rR) \otimes_{\rA} \mM  =
\rQ(\bar\rR) \otimes_{\rQ(\bar\rA)} \bigl(\rQ(\bar\rA)
\otimes_{\rA} \mM\bigr)
\lar \rQ(\bar\rR) \otimes_{\rR} \bigl(\rR  \otimes_{\rA} \mM\bigr)
 \lar
 \rQ(\bar\rR) \otimes_{\rR} \bigl(\rR  \boxtimes_{\rA} \mM\bigr)
$$
is an epimorphism. Moreover, one can show that the canonical morphism
$$\eta_\mM:  \rQ(\bar\rA) \otimes_{\rA} \mM \to
\rQ(\bar\rR) \otimes_{\rA} \mM \stackrel{\theta_\mM}\lar
\rQ(\bar\rR) \otimes_{\rR} \bigl(\rR  \boxtimes_{\rA} \mM\bigr)
$$
is a monomorphism provided  $\mM$ is Cohen-Macaulay.

\begin{definition}\label{D:triples}
In the notations of this section,  consider the following
\emph{category of triples} $\Tri(\rA)$. Its objects
are triples $(\widetilde\mM, V, \theta)$, where $\widetilde\mM$
is a Cohen-Macaulay $\rR$--module, $V$ is a Noetherian
$\rQ(\bar\rA)$--module and
$\theta:  \rQ(\bar\rR) \otimes_{\rQ(\bar\rA)} V \to
\rQ(\bar\rR) \otimes_\rR \widetilde\mM$ is an epimorphism
of $\rQ(\bar\rR)$--modules such that the induced morphism
of $\rQ(\bar\rA)$--modules
$$V \to \rQ(\bar\rR) \otimes_{\rQ(\bar\rA)} V
\stackrel{\theta}\lar \rQ(\bar\rR) \otimes_\rR \widetilde\mM $$
is an monomorphism.
A morphism between two triples $(\widetilde\mM, V, \theta)$
and $(\widetilde\mM', V', \theta')$ is given by a pair $(F, f)$, where
$F: \widetilde\mM \to \widetilde\mM'$ is a morphism of $\rR$--modules and
$f: V \to V'$ is a morphism of $\rQ(\bar\rA)$--modules such that
the following diagram
$$
\xymatrix
{
\rQ(\bar\rR) \otimes_{\rQ(\bar\rA)} V \ar[rr] \ar[d]_{1 \otimes f} & &
\rQ(\bar\rR) \otimes_\rR \widetilde\mM \ar[d]^{1 \otimes F}\\
\rQ(\bar\rR) \otimes_{\rQ(\bar\rA)} V' \ar[rr] & &
\rQ(\bar\rR) \otimes_\rR \widetilde\mM'
}
$$
is commutative.
\end{definition}

\noindent
The definition is motivated  by the following
theorem.

\medskip

\begin{theorem}[Burban-Drozd]\label{T:BurbanDrozd} Let $\kk$ be an algebraically
closed field and $(\rA, \idm)$ a reduced analytic two-dimensional
Cohen-Macaulay ring which is a non-isolated singularity.
Then in the notations of Definition \ref{D:triples} we have:
the functor
$
\mathbb{F}: \CM(\rA) \to \Tri(\rA)
$
mapping a Cohen-Macaulay module $\mM$ to the triple
$\bigl(\rR\boxtimes_\rA \mM, \rQ(\bar\rA) \otimes_\rA M,
\theta_\mM\bigr)$ is an equivalence of categories.

Moreover, the full
subcategory $\CM^{\mathsf{lf}}(\rA)$ consisting of Cohen-Macaulay
modules which are locally free on the punctured spectrum of $\rA$, is equivalent
to the full subcategory $\Tri^{\mathsf{lf}}(\rA)$ consisting
of those triples $(\widetilde\mM, V, \theta)$ for which the morphism
$\theta$ is an isomorphism.
\end{theorem}

\noindent
The details
of the proof of this theorem will appear in a forthcoming
joint paper of both authors \cite{BurbanDrozd}.
Next, we need
an explicit description of a functor
$\mathbb{G}: \Tri(\rA) \to \CM(\rA)$ (or, at least its description
on objects), which is
 quasi-inverse to $\mathbb{F}$. The construction is as follows.
Let $\mT = (\widetilde\mM, V, \theta)$ be an object
 of $\Tri(\rA)$. Then one can find an $\bar\rA$--module
 $\mX$, a morphism of $\bar\rA$--modules
 $\phi: \mX \to \bar\rR \otimes_\rR \widetilde\mM$
 and an isomorphism  $\psi: \rQ(\bar\rR) \otimes_{\bar\rA} \mX
 \to \rQ(\bar\rR) \otimes_{\rQ(\bar\rA)} V$ such that
 the following diagram
 $$
 \xymatrix
 {
 \rQ(\bar\rR) \otimes_{\bar\rA} \mX \ar[rr]^{1 \otimes \phi}
 \ar[dr]_\psi
 & & \rQ(\bar\rR) \otimes_{\rR} \widetilde\mM \\
  & \rQ(\bar\rR) \otimes_{\rQ(\bar\rA)} V \ar[ur]_\theta
 }
 $$
 is commutative.
 Consider the following commutative diagram with exact rows
 in the category of $\rA$--modules:
 $$
 \xymatrix
 {
 0 \ar[r] & I\widetilde\mM \ar[r] \ar[d]_{=}
 & \mM \ar[r] \ar[d] & \mX \ar[r] \ar[d]^{\phi} &
 0 \\
0 \ar[r] & I\widetilde\mM \ar[r] & \widetilde\mM \ar[r] &
\bar{\rR}\otimes_\rR \widetilde\mM \ar[r] &
 0.
}
 $$
 Then $\mathbb{G}(\mT) \cong M^\dagger$.

\medskip
The aim of this section is to apply Theorem \ref{T:BurbanDrozd}
for a classification of indecomposable Cohen-Macaulay modules
over non-isolated singularities $\kk\llbracket x,y,z\rrbracket /xy$ and
$\kk\llbracket x,y,z\rrbracket /(x^2y - z^2)$.

\medskip
\begin{theorem}\label{T:Ainfty}
Let $\rA = \kk\llbracket x,y,z\rrbracket /xy$. Then the indecomposable non-free
Cohen-Macaulay
$\rA$--modules are described by the following matrix factorizations:
\begin{itemize}
\item $\mM(x,y)$ and $\mM(y,x)$. Let
$\rR = \rR_1 \times \rR_2 = \kk\llbracket x,z_1\rrbracket  \times \kk\llbracket y, z_2\rrbracket $ be  the
normalization of $\rA$ then $\mM(x,y) \cong \rR_2$ and
$\mM(y,x) \cong \rR_1$.
\item $\mM(\varphi_n, \psi_n)$ and $\mM(\psi_n, \varphi_n)$, where
$$
\varphi_n =
\left(
\begin{array}{cc}
y  & z^n \\
0  & x
\end{array}
\right)
\quad
\mbox{and}
\quad
\psi_n =
\left(
\begin{array}{cc}
x & -z^n \\
0 & y
\end{array}
\right), \quad n \ge 1.
$$
\end{itemize}
\end{theorem}

\noindent
\emph{Proof}. Let
$\pi: \kk\llbracket x,y,z\rrbracket /xy \to \kk\llbracket x,z_1\rrbracket  \times \kk\llbracket y,z_2\rrbracket $ be the
normalization map, where $\pi(z) = z_1 + z_2$. Then,
 keeping the notations
of  Theorem \ref{T:BurbanDrozd}, we have:
\begin{itemize}
\item $I = \Ann_\rA(\rR/\rA) = (x, y)$;
\item $\bar\rA = \rA/I = \kk\llbracket z\rrbracket $, $\bar{\rR} =
\kk\llbracket z_1\rrbracket  \times \kk\llbracket z_2\rrbracket $.
\end{itemize}
Since the categories $\CM(\rA)$ and $\Tri(\rA)$ are equivalent, it suffices
to describe indecomposable objects of the category of triples.

Let $\mT = (\widetilde\mM, V, \theta)$ be an object of
$\Tri(\rA)$. Since the semi-local ring $\rR$ is regular,
by Corollary \ref{C:CMreg} any
Cohen-Macaulay $\rR$-module $\widetilde\mM$
has the form $\rR_1^p \oplus \rR_2^q$. Moreover, since
$\rQ(\bar\rA) = \kk((z))$ is a field, we have: $V = \kk((z))^n$.
Note that $\bar{\rR} \otimes_\rR \widetilde\mM
= \kk((z_1))^p  \oplus \kk((z_2))^q$.
Therefore, the gluing  morphism
$\theta: \rQ(\bar\rR) \otimes_{\rQ(\bar\rA)} V \to
\rQ(\bar\rR) \otimes_{\rR} \widetilde\mM$
is given by a pair of matrices $\bigl(\theta_1(z_1), \theta_2(z_2)\bigr) \in
\Mat_{p \times n}\bigl(\kk((z_1))\bigr) \times
\Mat_{q \times n}\bigl(\kk((z_2))\bigr)$. Additional assumptions
on  $\theta$ imply  that
\begin{itemize}
\item
both matrices $\theta_1(z_1)$
and $\theta_1(z_2)$ have full row rank;
\item the matrix
$\displaystyle 
\left(
\begin{array}{c}
\theta_1(z) \\
\hline
\theta_2(z))
\end{array}
\right) \in \Mat_{p + q, n}\bigl(\kk((z)))\bigr)
$
has full
column rank.
\end{itemize}

\noindent
Note that $\Aut_\rR(\rR_1^p) = \GL(p, \rR_1)$ and
$\Aut_\rR(\rR_2^q) = \GL(q, \rR_2)$. Hence,
$$
\Bigl(\rR_1^p \oplus \rR_2^q, \kk((z))^n,
\bigl(\theta_1(z_1), \theta_2(z_2)\bigr)\Bigr)
\cong
 \Bigl(\rR_1^p \oplus \rR_2^q, \kk((z))^n,
 \bigl(\theta'_1(z_1), \theta'_2(z_2)\bigr)\Bigr)$$
 in the category $\Tri(\rA)$ if and only if there exists an element
 $$\bigl(F_1(z_1), F_2(z_2), f(z)\bigr) \in \Bigl(\bigl(\GL(p, \kk\llbracket z_1\rrbracket \bigr)
\times \GL\bigl(q, \kk\llbracket z_2\rrbracket \bigr) \times
 \GL\bigl(n, \kk((z))\bigr)\Bigr)$$ such that
 $$
 \theta'_1(z_1) = F_1^{-1}(z_1) \theta_1(z_1) f(z_1), \quad
 \theta'_2(z_2) = F_2^{-1}(z_2) \theta_2(z_2) f(z_2).
 $$
 Observe that the obtained matrix problem is almost equivalent to
 the problem of classification of indecomposable representations
 of the quiver $\xymatrix{\bullet & \ar[l] \bullet \ar[r] & \bullet}$
 over the field $\kk((z))$ and the following lemma is true.

 \begin{lemma}
 The indecomposable objects of the category $\Tri(\rA)$ are the following:
 \begin{enumerate}
 \item $\Bigl(\rR_1, \kk((z)), \bigl((1), (\emptyset)\bigr)\Bigr)$ and
 $\Bigl(\rR_2, \kk((z)), \bigl((\emptyset), (1)\bigr)\Bigr)$;
 \item $\Bigl(\rR, \kk((z)), \bigl((1), (1)\bigr)\Bigr)$;
 \item $\Bigl(\rR, \kk((z)), \bigl((1), (z_2^n)\bigr)\Bigr)$ and
 $\bigl(\rR, \kk((z)), \bigl((z_1^n), (1)\bigr)\Bigr)$, $n \ge 1$.
 \end{enumerate}
 \end{lemma}

 Note that by Theorem \ref{T:BurbanDrozd} the triples
 of type (2) and (3) correspond to Cohen-Macaulay modules
 which are locally free on the punctured spectrum.
 It is not difficult to see that the indecomposable triples of
 type (1) correspond to both components $\rR_1$ and $\rR_2$ of the
 normalization $\rR$ and the indecomposable triple of type
 (2) is exactly the regular module $\rA$. An interesting problem
 is to describe matrix factorizations corresponding to the triples of type
 (3).

 Let $\mT = \Bigl(\rR, \kk((z)), \bigl((z_1^n), (z_2^m)\bigr)\Bigr)$, 
 where
 $n, m \ge 1$. Note that
 $$
 e_1 = \frac{x}{x+y}
 \quad
 \mbox{and}
 \quad
 e_2 = \frac{y}{x+y}
 $$
 are idempotents in $\rQ(\rA)$ and
 $$
 z_1 = \frac{xz}{x+y}
 \quad
 \mbox{and}
 \quad
 z_2 = \frac{yz}{x+y}
 $$
where $e_1, e_2 \in \rR$ are viewed as elements of $\rQ(\rA) = \rQ(\rR)$.
 Let $\mX = \kk\llbracket z\rrbracket $ and $\phi = (z_1^n, z_2^m)$.
 Consider the pull-back diagram
 $$
 \xymatrix
 {
 0 \ar[r] & I\rR \ar[r] \ar[d]_{=}
 & \mM \ar[r] \ar[d] & \mX \ar[r] \ar[d]^{\phi} &
 0 \\
0 \ar[r] & I\rR \ar[r] & \rR \ar[r] &
\bar{\rR}  \ar[r] &
 0
}
 $$
then $\mathbb{G}(\mT) \cong  \mM^\dagger$. Note that
$
\mM^\dagger =
\bigl(I, z_1^n + z_2^m\bigr)_{\rA}^\dagger,
$
where we view $\bigl(I, z_1^n+  z_2^m\bigr)_{\rA}$ as a submodule
of $\rQ(\rA)$. Since the element $x + y$ is a non-zero divisor in
$\rA$, we have:
$$
\bigl(I, z_1^n + z_2^m\bigr)_{\rA} \cong
(x + y)^{n+m}\bigl(I, z_1^n + z_2^m\bigr)_{\rA} \cong
\bigl(x^{n+m+1}, y^{n+m+1}, x^{n+m} z^n + y^{n+m} z^m\bigr)_{\rA}.
$$
Note that  without loss of generality we may assume that
either $n = 0$ or $m = 0$. For  $m = 0, n \ge 1$ we have:
$
\mM_n :=
\bigl(x^{n+1}, y^{n+1}, x^n z^n + y^n\bigr)_\rA =
\bigl(- x^{n+1},   x^n z^n + y^n\bigr)_\rA.
$
The minimal free resolution of $\mM_n$ is
$$
\dots \lar
\rA^2
 \xrightarrow{
 \left(
 \begin{array}{cc}
 y & z^n \\
 0 & x
 \end{array}
\right)
 }
 \rA^2
 \xrightarrow{
 \left(
 \begin{array}{cc}
 x & -z^n \\
 0 & y
 \end{array}
\right)}
 \rA^2
 \xrightarrow{
 \left(
 \begin{array}{cc}
 y & z^n \\
 0 & x
 \end{array}
\right)
 }
 \rA^2
 \lar \mM_n \to 0.
$$
In particular, $\mM_n \cong  \mM_n^\dagger$ is already Cohen-Macaulay
and $\mM_n = \mM(\varphi_n, \psi_n)$.
\qed

\medskip
\begin{remark}
Using Lemma \ref{L:rankofmf} and keeping the notations of Theorem \ref{T:Ainfty}
it is easy to see that for all $n \ge 1$ the
modules
$\mM(\varphi_n, \psi_n)$ are locally free of rank one
on the punctured spectrum.
\end{remark}

\medskip
\begin{remark}
The statement of Theorem \ref{T:Ainfty} remains
 true  in the case of  an   analytical algebra
$\kk\{x,y,z\}/xy$ with respect to an arbitrary valuation
of $\kk$.
\end{remark}

\medskip
\begin{theorem}\label{T:CMoverDinfty}
Let $\rA = \kk\llbracket x,y,z\rrbracket /(x^2 y  - z^2)$ be the coordinate ring of a surface
singularity of type $D_\infty$.
Then the indecomposable non-free Cohen-Macaulay
$\rA$--modules are given  by the following matrix factorizations:
\begin{enumerate}
\item
$\mM(\alpha^+,\alpha^-)$
where
$$
\alpha^+  =
\left(
\begin{array}{cc}
z  & xy \\
x  & z
\end{array}
\right)
\quad
\mbox{and}
\quad
\alpha^-  =
\left(
\begin{array}{cc}
- z  & xy \\
x  &  - z
\end{array}
\right).
$$
Observe  that $\mM(\alpha^+, \alpha^-) \cong \mM(\alpha^-,
\alpha^+) \cong \rR$, where $\rR$ is the normalization of $\rA$.
\hspace{0.2cm}
\item $\mM(\beta^+, \beta^-)$,
 where
$$
\beta^+  =
\left(
\begin{array}{cc}
x^2  & z \\
z  & y
\end{array}
\right)
\quad
\mbox{and}
\quad
\beta^- =
\left(
\begin{array}{cc}
y & -z \\
-z  & x^2
\end{array}
\right).
$$
\hspace{0.2cm}
\item
$\mM(\gamma_m^+, \gamma^{-}_m) \quad (m \ge 1)$,  where
$$
\gamma_m^+ =
\left(
\begin{array}{cccc}
z & xy & 0 & - y^{m+1} \\
x & z & y^m & 0 \\
0 & 0 & z & xy \\
0 & 0 & x & z
\end{array}
\right)
\quad
\mbox{and}
\quad
\gamma_m^- =
\left(
\begin{array}{cccc}
- z & - xy & 0 &  y^{m+1} \\
x & z & y^m & 0 \\
0 & 0 & - z & - xy \\
0 & 0 & x & z
\end{array}
\right).
$$
\hspace{0.2cm}
\item
$\mM(\delta_m^+, \delta_m^-) \quad (m \ge 1)$,  where
$$
\delta_m^+ =
\left(
\begin{array}{cccc}
z & xy & - y^{m} & 0  \\
x & z &  0 & y^{m}  \\
0 & 0 & z & xy \\
0 & 0 & x & z
\end{array}
\right)
\quad
\mbox{and}
\quad
\delta_m^- =
\left(
\begin{array}{cccc}
- z & - xy &  - y^{m} & 0  \\
x &  z &  0 & - y^{m}  \\
0 & 0 & - z & - xy \\
0 & 0 & x & z
\end{array}
\right).
$$
\end{enumerate}
Note that in all four cases we have: $\mM(\phi, \psi) \cong \mM(\psi, \phi)$.
Moreover, the indecomposable modules of types (2), (3) and (4) are
locally free on the punctured spectrum.
\end{theorem}

\noindent
\emph{Proof}. In the notations of Theorem \ref{T:BurbanDrozd}
we have  $$\pi: \rA = \kk\llbracket x,y,z\rrbracket /(x^2 y - z^2)
 \to \kk\llbracket x, t\rrbracket  = :\rR,$$
 where $\pi(x) = x, \pi(y) = t^2$ and $\pi(z) = tx$, is the
 normalization of $\rA$. It is easy to see that
 $$I = \Ann_\rA(\rR/\rA)  =
 (x, z)\rA
 = x\rR
 $$
 is the conductor ideal. Hence, $\bar\rA = \rA/I =
 \kk\llbracket y\rrbracket  = \kk\llbracket t^2\rrbracket $ and
 $\bar\rR = \rR/I = \kk\llbracket t\rrbracket $.

Let $\mT = (\widetilde\mM, V, \theta)$ be an object of the category
of triples
$\Tri(\rA)$. Since the  ring $\rR$ is regular,
by Corollary \ref{C:CMreg} any
Cohen-Macaulay $\rR$-module $\widetilde\mM$
has the form $\rR^n, \, n \ge 1$. Moreover, because
$\rQ(\bar\rA) = \kk((t^2))$ is a field, we have:
$V = \kk((t^2))^m, \, m \ge 1$.
Note that $\bar{\rR} \otimes_\rR \widetilde\mM
= \kk\llbracket t\rrbracket ^n$.
Therefore, the gluing  morphism
$\theta: \rQ(\bar\rR) \otimes_{\rQ(\bar\rA)} V \to
\rQ(\bar\rR) \otimes_{\rR} \widetilde\mM$
is simply
given by a  matrix
 $\theta(t) \in
\Mat_{m \times n}\bigl(\kk((t))\bigr)$. Additional constrains
on morphism $\theta$ imply  that
\begin{itemize}
\item
the matrix $\theta(t)$ has full row rank
\item if $\theta(t)=\theta_0+t\theta_1$, where $\theta_0,\theta_1 \in
\Mat_{m \times n}\bigl(\kk((t^2))\bigr)$, then the matrix 
$\displaystyle 
\left(
\begin{array}{c}
\theta_1 \\
\hline
\theta_2
\end{array}
\right)
$
 has full column rank.
\end{itemize}

\noindent
Note that $\Aut_\rR(\rR^m) = \GL(m, \rR)$. Moreover,
$F \in \Mat_{m\times m}(\rR)$ is invertible if and only if
$F(0)$ is invertible.   Therefore, the existence of an isomorphism
$$
\bigl(\rR^m, \kk((t^2))^n, \theta(t)\bigr)
\cong
\bigl(\rR^m, \kk((t^2))^n, \theta'(t)\bigr)
 $$
 in the category $\Tri(\rA)$ is equivalent to the existence
 of  an element
 $$\bigl(F(t), f(t^2)\bigr) \in
 \Bigl(\GL\bigl(m, \kk\llbracket t\rrbracket \bigr) \times 
 \GL\bigl(n, \kk((t^2))\bigr)\Bigr)$$
 such that
 \begin{equation}\label{E:transruleD}
 \theta' = F^{-1} \theta f.
 \end{equation}

\noindent
 A classification of indecomposable Cohen-Macaulay modules over
 $\rA$ follows  from the following lemma.

\begin{lemma}\label{L:keylemmaD}
Let a matrix $\theta \in \Mat_{m \times n}\bigl(\kk((t))\bigr)$ be of  full
row rank. Then applying the transformation rule (\ref{E:transruleD})
we can decompose $\theta$ to a direct sum of the following matrices:
$$(1), \, \, 
(t), \, \, 
(1 \,\, t)
$$
and
$$
\left(
\begin{array}{cc}
1 & t \\
t^d & 0
\end{array}
\right), \quad d \ge 1.
$$
\end{lemma}

\noindent
Assume Lemma \ref{L:keylemmaD} is proven. Since the isomorphism
classes of indecomposable Cohen-Macaulay $\rA$--modules
stand in bijection with the equivalence classes of matrices
$\theta(t)$ modulo the transformation rule (\ref{E:transruleD}),
the problem reduces to a description of modules corresponding to the
canonical forms listed above.

It is clear that $\theta = 1$ corresponds to the free module
$\rA$. In a similar way,  it is not difficult to show that
$\rR$ corresponds to $\theta = (1\,\, t)$.
Now let  us describe  modules corresponding to the matrices
$$
\theta =
\left(
\begin{array}{cc}
1 & t \\
t^d & 0
\end{array}
\right), \quad d \ge 1.
$$
Denote  $\mX = \kk\llbracket t^2\rrbracket  \oplus \kk\llbracket t^2\rrbracket $, put  $\phi = \theta$ and
 consider the pull-back diagram
 $$
 \xymatrix
 {
 0 \ar[r] & I\rR^2 \ar[r] \ar[d]_{=}
 & \mM \ar[r] \ar[d] & \mX \ar[r] \ar[d]^{\phi} &
 0 \\
0 \ar[r] & I\rR^2 \ar[r] & \rR^2 \ar[r] &
\bar{\rR}^2  \ar[r] &
 0.
}
 $$
Then $\mathbb{G}(\mT) \cong  \mM^\dagger$. Since
$t = \frac{\displaystyle z}{\displaystyle x}$, we have
$$
\begin{array}{r@{\ }c@{\ }l}
\mM & = &
\left(I(\rR^2),
\left(
\begin{array}{c}
1 \\
\frac{\displaystyle z^d}{\displaystyle x^d}
\end{array}
\right),
\left(
\begin{array}{c}
\frac{\displaystyle z}{\displaystyle x} \\
0
\end{array}
\right)
\right)_{\rA}  
\\
& = &
\left(
\left(
\begin{array}{c}
x \\
0
\end{array}
\right),
\left(
\begin{array}{c}
z \\
0
\end{array}
\right),
\left(
\begin{array}{c}
0 \\
x
\end{array}
\right),
\left(
\begin{array}{c}
0 \\
z
\end{array}
\right),
\left(
\begin{array}{c}
1 \\
\frac{\displaystyle z^d}{\displaystyle x^d}
\end{array}
\right),
\left(
\begin{array}{c}
\frac{\displaystyle z}{\displaystyle x} \\
0
\end{array}
\right)
\right)_{\rA},
\end{array}
$$
where $\mM$ is considered as a submodule of $\rQ(\rA)^2$.
Note that$$\left(\begin{array}{c} z \\ 0 \end{array}\right) =
x \left(\begin{array}{c} \frac{\displaystyle z}{\displaystyle x} \\ 0 \end{array}\right),
$$
hence  the second   generator can be omitted.
Since  $\rA$ is an integral domain, we have:
$$
\left(
\begin{array}{cc}
x  & 0 \\
0 & x^d
\end{array}
\right) \mM \cong \mM.
$$
This implies:
$$
\mM \cong
\left(
\left(
\begin{array}{c}
x^2 \\
0
\end{array}
\right),
\left(
\begin{array}{c}
0 \\
x^{d+1}
\end{array}
\right),
\left(
\begin{array}{c}
0 \\
x^d z
\end{array}
\right),
\left(
\begin{array}{c}
x \\
z^d
\end{array}
\right),
\left(
\begin{array}{c}
z \\
0
\end{array}
\right)
\right)_{\rA}.
$$
Next, we show that the first generator can be expressed via the
remaining four. First note that
$$
\left(\begin{array}{c} x^2 \\ 0 \end{array}\right) =
x
\left(\begin{array}{c}
x  \\ z^d \end{array}
\right) -
\left(\begin{array}{c}
0  \\
x z^d
\end{array}
\right).
$$
Now distinguish two cases.
\begin{itemize}
\item Let $d = 2m, m \ge 1$, be even. Since
$x z^{2m} = x^{2m+1} y^{2m} = x^{d+1} y^m$, we have:
$$
\left(\begin{array}{c}
0  \\
x z^{2m}
\end{array}
\right)
=
y^m
\left(\begin{array}{c}
0  \\
x^{2m+1}
\end{array}
\right) =
y^m
\left(\begin{array}{c}
0  \\
x^{d+1}
\end{array}
\right).
$$
\item Similarly, if $d = 2m +1, m \ge 0$, is odd, we have
$xz^d = x z^{2m+1} = x  x^{2m} y^m z = y^m x^d z$ and
$$
\left(\begin{array}{c}
0  \\
x z^d
\end{array}
\right)
=
y^m
\left(\begin{array}{c}
0  \\
x^d z
\end{array}
\right).
$$
\end{itemize}

\noindent
Thus, we obtain:
$$
\mM \cong
\left(
\left(
\begin{array}{c}
0 \\
x^{d+1}
\end{array}
\right),
\left(
\begin{array}{c}
0 \\
x^d z
\end{array}
\right),
\left(
\begin{array}{c}
x \\
z^d
\end{array}
\right),
\left(
\begin{array}{c}
z \\
0
\end{array}
\right)
\right)_{\rA} \subset \rA^2.
$$
It turns out that this module is already Cohen-Macaulay. Indeed,
a straightforward calculation shows that $\syz(\mM) = \rho$, where
$$
\rho =
\left(
\begin{array}{cccc}
z & xy & 0 & -y^{m+1} \\
-x & -z & -y^m & 0 \\
0 & 0 & z & xy \\
0 & 0 & -x & -z
\end{array}
\right)
$$
in the case  $d = 2m, m \ge 1$, is even, and
$$
\rho =
\left(
\begin{array}{cccc}
z & xy & - y^{m+1} & 0 \\
-x & -z & 0  & -y^{m+1} \\
0 & 0 & z & xy \\
0 & 0 & -x & -z
\end{array}
\right) \sim
\left(
\begin{array}{cccc}
z & xy & - y^{m+1} & 0 \\
x & z & 0  & y^{m+1} \\
0 & 0 & z & xy \\
0 & 0 & x & z
\end{array}
\right)
$$
if   $d = 2m + 1, m \ge 0$, is odd.
Moreover, we have  the following equalities
$$
\left(
\begin{array}{cccc}
z & xy & 0 & -y^{m+1} \\
x & z & y^m & 0 \\
0 & 0 & z & xy \\
0 & 0 & x & z
\end{array}
\right)
\left(
\begin{array}{cccc}
z & xy & 0 & -y^{m+1} \\
-x & -z & -y^m & 0 \\
0 & 0 & z & xy \\
0 & 0 & -x & -z
\end{array}
\right) = (x^2 y - z^2)I_4
$$
and
$$
\left(
\begin{array}{cccc}
z & xy & y^{m+1} & 0  \\
-x & -z &  0 & y^{m+1} \\
0 & 0 & z & xy \\
0 & 0 & -x & -z
\end{array}
\right)
\left(
\begin{array}{cccc}
z & xy & -y^{m+1} & 0  \\
x & z & 0 & -y^{m+1} \\
0 & 0 & z & xy \\
0 & 0 & x & z
\end{array}
\right) = (x^2 y - z^2)I_4
$$
which imply that the $\rA$-module $\mM$ has a 2--periodic resolution,
hence is Cohen-Macaulay.

It remains to identify the Cohen-Macaulay $\rA$--module corresponding
to the canonical form $\theta = (t)$. Since $t = \frac{\displaystyle z}{\displaystyle x}$, the
corresponding Cohen-Macaulay module is
$$
\mM = \Bigl(x, z, \frac{z}{x}\Bigr)_\rA \cong
\Bigl(x^2, z\Bigr)_\rA.
$$
It is easy to see that the corresponding matrix factorization
is
$\mM(\beta^+, \beta^-)$,
 where
$$
\beta^+  =
\left(
\begin{array}{cc}
x^2  & z \\
z  & y
\end{array}
\right)
\quad
\mbox{and}
\quad
\beta^- =
\left(
\begin{array}{cc}
y & -z \\
-z  & x^2
\end{array}
\right).
$$
To complete the proof of Theorem \ref{T:CMoverDinfty}
it remains to prove Lemma \ref{L:keylemmaD}.

\medskip
\noindent
\emph{Proof of Lemma \ref{L:keylemmaD}.}
It is a straightforward verification that all listed
canonical forms can not be further split
 and are pairwise non-isomorphic.

We prove using the induction on  the size, that
any matrix $\theta \in \Mat_{m\times n}\bigl(\kk((t))\bigr)$
 of full row rank can
be decomposed into a direct sum of blocks listed in
the formulation of Lemma \ref{L:keylemmaD}.
 The case $m = 1$ is clear.
Using the transformation
rule (\ref{E:transruleD}) the matrix
  $\theta$ can be reduced to the following
form:
$$
\left(
\begin{array}{ccccc}
t^{d_1} I_{s_1} & A_{1,2} & \dots & A_{1,l} & A_{1, l+1} \\
0 & t^{d_2} I_{s_2} & \dots & A_{2,l} & A_{2, l+1} \\
\vdots & \vdots & \ddots & \vdots & \vdots \\
0 & 0 & \dots & t^{d_l} I_{s_l} & A_{l,l+1}
\end{array}
\right),
$$
where $d_1 < d_2 < \dots < d_l$ and
$A_{i,j} \in \Mat_{s_i \times s_j}(t^{d_i+ 1}\bar\rR)$ for
$i \le j \le l$ and
$\bar\rR = \kk\llbracket t\rrbracket $.
\begin{itemize}
\item If all matrices $A_{1,j} = 0$ for $2 \le j \le l+1$,
either $(1)$ or $(t)$ is a direct summand of $\theta$ and we
can proceed by induction.
\item Assume $A_{1,j} \ne 0$ and $A_{1,p} = 0$ for $2 \le p < j\le l$.
It is easy to see that $A_{1,j}$ can be transformed to the form
$$
A'_{1,j} =
t^{d_1 + 1}
\left(
\begin{array}{cc}
I_{s_1'} & 0 \\
0 & 0
\end{array}
\right),
$$
where $s_1'$ can be also $0$. If $s_1' = 0$, we return to the previous
step, otherwise we can split the matrix
$$
\left(
\begin{array}{cc}
t^{d_1} & t^{d_1 +1} \\
0       & t^{d_j}
\end{array}
\right)
$$
as a direct summand. If $d_j = d_2 = d_1 +1$, this matrix splits
into a direct sum of two $1 \times 1$--matrices, otherwise it is equivalent
to a matrix of the form
$$
\left(
\begin{array}{cc}
1 & t \\
t^d & 0
\end{array}
\right)
$$
for a proper $d \ge 1$.
\item Finally, if $A_{1, j} = 0$ for $2 \le j \le l$ and
$A_{1, l+1} \ne 0$, then we can split either
$(1)$, or $(t)$ or $(t^{d_1} \, \, t^{d_1+ 1}) \sim (1 \,\, t)$ as a
direct summand of $\theta$.
 This concludes the induction step and proves
Lemma \ref{L:keylemmaD}, hence Theorem \ref{T:CMoverDinfty}.
\end{itemize}
\qed

\medskip
\begin{remark}
In the notations of Theorem \ref{T:CMoverDinfty} the rank
of the Cohen-Macaulay modules $\mM(\alpha^+, \alpha^-)$ and
$\mM(\beta^+, \beta^-)$ is one, whereas
$\mM(\gamma_m^+, \gamma_m^-)$ and $\mM(\delta_m^+, \delta_m^-)$
have ranks two.
\end{remark}

\medskip
\begin{remark}
Note that for $\rA = \kk\llbracket x,y,z\rrbracket /(x^2 y - z^2)$
the stable category $\underline{\CM}(\rA)$ has the following
interesting property:
it is a triangulated category with a shift functor $\mathbb{T}$ such
that $\mathbb{T}(\mM) \cong \mM$ for any object $\mM$ of
$\underline{\CM}(\rA)$.  
Moreover, the stable category
$\underline{\CM}^{\mathsf{lf}}(\rA)$
of Cohen-Macaulay
$\rA$--modules, which are locally free on the punctured spectrum, is
a $\Hom$--finite triangulated category having the same property.
\end{remark}

\medskip

\section{Geometric McKay correspondence for simple surface
singularities}

The main goal of this section is to give a geometric description
of indecomposable Cohen-Macaulay modules on simple surface
singularities. Throughout this section the base field
$\kk$ is equal to $\mathbb{C}$ and we work either in the category
of complex analytic spaces or in the category of algebraic schemes.
Let us first recall some basic results
on resolutions of singularities.

Let $(\sX, o)$ be the   germ of a normal surface singularity
and  $\pi: (\widetilde{\sX}, \sE)
\to (\sX, o)$ a resolution of singularities,  i.e.

\begin{itemize}
\item $\widetilde{\sX}$ is a  germ of a smooth surface;
\item $\pi$ is a proper morphism of germs of complex-analytic spaces;
\item
$\pi: \widetilde\sX\setminus\sE \to \sX\setminus o$ is an isomorphism
and $\widetilde\sX\setminus\sE$ is dense in $\widetilde\sX$. In particular,
$\pi$ is birationally an isomorphism.
\end{itemize}

In this case the \emph{exceptional fiber}
$\sE = \pi^{-1}(o)$  is a  complex projective curve
(possibly singular). Remind that $\sE$ is always connected (it follows
from Zariski's Main Theorem, see \cite[Corollary III.11.4]{Hartshorne}).

\begin{definition}
A resolution of singularities
$\pi: \widetilde\sX \to \sX$ is called
\emph{minimal} if $\widetilde{\sX}$ does not contain contractible
curves, i.e. smooth rational
projective curves with self-intersection $-1$.
\end{definition}

\medskip
\begin{remark} A minimal resolution of singularities has the following
universal property. If   $\pi': (\widetilde\sX', \sE') \to (\sX, o)$ is any other  resolution
of singularities,
then there exists a unique morphism  of  germs of
complex-analytic spaces
$f: (\widetilde\sX', \sE') \to
(\widetilde\sX, \sE)$ such that the diagram
$$
\xymatrix
{
\widetilde\sX' \ar[rr]^f \ar[rd]_{\pi'}
& & \widetilde\sX \ar[ld]^\pi \\
 & \sX &
}
$$
is commutative, see \cite[Section V.5]{Hartshorne}, \cite{Grauert} and
\cite{Laufer}.
\end{remark}

\medskip
\noindent
The main topological invariant of a normal surface singularity
$(\sX, o)$ is the intersection matrix of its minimal resolution.

\medskip
\begin{definition}
Let $\pi: (\widetilde\sX, \sE) \to (\sX, o)$ be a minimal
resolution of a normal surface singularity $(\sX, o)$ and
$\sE = \cup_{i= 1}^n \sE_i$ the  decomposition
of the exceptional divisor
$\sE$ into a union of  irreducible components.  Then
$M(\sX) = (m_{i,j}) = (E_i \cdot E_j)_{i,j}
\in \Mat_{n \times n}(\mathbb{Z})$ is called
 the
\emph{intersection
matrix}.
\end{definition}

\medskip
\begin{proposition}[Mumford]
Let $(\sX,  o)$ be a normal surface singularity, then
its  intersection matrix $M(\sX)$ is non-degenerate and
negatively definite.  
\end{proposition}

\noindent
For a proof of this result, see for example \cite[Theorem 4.4]{Laufer}.

\medskip
The key technique to resolve normal surface singularities is provided by the following construction.
Assume for simplicity of notation that we work in the category of algebraic schemes over an algebraically
closed  field $\kk$.
Consider the scheme $\widetilde{\mathbb{A}}^3$ defined as follows:
$$
\mathbb{A}^3 \times \mathbb{P}^2 \supseteq \widetilde{\mathbb{A}}^3 :=
\Bigl\{\bigl((x_1, x_2, x_3), (y_1: y_2: y_3)\bigr) \, \big| \, x_i y_j = x_j y_i, \quad 1 \le i, j \le 3\Bigr\},
$$
where $\mathbb{A}^3$ is the three-dimensional affine space over $\kk$ and $\mathbb{P}^2$ is
the two-dimensional
projective space over $\kk$. Note that $\widetilde{\mathbb{A}}^3$ is an algebraic scheme
of Krull dimension  three (although not affine) and we have a morphism of schemes
$\pi: \widetilde{\mathbb{A}}^3 \to \mathbb{A}^3$ defined  as the composition of the inclusion
$\widetilde{\mathbb{A}}^3 \hookrightarrow  \mathbb{A}^3 \times \mathbb{P}^2$ and the canonical projection
$\mathbb{A}^3 \times \mathbb{P}^2 \to \mathbb{A}^3$. Moreover, the following properties hold:

\medskip

\noindent
1. The morphism $\pi$ is projective (hence proper).

\noindent
2.  Let $o := (0, 0, 0) \in \mathbb{A}^3$. Then
$\pi^{-1}(o) =
\Bigl\{\bigl((0, 0, 0), (y_1: y_2: y_3)\bigr)\Bigr\} \cong \mathbb{P}^2
$
and for any other point $x \in \mathbb{A}^3\setminus\{o\}$ the preimage $\pi^{-1}(x)$ is a single point.
Hence, $\pi: \widetilde{\mathbb{A}}^3 \setminus \pi^{-1}(o) \lar  \mathbb{A}^3\setminus\{o\}$
is an isomorphism and $\pi$ itself is a birational isomorphism.

\noindent
3. The projective space $\mathbb{P}^2$ has three affine charts
$$U_i = \bigl\{(y_1:y_2:y_3) \, \big| \, y_i \ne 0\bigr\}, 1 \le i \le 3.$$
 Moreover, we have
isomorphisms
$\phi_i: \widetilde{\mathbb{A}}^3_i := \widetilde{\mathbb{A}}^3\cap \bigl(\mathbb{A}^3 \times U_i \bigr) \stackrel{\cong}\lar
\mathbb{A}^3$, $1 \le i \le 3$.
 For example, take $i = 1$ and denote 
 $\tilde{y}_2 =  \frac{y_2}{y_1}$ and
$\tilde{y}_3 = \frac{y_3}{y_1}$. Then
$$
\widetilde{\mathbb{A}}^3_1 = \Bigl\{\bigl((x_1, x_2, x_3), (1: \tilde{y}_2: \tilde{y}_3)\bigr) \, \mid \,
x_2 = x_1 \tilde{y}_2, \quad x_3 = x_1 \tilde{y}_3 \Bigr\}
$$
and the isomorphism $\phi_1: \widetilde{\mathbb{A}}^3_1 \to \mathbb{A}^3$ is given by the
formula  $$\bigl((x_1, x_2, x_3), (1: \tilde{y}_2: \tilde{y}_3)\bigr) \mapsto (x_1, \tilde{y}_2,  \tilde{y}_3).
$$
In particular, the morphism
$$
\pi_1: \Spec\bigl(\kk[x_1,\tilde{y}_2,  \tilde{y}_3]\bigr) = \mathbb{A}^3
\stackrel{\phi^{-1}_1}\lar \widetilde{\mathbb{A}}^3_1 \hookrightarrow
\widetilde{\mathbb{A}}^3  \stackrel{\pi}\lar  \mathbb{A}^3 =
\Spec\bigl(\kk[x_1, x_2, x_3]\bigr)
$$
is given by the formulae: $x_1 \mapsto x_1, x_2 \mapsto x_1 \tilde{y}_2$ and
$x_3 \mapsto x_1 \tilde{y}_3$.

\medskip
Let $f \in \kk[x,y,z]$ be a polynomial such that $\sX := V(f) \subseteq \mathbb{A}^3$
is a normal surface and
$o \in \sX$  is a singular point. Let
$\widetilde{\sX} := \overline{\pi^{-1}(\sX \setminus\{o\}} \subseteq \widetilde{\mathbb{A}}^3$.
From  the  commutative diagram
$$
\xymatrix
{\pi^{-1}(\sX \setminus\{o\}) \ar@{^{(}->}[rr]  \ar@{_{(}->}[d] & &
\widetilde{\mathbb{A}}^3 \setminus \pi^{-1}(o) \ar@{^{(}->}[d] \\
\sX \setminus\{o\} \ar@{^{(}->}[rr] & & \mathbb{A}^3 \setminus\{o\}
}
$$
we obtain an induced  morphism
 $\pi_\sX:
\widetilde{\sX} \to \sX$, which  is proper and birational. Moreover,
the exceptional fiber $\sE = \pi_\sX^{-1}(o) = \widetilde{\sX} \cap \pi^{-1}(o) \subseteq
\mathbb{P}^2$ is a closed subscheme of a projective plane,
hence it is a projective curve.
The constructed morphism $\pi_\sX: \widetilde\sX \to \sX$ is called
the \emph{blowing-up}
of the surface $\sX$ at the  singular point  $o$.

\medskip
\noindent
We illustrate the technique of resolutions of surface singularities on two examples.

\begin{example}\label{E:resol1}
Consider a surface singularity of type
$A_1$ given by the equation $x^2 + y^2 + z^2 = 0$.
Taking the chart $u \ne 0$ of the scheme
  $\widetilde{\mathbb{A}}^3 = \Bigl\{\bigl((x,y,z), (u:v:w)\bigr) \, \mid \,
xv = yu, xw = zu, yw =zv \Bigr\}$, we get the following
morphisms of affine schemes:
$$
\mathbb{A}^3 \stackrel{\pi}\longleftarrow \widetilde{\mathbb{A}}^3_1
:= \Bigl\{\bigl((x,y,z), (1:v:w)\bigr) \, \mid \,
xv = y,  yw =z \Bigr\}  \stackrel{\phi_1}\lar \mathbb{A}^3,
$$
where $\pi\Bigl(\bigl((x,y,z), (1:v:w)\bigr)\Bigr) =(x,y,z)$ and
$\phi_1\Bigl(\bigl((x,y,z), (1:v:w)\bigr)\Bigr) = (x, xv, xw)$. Then  the morphism
$\pi_1^*  = (\pi \circ \phi_1^{-1})^*:  \kk[x, y, z] \to \kk[x, v, w]$
 maps
$x^2 + y^2 + z^2$ to $x^2 + x^2 v^2 + x^2 w^2$. It is easy to see that
$$\phi_1\bigl(\widetilde\sX \cap \widetilde{\mathbb{A}}^3_1\bigr) =
\phi_1\bigl(\overline{\pi^{-1}(\sX\setminus\{o\})} \cap \widetilde{\mathbb{A}}^3_1\bigr) \subseteq
\mathbb{A}^3 = \Spec\bigl(\kk[x,v,w]\bigr)$$
 is given by the equation $1 + v^2 + w^2 = 0$, which is a smooth surface.

Since the fiber $\pi^{-1}(o)$ of the morphism
$\pi: \widetilde{\mathbb{A}}^3 \to \mathbb{A}^3$
in the local chart $\phi_1(\widetilde{\mathbb{A}}^3)$ in the chart 
$\widetilde{\mathbb{A}}^3_1$  is just $\mathbb{A}^2 =
\{(0, v, w)\} = V(x)$, the exceptional curve  $\sE = \widetilde{\sX} \cap \pi^{-1}(o)$
is $V(x, 1+ v^2 + w^2) \cong \mathbb{A}^1$.

The description $\pi$ in other charts $v \ne 0$ and $w \ne 0$ is completely symmetric.
Hence, the blowing-up $\pi: \widetilde{\sX} \to \sX$ is  already a  resolution of singularities
and the exceptional fiber $\sE = \pi^{-1}(o)$ is isomorphic to 
$\mathbb{P}^1$.

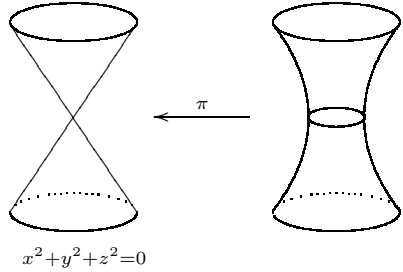
\begin{center}
\hspace{0.7cm}
\xy /r0.2pc/:
{\POS(10,30)\ellipse(10,3)_,=:a(360){-}};
{\POS(10,0)\ellipse(10,3){.}\ellipse(10,3)__a{-}};
{\ar@{-}(0,29.7);(20,0.1)}{\ar@{-}(20,29.7);(0,0.1)}
\POS(0,0);
{(50,30)\ellipse(10,3)_,=:a(360){-}};
\POS(50,0)="a";
{"a"\ellipse(10,3)=:a(200){-}\ellipse(10,3)a{.}}
\POS(50,15)="b";
{"b"\ellipse(4.5,1.5),a(200){-}\ellipse(4.5,1.5)a{.}}
{\ar@{-}@/^13.5pt/(40,29.7);(40,0.2)}{\ar@{-}@/_13.5pt/(60,29.7);(60,0.2)}
{\ar_{\pi}(35,15);(20,15)}
\POS(5,-10)
{\ar@{}^{x^2+y^2+z^2=0}(-10,-10),(10,-10)} 
\endxy
\vspace{0.5cm}
\end{center}

The next problem is to compute the self-intersection number of the exceptional divisor
$\sE$. To do  this,  we use the following trick.
Since the constructed morphism $\pi: \widetilde{\sX} \to \sX$ is a birational isomorphism, it
induces an isomorphism of the fields or rational functions
$\pi^*: \kk(\sX) \to \kk(\widetilde\sX)$. Moreover, for any function
$f \in \kk[\sX]$ we have:
$
\mathsf{div}\bigl(\pi^{*}(f)\bigr)\cdot \sE = 0.
$
In our particular case, take $f$ to be the class of $y$ in the ring
 $\kk[x,y,z]/(x^2 + x^2 + z^2)$.
Consider  the chart  $u \ne 0$ in $\widetilde{A}^3$ and the following commutative diagram
$$
\xymatrix
{
 & \widetilde{\sX} \cap \widetilde{\mathbb{A}}^3_1  \ar[r]  \ar@{_{(}->}[d] &
\sX  \ar@{^{(}->}[d] \\
\mathbb{A}^3 & \widetilde{\mathbb{A}}^3_i \ar[r]^{\pi} \ar[l]_{\phi_1} & \mathbb{A}^3
}
$$
where $\pi_1 = \pi \circ \phi_1^{-1}: \mathbb{A}^3 \to \mathbb{A}^3$ is given
by $\pi_1(x,v,w) = (x, xv, xw)$. It is easy to see 
that  $V\bigl(\pi^*(f)\bigr) \subseteq \widetilde\sX$  
is given by the equation $xv = 0$.
Moreover, we have:
$$\mathsf{div}\bigl(\pi^*(f)\bigr) = \sE + \sC_1 + \sC_2,$$
 where
$\sE = V(x)$ is the exceptional fiber   and two other components are $C_1 = V(v, w - \sqrt{-1})$
and  $C_2 = V(v, w + \sqrt{-1})$.  Since $\sE$ intersect $\sC_1$ and $\sC_2$
transversally exactly at one point, we get:
$$
0 = \mathsf{div}\bigl(\pi^{*}(f)\bigr) \cdot \sE = (\sE + \sC_1 + \sC_2) \cdot \sE =
\sE^2 + \sC_1 \cdot \sE + \sC_2 \cdot \sE =  \sE^2 + 2,
$$
hence $\sE^2 = -2$.
\end{example}

\medskip
In Example \ref{E:resol1} a minimal resolution of an $A_1$-singularity was achieved by a single blow-up.
Of coarse, it does not reflect the situation one has in general. The following example shows various tricks
and pitfalls one can meet by resolving normal surface singularities.

\medskip

\begin{example}\label{E:resol2}
Let $\sX \subseteq \mathbb{A}^3$ be the surface given by  the equation $x^2 + y^3 + z^7 + y^2 z^2 = 0$.
At the point $o = (0,0,0)$ is has the  so-called $T_{2,3,7}$--singularity. Let us show how this
normal singularity can be resolved.

\medskip
\noindent
\emph{Step 1}. Recall that $\widetilde{\mathbb{A}}^3 =
\Bigl\{\bigl((x,y,z), (u:v:w)\bigr) \, \big| \,
xv = uy,  yw =uz \Bigr\}$. Let
$q: \widetilde{\mathbb{A}}^3 \to \mathbb{A}^3$ be the projection morphism
$q\bigl((x,y,z), (u:v:w)\bigr) = (x, y, z)$ and
 $\sY = \overline{\pi^{-1}(X\setminus\{o\})} \stackrel{q}\lar \sX$
 be the blowing-up of $\sX$
at the point $o$.
Consider the chart $w \ne 0$ and denote $\tilde{u} =
\frac{\displaystyle u}{\displaystyle w}$, $\tilde{v} = 
\frac{\displaystyle v}{\displaystyle w}$.
In these coordinates,  the surface $\sY$ is given by the equation
$
\tilde{u}^2 + z \tilde{v}^3 + z^5 + z^2 \tilde{v}^2 = 0.
$
At the point $(0,0,0)$ it has an isolated singularity.
The exceptional divisor
$\sE_1 = q^{-1}(o)$ is given in this chart by the equations $ u = z = 0$.
It is not difficult to check that in two other charts $u \ne 0$ and $v \ne 0$ the surface $\sY$ is smooth
and that $\sE_1 \cong \mathbb{P}^1$.

\medskip

\noindent
\emph{Step 2}. For the sake of simplicity of the notation, denote again $x = \tilde{u}$ and  $y = \tilde{v}$.
The surface singularity $\sY$ is defined by the equation
$$
x^2 + zy^3 + z^2 y^2 + z^5 = 0
$$
and $o = (0, 0, 0)$ is its unique singular point. We also have to keep track of the equation
of the exceptional fiber $\sE_1 = V(x, z)$. Consider again the blowing-up
 $q: \widetilde{\sY} \to \sY$  of the surface $\sY$
at the point $o$ obtained from the morphism
$\widetilde{\mathbb{A}}^3  = \Bigl\{\bigl((x,y,z), (u:v:w)\bigr) \, \big|  \,
xv = uy,  yw =uz \Bigr\}  \to \mathbb{A}^3$.
This time, however, we have to look at two charts.

\medskip
\noindent
\emph{Case 1}. Consider first the chart $w  \ne 0$ and denote 
$u_1 = \frac{\displaystyle u}{\displaystyle w}$, 
$v_1 = \frac{\displaystyle v}{\displaystyle w}$. 
In these coordinates, the surface
$\widetilde{\sY}$ is given  by the equation
$$
u_1^2 + z^2 v_1^3 + z^4 v_1^2 + z^5 = 0.
$$
Note that $\widetilde{\sY}$  is not normal! Indeed, in this chart the singular locus of
$\widetilde{Y}$ is given by the equations $z = u_1 = 0$. Moreover, the singular locus coincides
with the exceptional fiber $\sE_2 =  p^{-1}(o)$! Note also, that
$q^{-1}(\sE_1)$ does not belong to this chart.

\medskip
\noindent
\emph{Case 2}.  Consider another chart  $v \ne 0$
and denote $u_2 = \frac{\displaystyle u}{\displaystyle v}$, 
$w_2 = \frac{\displaystyle w}{\displaystyle v}$.
In these coordinates, the surface $\widetilde{\sY}$  is given by the equation
$$
u_2^2 + y^2 w_2 + y^2 w_2^2 + y^3 w_2^5 = 0.
$$
In this chart, the exceptional fiber $\sE_2 = V(u_2, y)$ coincides with
the singular locus of $\widetilde{\sY}$  and $q^{-1}(\sE_1) = V(u_2, w_2)$.

\medskip
\noindent
\emph{Case 3}. In the third chart $v \ne 0$,  the surface
 $\widetilde{\sY}$ is given by the following  equation:
$$
1 + x^2 v_3^3 w_3 + x^2 v_3^2 w_3^2 + x^3 w_3^5 = 0.
$$
It is easy to see that  in this chart $\widetilde{\sY}$ is smooth.

\medskip
\noindent
Summing everything up, we get that $\widetilde{\sY} \subset \widetilde{\mathbb{A}}^3
\subseteq \mathbb{A}^3 \times \mathbb{P}^2$ is a surface with one-dimensional singular locus,
$$
\sE_1 = \Bigl\{\bigl((0,y,0), (0:1:0)\bigr) \, \big| \,   y \in \mathbb{A}^1 \Bigr\}
$$
is the strict transform of the first exceptional divisor and
$$
\sE_2 = \Bigl\{\bigl((0,0,0), (0:v:w)\bigr) \, \big| \, 
(v:w) \in \mathbb{P}^2 \Bigr\}
$$
is the the second exceptional divisor, which coincides with  the singular locus
of $\widetilde\sY$.

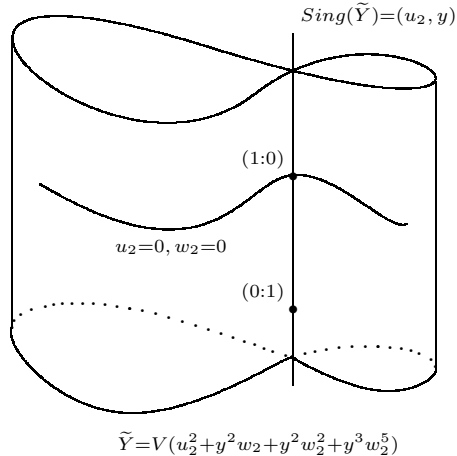
\begin{center}
\hspace{1cm}
\xy 0;/r3pc/:="p",
"p";p+(.5,-.5)="c",
,"p"+(.825,-.25)="x","c"
,{\xycompile{\ellipse~c(5)`"x"{-}}}
,{\xycompile{\ellipse~c(-5)`"x"{-}}}
,"c"+(.4,.235)="t",
,"t"+(1.5,-0.1)="r",
,"p"+(-2.05,+0.05)="l",
,"p"+(0,-3)="p1",
,"x"+(0,-3)="x1",
,"c"+(0,-3)="c1",
,"t"+(0,-3)="t1",
,"r"+(0,-3)="r1",
,"l"+(0,-3)="l1",
,"t"+(0,-1.5)="t2", "t2"+(0,0.4)*\dir{*}*+!DR{\scriptstyle{(1:0)}},
                    "t2"+(0,-1)*\dir{*}*+!DR{\scriptstyle{(0:1)}},
,"r"+(0,-1.5)="r2",
,"l"+(0,-1.5)="l2",
,"c"+(0,-1.5)="c2",
%
,"l2"+(0.3,0);"r2"+(-0.3,0)
**\crv{"c2"+(-0.5,-0.5)&"p"+(0.7,-0.5) &"r2"+(-0.5,-0.1)},
"l2"+(1,-0.5)*+!UL{\scriptstyle{u_2 = 0, \, w_2 =0}}
%
,"r";"r1"**\dir{-},
,"l";"l1"**\dir{-},
,"t"+(0,0.4)*+!DL{\scriptstyle{Sing(\widetilde{Y})=(u_2, \, y)}},;"t1"+(0,-0.3)**\dir{-},
"t1";"l1" **\crv{"t1"+(-0.1,0)&"c1"+(-1.5,-1)&"l1"+(-0.1,-0.2) },
**\crv{~*=<4pt>{.}"t1"+(-0.1,0)&"p1"+(-1,0.5)&"l1"+(-0.0,0.2)},
"t1";"r1" **\crv{"t1"+(-0.03,0)&"r1"+(-0.5,-0.5)&"r1"+(+0.05,-0.1) },
**\crv{~*=<4pt>{.}"t1"+(0.1,0.05)&"r1"+(-0.5,0.3)&"r1"+(-0.,0.1) },
"t1"+(1.2,-0.7)*+!UR{\scriptstyle{\widetilde{Y} =
V(u_2^2 + y^2 w_2 + y^2 w_2^2 + y^3 w_2^5)}}
\endxy
\end{center}

\medskip
\noindent
\emph{Step 3}. We resolve singularities of $\widetilde{\sY}$ using the \emph{normalization} morphism
$n: \overline\sY \to \widetilde\sY$.

\medskip
\noindent
\emph{Case 1}. Consider the ring homomorphism
$\kk[u_1, v_1, z] \to \kk[t_1, v, z]$ given by 
$u_1 \mapsto t_1 z, v_1 \mapsto v_1$ and $z \mapsto z$.  
 This induces the following ring homomorphism
of coordinate algebras:
$$
\kk[u_1, v_1, z]/(u_1^2 + z^2 v_1^3 + z^4 v_1^2 + z^5) \lar
\kk[t_1, v_1, z]/(t_1^2 + v_1^3 + v_1^2 + z).
$$
Note that $\overline\sY = V(t_1^2 + v_1^3 + v_1^2 + z) \subseteq \mathbb{A}^3_{(t_1, v_1, z)}$ is
a smooth surface and the morphism of algebraic schemes $n: \overline\sY \to \widetilde\sY$
is the normalization map. Next, 
the preimage of $\widetilde{\sE}_2$ under $n$ is given by the equations
$$
t_1^2 + v_1^3 + v_1^2 = 0, \quad z = 0.
$$
This is a nodal cubic curve, having an $A_1$-singularity at the point $t_1 = v_1 = z = 0$.
Note that $n(t_1, v_1, 0) = v_1$, hence $n: \widetilde{\sE}_2 \to \sE_2$ is a ramified
covering of order two.

\medskip
\noindent
\emph{Case 2}.
In a similar way,  consider the morphism of affine spaces
$n: \mathbb{A}^3_{(u_2, w_2, y)} \lar \mathbb{A}^3_{(t_2, w_2, y)}$ given by the formula
$(u_2, w_2, y) \mapsto (u_2 y, w_2, y)$. This induces the following ring homomorphism
of coordinate algebras:
$$
\kk[u_2, w_2, y)]/(u_2^2 + y^2 w_2 + y^2 w_2^2 + y^3 w_2^5)
\lar
\kk[t_2, w_2, y]/(t_2^2 + w_2 + w_2^2 + y w_2^5).
$$
Again,  $\overline\sY = V(t_2^2 + w_2 + w_2^2 + y w_2^5) \subseteq \mathbb{A}^3_{(t_2, w_2, y)}$ is
a smooth surface and the morphism of algebraic schemes $n: \overline\sY \to \widetilde\sY$
is the normalization map.
In this chart, the preimage $\widetilde\sE_2$  of the exceptional divisor
$\sE_2$  is given by the equation
$$
t_2^2 + w_2 + w_2^2 = 0, \quad y = 0.
$$
Moreover, the morphism $n$ induces an isomorphism of $\sE_1$ on its preimage $\widetilde{\sE}_1$, given by the
equation
$$
t_2 = w_2 = 0.
$$

\medskip
\noindent
\emph{Step 4}. So far, we have constructed a smooth surface $\overline\sY$ and
a birational isomorphism $\pi: \overline\sY \to \sX$, given
by a sequence of  projective birational isomorphisms
$$\overline\sY \stackrel{n}\lar  \widetilde\sY \stackrel{p}\lar  \sY \stackrel{q}\lar  \sX.
$$
Moreover, $\pi^{-1}(o) = \widetilde{\sE}_1 \cup \widetilde{\sE}_2$ is a reducible curve with
two rational components intersecting transversally, where $\widetilde{\sE}_1 \cong \mathbb{P}^1$ and
$\widetilde{\sE}_2$ is isomorphic to a plane projective nodal cubic curve.

\medskip
Our next goal is to determine the self-intersection numbers of 
the divisors $\widetilde{\sE}_1$
and $\widetilde{\sE}_2$. To do  this, take the function
$$f = x \in \kk[X] = \kk[x,y,z]/(x^2 + y^3 + z^7 + y^2 z^2).$$
 Then a straightforward computation shows
that the function $\tilde{f} := \pi^*(f) \in \kk[\overline{\sY}]$ is given by the following formulae.

\medskip
\noindent
\emph{Case 1}. In the coordinates from Step 3 (Case 1) we have: $\tilde{f} =
\pi^*(f) = z^3 t_1$.
Note that  the order of $\tilde{f}$ in the local ring
$\kk[x,y,z]_{\langle z \rangle}/(t_1^2 + v_1^3 + v_1^2 + z)$ is three, where $\langle z\rangle \subset
\kk[x,y,z]$ is the prime ideal generated by $z$.  Moreover, the exceptional fiber
$\widetilde{\sE}_2 = V(t_1^2 + v_1^3 + v_1^2, z)$ intersects  the curve given by the
equation  $t_1 = 0$ at two points of $\mathbb{A}^3_{(t_1, v_1, z)}$. The first of them is
$(0, 0, 0)$, which is the singular point of $\widetilde{\sE}_2$. The second is
$(0, -1, 0)$. Note that the intersection of these two curves is transversal
at the second point  and has  multiplicity two at the first one.
Hence, in this chart
$$
\mathsf{div}(\tilde{f}) = 3 \sE_2 + \sC
$$
and $\sC \cdot \sE_2 = 2 +1 = 3$.

\medskip
\noindent
\emph{Case 2}. In the coordinates from Step 3 (Case 2)  we have: $\tilde{f} = \pi^*(f) =  y^3 t_2 w_2$.
Recall that the exceptional divisor  $\widetilde{\sE}_1$ in this chart  is
described as follows:
$$\widetilde\sE_1 =  V(t_2, w_2) \subseteq
\overline{\sY} = V(t_2^2 + w_2 + w_2^2 + y w_2^5) \subseteq \mathbb{A}^3_{(t_2, w_2, y)}.
$$
Moreover, the second exceptional divisor $\widetilde{\sE}_2$ is given by the formula
$$
\widetilde{\sE}_2 = V(y, t_2).
$$
Since $t_2 = w_2(1 + w_2 + y w_2^4)$, the curve $\widetilde\sE_1$ is locally given by the equation
$t_2 = 0$, and $\tilde{f}$ vanish on $\widetilde{\sE}_2$ with multiplicity three.
 Note that the strict transform of $f$ is given by the polynomial
$g= 1 + w_2 + y w_2^4$. Observe that  $V(g)$ and $\sE_1$ do not intersect,
whereas
$$
V(g) \cap \sE_1 = (0, -1, 0) \in \mathbb{A}^3_{(t_2, w_2, y)}.
$$
However,  it is the same point as the one found in Case 1! Summing everything up, we obtain:
\begin{itemize}
\item We have constructed
a projective birational isomorphism $\pi: \overline\sY \to \sX$,
 whose exceptional locus
consists of the  union of two rational curves $\widetilde{\sE}_1$ and
$\widetilde{\sE}_2$. Moreover,
 $\widetilde{\sE}_1 \cong \mathbb{P}^1$ and
$\widetilde{\sE}_2$ is isomorphic to a plane nodal cubic curve.
 These  curves intersect transversally
at one smooth point.
\item We have: $\mathsf{div}(\tilde{f}) =
\mathsf{div}\bigl(\pi^*(f)\bigr) = 3 \widetilde{\sE}_1 + 3 \widetilde{\sE}_2 + \sC.$
\item  The curve $\sC$ does not intersect $\widetilde{\sE}_1$. Moreover, it intersects
 $\widetilde{\sE}_2$ with multiplicity two
at  its  singular point  and transversally at another point (which is smooth).
All-together, this implies:
$$
\left\{
\begin{array}{ccc}
\sC \cdot \widetilde{\sE}_1 & =  & 0 \\
\sC \cdot \widetilde{\sE}_2 & =  & 3 \\
\widetilde{\sE}_1 \cdot \widetilde{\sE}_2 & =  & 1. \\
\end{array}
\right.
$$
\end{itemize}

\noindent
These computations imply that
$$
0 = \mathsf{div}\bigl(\pi^*(f)\bigr) \cdot \widetilde{\sE}_1 = 3\widetilde{\sE}_1^2 + 3 = 0,
$$
hence $\widetilde{\sE}_1^2 = -1$. In a similar way,
$$
0 = \mathsf{div}\bigl(\pi^*(f)\bigr) \cdot \widetilde{\sE}_2 = 3\widetilde{\sE}_1^2 + 6 = 0,
$$
hence $\widetilde{\sE}_1^2 = -2$.

\medskip
\noindent
\emph{Step 5}.
We have constructed a resolution of singularities $\pi: \overline{\sY} \to \sX$,
whose exceptional divisor is $\widetilde\sE = \widetilde{\sE}_1 \cup \widetilde{\sE}_2$. However,
since $\sE_1 \cong \mathbb{P}^1$ and $\widetilde{\sE}_1^2 = -1$, this resolution is not minimal! Hence, using Castelnuovo's
theorem (see \cite[Theorem V.5.7]{Hartshorne})
we can blow down $\widetilde{\sE}_1$ and obtain a commutative diagram
$$
\xymatrix
{
\overline{\sY} \ar[rr]^\phi \ar[dr]_\pi & & \widetilde{\sX} \ar[dl]^{\pi'} \\
& \sX &
}
$$
where the surface $\widetilde{\sX}$ is smooth and
$\phi: \overline{\sY} \setminus \widetilde{\sE}_1 \lar \widetilde{\sX}
\setminus \phi(\widetilde{\sE}_1)$ is an isomorphism.
Note that the morphism $\phi$ induces a ring homomorphism of  Chow groups
$\phi^*:  A^*(\widetilde{\sX}) \to A^*(\overline{\sY})$. Hence, if $\sD = \phi(\widetilde{\sE}_2)$ is the exceptional divisor of the resolution of singularities $\pi': \widetilde{\sX} \to \sX$, then
$$
\sD^2 = (\widetilde{\sE}_1 + \widetilde{\sE}_2)^2 = \widetilde{\sE}_1^2 + 2
\widetilde{\sE}_1 \cdot \widetilde{\sE}_2 + \widetilde{\sE}_1^2 = -1.
$$

\noindent
Summing everything up, the singularity
$(\sX, o) = \bigl(V(x^2 + y^3 + z^7 + y^2 z^2), o\bigr) \subseteq (\mathbb{A}^3, o)$
has a minimal resolution $\pi': (\widetilde\sX, \sD) \to \sX$ such that
$\sD$ is isomorphic to the plane projective cubic curve
$zy^2 = x^3 + x^2 z$. Moreover, the self-intersection index of $\sD$ is $-1$.
\qed

\end{example}

\medskip
\begin{definition}\label{D:resolofsing}
Let $(\sX, o)$ be the germ of a surface singularity and
$\pi: (\widetilde\sX, \sE) \to (\sX, o)$ its minimal resolution.
A  \emph{fundamental cycle} is the minimal
 cycle $\sZ = m_1 \sE_1 +
m_2 \sE_2 + \dots + m_n \sE_n$,  such that
$m_i > 0$ and $(\sZ, E_i) \le 0$ for all $1 \le i \le n$. By a result
of Artin, see \cite[Proposition 2]{Artin}, a fundamental cycle \emph{exists} and is \emph{unique}.
\end{definition}

\medskip
\begin{definition}
In the notations of Definition \ref{D:resolofsing}
the singularity $(X, o)$ is called \emph{rational} if
$R^1\pi_{*}\kO_{\widetilde\sX} = 0$ (it is equivalent to
$H^1(\kO_{\widetilde\sX}) = 0$).
Note that since $\sX$ is normal, we have
$\pi_*\kO_{\widetilde\sX} =
\kO_\sX$.
\end{definition}

\medskip
\begin{proposition}
Let $(\sX, o)$ be a rational normal surface singularity,
$\pi: (\widetilde\sX,\sE) \to (\sX, o)$ its
minimal resolution and
$\sE_1, \sE_2, \dots, \sE_n$  be the  irreducible components of
$\sE$. Then we have:
\begin{itemize}
\item All components $\sE_i$ are smooth and rational,
i.e. $\sE_i \cong \mathbb{P}^1$;
\item For all $i \ne j \ne l$ we have:
$\sE_i \cap \sE_j \cap \sE_l = \emptyset$.
\item For all $1 \le i \ne j \le n$ we have:
$\sE_i \cdot \sE_j  \in \{0, 1\}$, i.e. $\sE$ is a configuration
of projective lines intersecting transversally.
\item
Moreover $\sE$ has no cycles, i.e. it is a tree of projective lines.
\end{itemize}
\end{proposition}

\noindent
A proof of this Proposition can be found in \cite[Lemma 1.3]{Brieskorn}.

\medskip
\begin{remark}
Note that there exist non-rational normal surface singularities
such that
$\sE$ is a tree of projective lines, see \cite[Proposition 3.5]{Laufer}.
\end{remark}

\medskip
\noindent
Knowing the intersection matrix $M(\sX)$ of a normal
surface singularity $(\sX, o)$ one can determine  its local
fundamental group. In particular, in the case of rational normal
surface singularities we have:

\begin{theorem}[Brieskorn]
The local fundamental group $G = \pi_1(\sX, o)$ of a rational normal surface singularity
$(\sX, o)$ can be presented by the following generators and relations:
$$
G =  \Bigl\langle g_1, g_2,\dots,g_n\, \Big| \,
g_i g_j^{m_{ij}} = g_j^{m_{ij}} g_i, 1 \le i, j \le n, 
g_1^{m_{i1}} g_2^{m_{i2}}
\dots g_n^{m_{in}} = e, 1 \le i \le n\Bigr\rangle,
$$
where $M = M\bigl(\sX, o\bigr) = (m_{ij}) \in \GL_n(\mathbb{Z})$ is the intersection matrix of $(\sX, o)$.
\end{theorem}

\noindent
For a proof of this theorem we refer to \cite[Lemma 2.7]{Brieskorn} and
\cite[Chapter 1]{Mumford}.

\medskip
\begin{lemma}\label{L:resolofratsing}
Let $(\sX, o)$ be a rational singularity, $\sE = \cup_{i=1}^n \sE_i$.
 Then  $\Pic(\widetilde\sX) \cong\ H^2(\widetilde{\sX}, \mathbb{Z})\cong \mathbb{Z}^n$.
\end{lemma}

\noindent
\emph{Proof}. Consider the exponential sequence
$
0 \to \mathbb{Z}_{\widetilde\sX} \to
\kO_{\widetilde\sX}  \stackrel{\exp}\lar
 \kO^*_{\widetilde\sX} \to 0.
$
Since $H^2(\kO_{\widetilde{\sX}}) =
H^1(\kO_{\widetilde{\sX}}) =  0$, taking cohomology  we get isomorphisms
$$
\Pic(\widetilde\sX) \cong H^1(\kO^*_{\widetilde\sX})
\stackrel{c_1}\lar H^2(\widetilde{\sX}, \mathbb{Z})
\stackrel{\phi}\lar \mathbb{Z}^n
$$
where $c_1(\kL)$ is the first Chern class of a line bundle
$\kL$ and the isomorphism $\phi$ is defined by the rule
$$\phi(D) = (D \cdot E_1, D \cdot E_2,\dots, D \cdot E_n) \in \mathbb{Z}^n.$$
\qed

\medskip

\begin{remark}
Lemma \ref{L:resolofratsing} essentially means   that the divisors
on  a  minimal resolution of a rational singularity can be described
\emph{topologically}.
\end{remark}

\medskip
\noindent
At this point we would like to have some examples
of rational normal surface singularities.

\begin{lemma}
Let $G \subset \GL_2(\mathbb{C})$ be a finite  small subgroup,
$\rA = \mathbb{C}\{x_1, x_2\}^G$  and $(X, o)$ the corresponding
complex germ. Then  the singularity $(X, o)$ is rational.
\end{lemma}

\noindent
For a proof of this Lemma, see \cite{Artin} and \cite{Brieskorn}.

\medskip
\noindent
The following proposition says when a rational normal surface singularity
is a complete intersection.

\begin{theorem}[Wahl]
Let $(\sX, o)$ be a complex rational normal surface singularity, which
is a complete intersection and $(\rA, \idm)$ be
the corresponding local analytic algebra. Then
$\rA \cong  \mathbb{C}\{x_1, x_2\}^G$, where $G \subset \SL_2(\mathbb{C})$
is a finite subgroup.
\end{theorem}

\noindent
For a proof of this theorem,  see \cite[Theorem 2.1]{Wahl}.

\medskip
In other words, normal rational singularities, which are complete
intersections, are exactly simple hypersurface singularities. Moreover,
let $\Delta = \Delta_n$ be the Dynkin type
of $(\sX, o)$, then
\begin{itemize}
\item For any irreducible component $\sE_i$ of $\sE$  we have:
$\sE_i^2 = -2$.
\item Moreover, the quadratic form on
$\mathbb{Z}^n \times \mathbb{Z}^n  \to \mathbb{Z}$ defined
by the intersection matrix $M(\sX)$ is equal to  $-q_\Delta$, where
$q_\Delta$ is Tits quadratic form attached to $\Delta$.
\item The fundamental cycle
$\sZ = m_1 \sE_1 + m_2 \sE_2 + \dots + m_n \sE_n$ is the
\emph{maximal root}  of
the   quadratic form of $M(\sX)$. In particular,  $\sZ^2  = -1$.
\end{itemize}

\medskip

\begin{lemma}[Artin]
More generally, the embedding dimension of a normal rational
surface singularity $(\sX, o) = \Specan(\rA,\idm)$ can be computed
by the formula
$$\edim(\sX, o) := \dim_{\mathbb{C}}(\idm/\idm^2) = - \sZ^2 + 1.$$
\end{lemma}

\noindent
For a proof of this formula,  see \cite[Corollary 6]{Artin}.

\medskip
\begin{remark}
The fundamental cycles of  the simple surface singularities 
 are described below.
Bullets in the underlying diagrams denote irreducible components,
they are connected by an edge if the corresponding components intersect.
Integers drawn over bullets denote coefficients of these components in the
formula for the fundamental  cycle. The obtained graph is called
the \emph{dual intersection graph}  of the isolated surface singularity.
 In all listed cases, if
$\Delta = \Delta_n$ is the Dynkin type of the quotient singularity,
then  the number of irreducible components of the exceptional
divisor of the minimal resolution is $n$.
\begin{itemize}
\item $A_{n}$--singularity $x^2 + y^{n+1} + z^2 \ (n \ge 1)$ has the following
dual intersection graph and the fundamental cycle:
$$
\xymatrix
{\stackrel{1}\bullet \ar@{-}[r] & \stackrel{1}\bullet \ar@{-}[r] & \stackrel{1}\bullet \ar@{-}[r] & \stackrel{1}\bullet &  \dots
& \stackrel{1}\bullet \ar@{-}[r] & \stackrel{1}\bullet
}
$$
\item In the case of $D_n$--singularities  $
x^2 y + y^{n-1} + z^2 \ (n \ge 4)$ we have
$$
\xymatrix@R=.5em
{
& & &  &
  &  &  \stackrel{1}\bullet\\
\stackrel{1}\bullet \ar@{-}[r] & \stackrel{2}\bullet \ar@{-}[r]
& \stackrel{2}\bullet & \dots & \stackrel{2}\bullet \ar@{-}[r] &
\stackrel{2}\bullet \ar@{-}[ru]
\ar@{-}[rd]  & \\
& & &  &  & &  \stackrel{1}\bullet
}
$$
\item $E_6$--singularity $x^4 + y^3 + z^2$:
$$
\xymatrix@R=1em
{\stackrel{1}\bullet \ar@{-}[r] & \stackrel{2}\bullet \ar@{-}[r] &
\stackrel{3}\bullet \ar@{-}[r] \ar@{-}[d]
&  \stackrel{2}\bullet \ar@{-}[r] & \stackrel{1}\bullet \\
& & \stackrel{2}\bullet  & &  \\
}
$$
\item $E_7$--singularity $x^3 y + y^3 + z^2$:
$$
\xymatrix@R=1em
{\stackrel{2}\bullet \ar@{-}[r] & \stackrel{3}\bullet \ar@{-}[r] &
\stackrel{4}\bullet \ar@{-}[r] \ar@{-}[d]
&  \stackrel{3}\bullet \ar@{-}[r] & \stackrel{2}\bullet
 \ar@{-}[r] & \stackrel{1}\bullet \\
& & \stackrel{2}\bullet  & & &  \\
}
$$
\item Finally, in the case of $E_8$--singularity $x^5  + y^3 + z^2$
we have:
$$
\xymatrix@R=1em
{\stackrel{2}\bullet  \ar@{-}[r] & \stackrel{4}\bullet
 \ar@{-}[r] &  \stackrel{6}\bullet \ar@{-}[r] \ar@{-}[d]
&  \stackrel{5}\bullet \ar@{-}[r] & \stackrel{4}\bullet
 \ar@{-}[r] & \stackrel{3}\bullet
\ar@{-}[r] & \stackrel{2}\bullet  \\
& & \stackrel{3}\bullet  & & &  & \\
}
$$
\end{itemize}
\end{remark}

\medskip
\begin{theorem}[Gonzalez-Springer--Verdier, Artin--Verdier]\label{T:McKayCorr}
Let $G \subset \SL_2(\mathbb{C})$ be a finite subgroup,
$\rA = \mathbb{C}\{x_1, x_2\}^G$ and $(X, o) = \Specan(\rA)$
 the corresponding quotient singularity,
$\pi: (\widetilde\sX, \sE) \to (X, o)$ its  minimal resolution,
$\sE = \sE_1 \cup \sE_2 \cup \dots \cup \sE_n$  the decomposition
of $\sE$ into a union of irreducible components and
$\sZ$ the fundamental  cycle.
Let
$\mM$ be an indecomposable Cohen-Macaulay $\rA$--module,
$\kM$ the corresponding coherent sheaf on $\sX$ and
$\widetilde\kM = \pi^*\kM/\mathsf{tors}(\pi^*\kM)$. Then the following
properties hold:
\begin{itemize}
\item $\widetilde\kM$ is a vector bundle on $\widetilde\sX$;
\item Let $\sD = c_1(\widetilde\kM)$. If $\mM \not\cong \rA$ then
there exists $1 \le i \le n$ such that $D = \sE_i^*$, i.e
$
D \cdot E_j = \delta_{i,j}.
$
Moreover, $\rank(\mM) = \sZ \cdot  \sD$.
\end{itemize}
\end{theorem}

\noindent
For a proof of this theorem, see \cite[Theorem 1.11]{ArtinVerdier}.

\medskip

\begin{remark}
Theorem \ref{T:McKayCorr} implies that a Cohen-Macaulay $\mM$
on a simple surface singularity $(\rA, \idm)$ is determined
by two discrete parameters: $\rank(\mM)$ and
$D = c_1(\widetilde\kM) \in H^2(\widetilde\sX, \mathbb{Z})
\cong \mathbb{Z}^n$, where $n$ is the number of irreducible
components of the exceptional divisor $\sE$.  The cohomology class
$D$ satisfies the following property: $D \cdot E_i \ge 0$ for all
$1 \le i \le n$.
\end{remark}

\medskip
\begin{corollary}
Let $G \subset \SL_2(\mathbb{C})$ be a finite subgroup,
$\rA = \mathbb{C}\{x_1, x_2\}^G$ and  $(\sX, o)$ the corresponding
complex singularity. Let  $\pi: (\widetilde\sX, \sE) \to
(\sX, o)$ be  a  minimal resolution of singularities.
Combining Theorem \ref{T:McKayCorr} and
Theorem \ref{T:CMquot}, we obtain a bijection
between the following three
sets:
\begin{enumerate}
\item Isomorphy classes of indecomposable non-free Cohen-Macaulay $\rA$-modules;
\item Isomorphy classes of non-trivial 
irreducible representations of $G$;
\item Irreducible components of the exceptional divisor $\sE$.
\end{enumerate}
\end{corollary}

\medskip
\begin{remark}
A certain generalization of Theorem \ref{T:McKayCorr} in the case
of arbitrary cyclic quotient singularities was obtained by
Wunram \cite{Wunram}. The general case of Cohen-Macaulay modules on
arbitrary quotient surface singularities
was studied in a work of Esnault \cite{Esnault}, see also \cite[Chapter 4]{Kahn2}.
\end{remark}

\medskip
\section{Cohen-Macaulay modules over minimally elliptic singularities and
vector bundles on genus one curves}

Two-dimensional Gorenstein quotient singularities are exactly
the two-dimensional hypersurface singularities of modality zero \cite{Arnold}.
This in particular means that
a deformation of a simple surface singularity is again a simple surface
singularity.

The next interesting class of surface singularities is formed by
singularities of modality one. The most interesting among them
are the so-called \emph{minimally elliptic} singularities.

\medskip

\begin{definition}[Laufer, see Definition 3.2 in   \cite{Laufer2}]
A normal surface singularity $(\sX, o)$ is called minimally elliptic
if it is Gorenstein and $H^1(\kO_{\widetilde\sX}) = \kk$, where
$\pi: (\widetilde\sX, \sE) \to (\sX, o)$ is
its minimal resolution.
\end{definition}

\medskip
\begin{theorem}[Laufer]
Let $(\sX, o)$ be a minimally elliptic singularity,
$\pi: (\widetilde\sX, \sE) \to (\sX, o)$ its minimal resolution of
singularities. Then $\sE$ is a configuration of projective lines with transversal
 intersections and a tree as its dual graph, except the following cases:
\begin{enumerate}
\item $\sE$ is an elliptic curve; then $(\sX, o)$ is called
\emph{simple elliptic}.
\item $\sE$ is a cycle of $n$ projective lines (a plane nodal
cubic curve $zy^2 = x^3 - x^2z$ for $n = 1$); then, following
Hirzebruch \cite{Hirzebruch}, the corresponding singularity is called
a \emph{cusp}.
\item $\sE$ is a cuspidal cubic
curve $zy^2 = x^3$, a tachnode curve $(yz - x^2)(yz + x^2) = 0$ or
three concurrent lines in a plane $xy(x - y) = 0$.
\end{enumerate}
\end{theorem}

\noindent
For a proof of this theorem, see \cite[Proposition 3.5]{Laufer2}.

\medskip
\begin{example}
The surface singularity, given by the equation $x^2 + y^3 + z^7 = 0$ is minimally elliptic.
In this case the exceptional divisor is isomorphic to a cuspidal plane
cubic curve $zy^2 = x^3$. 

\medskip
\noindent
The surface singularity $x^2 + y^3 + z^7 + xyz = 0$ is a cusp.
Its resolution of singularities was described in Example \ref{E:resol2}.
In particular, its exceptional divisor is a nodal 
plane cubic curve $zy^2 = x^3 + x^2 z$ having self-intersection 
index $-1$.  

\medskip
\noindent
Other examples of minimally elliptic singularities of small multiplicities
can be found in \cite[Tables 1, 2 and 3]{Laufer2}.
\end{example}

\medskip

\begin{theorem}[Laufer]
Let $(\sX, o)$ be a minimally elliptic singularity,
$(\widetilde\sX, \sE) \to (\sX, o)$ a minimal resolution of singularities
and $\sZ$ the fundamental cycle. Then
$$\edim(\sX, o)  = \mathsf{max}\{ - \sZ^2, 3\}.
$$
\end{theorem}

\noindent
For a proof of this theorem, see \cite[Theorem 3.13]{Laufer2}.

\medskip

\begin{theorem}
Let $(\sX, o)$ be a simple elliptic singularity,
$(\widetilde\sX, \sE) \to (\sX, o)$ its minimal resolution. Then
the analytic type of the singularity $(\sX, o)$ is determined
by the analytic type of the elliptic curve $\sE$
and its self-intersection number $\sE^2 = -b, \quad  b \ge 1$. Moreover,  $\widetilde\sX$ is locally
(in a neighborhood of $\sE$) isomorphic to the  total
space of the line  bundle $\big|\kO_\sE(-b[p_0])\big|$, where $p_0$ is
a point corresponding to the zero element of $\sE$ viewed as
an algebraic group.
\end{theorem}

\noindent
For a proof of this result, see \cite[Korollar 1.4]{Saito}.

\medskip

\begin{corollary}
All simple elliptic singularities are quasi-homogeneous.
\end{corollary}

\medskip

\begin{example}
Simple elliptic singularities of low multiplicities
are given by the following equations, see \cite[Satz 1.9]{Saito}.
\begin{itemize}
\item $\mathsf{El}(\sE, 1):  z^2 = y(y-x^2)(y - a x^2)$;
\item $\mathsf{El}(\sE, 2):  z^2 = xy(x-y)(x - ay)$;
\item $\mathsf{El}(\sE, 3):  z^2 = y(y-x)(y - a x)$;
\item $\mathsf{El}(\sE, 4):  z^2 = y(w - (a+1)y + ax), \quad  y^2 = xw$;
\end{itemize}
In all these cases $a \in
\mathbb{C}\setminus \{0, 1\}$ and
$$
j(\sE) = \frac{4}{27}\frac{(a^2 - a + 1)^3}{a^2(a-1)^3}
$$
is the $j$--invariant of the elliptic curve $\sE$. All other
simple elliptic singularities are no longer complete intersections.
\end{example}

\noindent

\begin{remark}
Following  Saito \cite{Saito}, there is an alternative notation
for those simple elliptic singularities which are complete
intersections. Namely, let $(\sX, o)$ be a simple elliptic singularity of type
$\mathsf{El}(\sE, i), \quad 1 \le i \le 4$ and
$\mathcal{X} \to \mathcal{B}$ its semi-universal deformation.
For  $t \in  \mathcal{B}\setminus \{0\}$ the  fiber
$\mathcal{B}_t$ (called \emph{Milnor fiber})
is homotopy equivalent to a bouquet of spheres and
$H^1(\mathcal{B}_t, \mathbb{Z}) \cong H^2(\mathcal{B}_t, \mathbb{Z})
\cong \mathbb{Z}^\mu,$ where $\mu$ is the Milnor number of $(\sX, o)$.
Moreover, the intersection form
$H^1(\mathcal{B}_t, \mathbb{Z}) \times
 H^2(\mathcal{B}_t, \mathbb{Z}) \to \mathbb{Z}$ is non-negatively
  definite with radical of rank two and
\begin{itemize}
\item for $\mathsf{El}(\sE, 1)$ it has Dynkin type $\widetilde{E}_8$;
\item for $\mathsf{El}(\sE, 2)$ it has Dynkin type $\widetilde{E}_7$;
\item for $\mathsf{El}(\sE, 3)$ it has Dynkin type $\widetilde{E}_6$;
\item for $\mathsf{El}(\sE, 4)$ it has Dynkin type $\widetilde{D}_4$.
\end{itemize}
\end{remark}

\begin{proposition}
Let $(\sX, o)  = \mathsf{El}(\sE, b)$ be a simple elliptic singularity
of analytic type $(\sE, b)$. Then its local fundamental group
has the following presentation:
$$
\pi_1(\sX, o) = H_b = \bigl\langle \alpha, \beta, \gamma \, \big|\,
[\alpha,\gamma] = 1, [\beta, \gamma]= 1, [\alpha, \beta] =
\gamma^b\bigr\rangle.
$$
In group theory,  it is called \emph{discrete Heisenberg group}.
\end{proposition}

\noindent
For a proof of this proposition we refer to \cite[Chapter 1]{Mumford} and \cite[Proposition 6.2]{Kahn2}.

\medskip

\begin{theorem}
The analytic type of a cusp singularity is determined by its intersection
matrix. Moreover, a cusp singularity is a complete intersection only
in the following cases:
\begin{itemize}
\item $T_{p,q,r}$--singularities
given by the equation $x^p + y^q + z^r - xyz, \quad \frac{1}{p} +
\frac{1}{q} +\frac{1}{r} < 1$.
\item $T_{p,q,r,t}$--singularities given by two equations
$x^p + y^q = uv, \quad u^r + v^t = xy$, where $p,q,r,t \ge 2$ and
$\max(p,q,r,t) \ge 3$.
\end{itemize}
\end{theorem}

\noindent
This theorem follows from a general result of Wahl \cite[Theorem 2.8]{Wahl} about syzygies
of minimally elliptic singularities.

\medskip
Inspired by results on the geometric McKay Correspondence
for quotient surface singularities, Kahn has proven the following
beautiful theorem on Cohen-Macaulay modules over minimally
elliptic singularities.

\begin{theorem}[Kahn]\label{T:Kahn}
Let $(\sX, o)$ be a minimally  elliptic singularity, $(\rA, \idm)$
the corresponding analytic algebra,
$\pi: (\widetilde\sX, \sE) \to (\sX, o)$ its minimal resolution
 and $\sZ$ the fundamental cycle. Then
the functor
$$
\mathbb{F}: \CM(\rA) \to \VB(\sZ), \quad \mathbb{F}(\mM) =
\pi^*(\mM)^{\vee\vee}|_{\sZ}
$$
preserves iso-classes of objects, i.e. $\mathbb{F}(\mM) \cong
\mathbb{F}(\mM')$ implies $\mM \cong \mM'$. Moreover, the
Cohen-Macaulay modules correspond to the following vector bundles.
\begin{itemize}
\item The regular  Cohen-Macaulay module $\rA$ corresponds
to the vector bundle $\kO_\sZ$.
\item Cohen-Macaulay modules without free direct summands
correspond to vector bundles of the form $\kV \oplus \kO^n$ where
\begin{enumerate}
\item $\kV$ is generically spanned by its global sections, i.e.
the cokernel of the evaluation map $H^0(\kV) \otimes_\kk \kO \to
\kV$ is a skyscraper sheaf;
\item we have: $H^1(\kV) = 0$;
\item finally, $n = h^0\bigl(\kV \otimes 
\kO_{\widetilde\sX}(\sZ)\big|_{\sZ}\bigr)$.
\end{enumerate}
\end{itemize}
\end{theorem}

\noindent
This result is proven in \cite[Theorem 2.1]{Kahn}.

\medskip
\begin{remark} The functor $\mathbb{F}$ from Theorem \ref{T:Kahn}
is neither exact nor faithful. It is not known whether it is full
or not.
\end{remark}

\medskip
\begin{remark}
The result of Kahn also remains valid
in  the case of complete rings over an algebraically closed
field $\kk$ of arbitrary characteristic, see \cite{DGK}.
\end{remark}

\noindent
Hence, a description of Cohen-Macaulay modules over minimally elliptic singularities reduces
to the classification of indecomposable vector bundles
(with some restrictions) on projective curves
 (maybe non-reduced) of arithmetic genus one. Note that if a minimally
elliptic surface singularity is either simple elliptic or cusp then
its fundamental  cycle $\sZ = \sE$ and we deal with vector bundles
on reduced  curves.

\medskip
\begin{theorem}[Atiyah]\label{T:Atiyah}
Let $\kk$ be an algebraically closed field and
$\sE$ be an elliptic curve over $\kk$, $\kV$ an indecomposable vector bundle
on $\sE$. Then $\kV$ is semi-stable and is
uniquely determined by its rank $r$, degree
$d$ and a point $x$ of the curve $\sE$. To be more precise, let
$\kV_0$ be the uniquely determined Jordan-H\"older factor  of
$\kV$, then $\det(\kV_0)= \kO_\sE\bigl((d-1)[x] + [x-x_0]\bigr)$, where
$x_0$ is the point corresponding to the zero
element of $\sE$ viewed as an algebraic group.
\end{theorem}

\noindent
For a proof, see \cite[Theorem 7]{Atiyah}.

\medskip

\begin{theorem}[Oda]
Let $\sE = \sE_\tau = \mathbb{C}/\langle 1, \tau\rangle$  be a complex
torus, $\kV$ an indecomposable holomorphic vector bundle
on $\sE$ of rank $r$ and degree $d$. If $h = \gcd(r,d)$, $r = r' h$ and
$d = d' h$  then there exists a line bundle $\kL$ on the elliptic curve
$\sE_{r'\tau}$ such
that
$$
\kV \cong  \pi_*(\kL) \otimes \kA_h,
$$
where $\pi: \sE_{r'\tau} \to \sE_{\tau}$ is an \'etale covering
of degree $r'$ and $\kA_h$ is the unipotent vector bundle of rank
$h$, i.e. the vector bundle recursively defined by the following
non-split sequences:
\begin{equation}\label{E:unipotent}
0 \to \kA_{h-1} \to \kA_{h} \to \kO \to 0, \quad h \ge 2, \quad
\kA_1 = \kO.
\end{equation}
Note that we automatically have: $\Ext^1_\sE(\kO, \kA_h) \cong
 H^1(\kA_h) =  \mathbb{C}$, hence the middle term
 of the short exact sequence (\ref{E:unipotent}) is uniquely determined.
\end{theorem}

\noindent
For a proof of this theorem, see \cite[Proposition 2.1]{Oda}.

\medskip
\begin{lemma}
Kahn's functor  $\mathbb{F}$ maps the fundamental module $\mD$ to the
Atiyah's bundle $\kA_2$.
\end{lemma}

\noindent
A proof goes along the same lines as  in \cite[Theorem 3.1]{Kahn}. Note that,
since a simple elliptic singularity $(\sX, o) = \Specan(\rA)$
is quasi-homogeneous, we have:
$\mD \cong (\Omega^1_\rA)^{**}$.

\medskip
\begin{remark}
Since cusp singularities are not quasi-homogeneous,  in that case
$\mD \not\cong (\Omega^1_\rA)^{**}$.
A description of the fundamental module $\mD$ in this case is due to Behnke, see
\cite[Section 5]{Behnke}.
\end{remark}

\noindent

\begin{theorem}[Kahn]
Let $(\sX, o)$ be a simple elliptic singularity, $(\rA, \idm)$ the corresponding local ring,
 $\pi: (\widetilde\sX, \sE) \to (\sX, o)$
its  minimal resolution and  $b = \mathsf{max}\{2, -\sE^2\}$ its
 multiplicity. Then the indecomposable Cohen-Macaulay
$\rA$-modules are parameterized by two discrete parameters:
$$\Bigl\{(r,d) \in \mathbb{Z}^2 \, \Big| \, 1 \le r \le d \le (b+1)r\Bigr\}$$
and one continuous parameter $\lambda \in \sE$.

Recall that $\widetilde\sX$ can be locally written  as the total space
of a line bundle on $\sE$, i.e.
$$\widetilde\sX  =  \big|\kO_\sE(-b[p_0])\big|,$$
 where $p_0 \in \sE$ is some fixed point, let
$p: \widetilde\sX \to \sE$ be the corresponding projection. Consider the commutative diagram
$$
\xymatrix
{
  & \sE & \\
\widetilde{\sX}\setminus\sE \ar@{^{(}->}[rr] \ar[ru]^s \ar[d]_q &
& \widetilde\sX \ar[lu]_p \ar[d]^\pi \\
\sX\setminus\{0\} \ar@{^{(}->}[rr]^i & & \sX
}
$$
Let $\kE(r,d,\lambda)$ be the indecomposable vector bundle on $\sE$ of rank $r$, degree $d$
and continuous parameter $\lambda$, constructed in Theorem \ref{T:Atiyah}.
 Then
$$
\mM(r,d,\lambda) =
\Bigl\{\Gamma\Bigl(i_* q_* s^* \kE(r,d,\lambda)\Bigr) \, \Big| \, 
 1 \le r \le d \le (b+1)r, \, \lambda
\in \sE
\Bigr\}
$$
is the complete list of indecomposable Cohen-Macaulay $\rA$-modules.
\end{theorem}

\noindent
For a proof, see  \cite[Proposition 5.18]{Kahn2}.

\medskip
\noindent
Unfortunally, Kahn's description of Cohen-Macaulay modules
over minimally elliptic singularities is not really explicit.
In some cases, however, one can find families of matrix
factorizations.

\begin{example}[Etingof-Ginzburg]
Consider the simple elliptic singularity of type $\widetilde{E}_6$ given by the
equation $w = x^3 + y^3 + z^3  + \tau xyz, \quad \tau \in
\mathbb{C}^*$.  
Let $\sE = \sE_\tau \subset \mathbb{P}^2$ be the corresponding elliptic curve and
$(a:b:c) \in \sE$ a point such that its components
$a, b$ and $c$ are non-zero. Then the matrix
$$
\phi =
\left(
\begin{array}{ccc}
ax & b y & c z \\
cz & ay &  b x \\
by & cz & az
\end{array}
\right)
$$
and its adjoint
$$
\phi'=
\left(
\begin{array}{ccc}
a^2 yz - bc x^2 &  c^2 xy - ab z^2  & b^2 xy  - ac y^2  \\
b^2 xy - ac z^2 &  a^2 yz - bc y^2  & c^2 yz  - ab x^2  \\
c^2 xz   - ab  y^2 &  b^2 yz  - ac x^2  & a^2 xy - bc z^2  \\
\end{array}
\right)
$$
satisfy the equality:
$$
\phi  \cdot \phi' = \phi' \cdot \phi  = - (abc) w\cdot  \mathsf{id}.
$$
In particular, the matrix factorization $M(\phi, \phi')$ defines
a family of indecomposable Cohen-Macaulay modules of rank one 
\cite[Example 3.6.5]{EtingofGinzburg}.
\end{example}

\medskip
Other examples of matrix factorizations in the case of the simple
elliptic singularity of type $\widetilde{E}_7$ given by the equation
$x^4  + y^4  + \tau x^2 y^2 + z^2 = 0$ were found by Knapp
\cite{Knapp}. Some examples of matrix factorizations for $x^3 + y^3 + z^3 = 0$
corresponding to Cohen-Macaulay modules of higher ranks, were found
by Laza, Pfister and Popescu \cite{LazaandK}.

\medskip
\medskip
Indecomposable vector bundles and torsion free sheaves on Kodaira cycles
of projective lines were classified by Drozd and Greuel, see  \cite[Theorem 2.12]{DGVB}.
In
the case the ground field $\kk$ is algebraically closed of characteristic
zero, their classification
can be presented in the following form, generalizing
Oda's description,  see \cite[Theorem 19]{surveySchwpkt}.

\begin{theorem}
Let $\sE = \sE_n$  be a Kodaira cycle of $n$ projective lines,
 $\kV$ an indecomposable holomorphic vector bundle
on $\sE$ of rank $r$ and degree $d$. Then
there exists an \'etale covering
$\pi: \sE_{nt} \to \sE_n$ of degree $t$, a line bundle $\kL \in
\Pic(\sE_{nt}) = \mathbb{Z}^{nt} \times \kk^*$ and a positive integer
$h \ge 1$
such
that $r = th$ and
$$
\kV \cong  \pi_*(\kL) \otimes \kA_h,
$$
where $\kA_h$ is defined  by short exact sequences (\ref{E:unipotent}) in the
same way as in the case  of  elliptic curves.
\end{theorem}

\medskip

\begin{corollary}
Simple elliptic and cusp singularities as well as their
quotients by a finite group of automorphisms
(the so-called simple
$\mathbb{Q}$--elliptic
and $\mathbb{Q}$--cusp singularities)  have tame Cohen-Macaulay
representation type.
\end{corollary}

\medskip
The following result of Drozd and Greuel implies that the remaining minimally elliptic singularities
are Cohen-Macaulay wild.

\begin{theorem}[Drozd--Greuel]
Let $\sE$ be a reduced projective curve.
\begin{enumerate}
\item If $\sE$ is a chain of projective lines, then $\VB(\sE)$ is of finite type.
\item If $\sE$ is a smooth elliptic curve, then $\VB(\sE)$ is tame of polynomial growth.
\item If $\sE$ is a cycle of projective lines, $\VB(\sE)$ is tame of exponential growth.
\item In all other cases, $\VB(\sE)$ has wild representation type.
\end{enumerate}
\end{theorem}

\noindent
For a proof of this theorem, see \cite[Theorem 2.8]{DGVB}.

\medskip
\begin{remark}
By a result of Kawamata \cite{Kawamata}, simple $\mathbb{Q}$--elliptic and
 $\mathbb{Q}$--cusp singularities  are precisely the \emph{log-canonical}
 singularities. Moreover, results of
 Mumford \cite{Mumford} and Karras \cite{Karras}
imply they have solvable local fundamental group.
 Other way around,  a classification of  Karras, based on an
 earlier result of Wagreich \cite{Wagreich}
 and completed by Kawamata's classification
 of rational log-canonical singularities  implies
 that a normal surface singularity with an infinite
 solvable local fundamental group
 is either simple $\mathbb{Q}$--elliptic or
 $\mathbb{Q}$--cusp.
 \end{remark}

\begin{remark}
The local fundamental group of the $E_8$--singularity is
binary-icosahedral, which is known to be not solvable.
\end{remark}

\section{Other results on Cohen-Macaulay representation type}

In this section,  we briefly mention some other results
and conjectures related to our
 study of  surface singularities of finite and tame 
 Cohen-Macaulay representation types,
see \cite{DrozdGreuel, semicont} for the definition of Cohen-Macaulay tame and
Cohen-Macaulay wild representation type.

The original motivation to classify the indecomposable  Cohen-Macaulay modules
over a Cohen-Macaulay  local ring originates  from the theory of
integral representations of finite groups. For example,
let $\mathbb{Z}_{(2)}$ be the ring
of $2$-adic integers and $G = \mathbb{Z}/8\mathbb{Z}$ be a cyclic group of order $8$.
Then the category
of finitely generated torsion free $\mathbb{Z}_{(2)}[G]$--modules
is tame \cite{Yakovlev}. Moreover, the underlying classification problem
is closely related to   the study of Cohen-Macaulay modules over
a minimally elliptic curve singularity of type $T_{36}$.

More generally, one can pose a question about the representation type of the
category of \emph{lattices} over a complete \emph{order}. However, in this framework tameness
results mainly concern the case of orders of Krull dimension one, see surveys
of Dieterich  \cite{Dieterich} and Simson \cite{Simson} for an overview.
Complete two dimensional  orders of finite lattice type were studied
by Artin \cite{ArtinOrders},  Reiten and van den Bergh \cite{Reiten}.
The following question is very natural.

\begin{question}\label{Q:latticetype}
Let $(\rA, \idm)$ be a complete Cohen-Macaulay  ring of Krull dimension two,
and $\Lambda$ be an $\rA$-order. When the category of lattices over $\Lambda$ (i.e.
the category of finitely generated left $\Lambda$-modules, which are Cohen-Macaulay
over $\rA$) is of tame representation type?
\end{question}

Some ``trivial''  examples of such orders are given by
skew group  orders $\rA*G$, where $(\rA, \idm)$ is either
a minimally elliptic or a cusp singularity and $G$ a finite group of automorphisms
of $\rA$. However, the general answer on question
\ref{Q:latticetype} is widely unknown.

Moreover, we also do not know the answer about the characterization
of Cohen-Macaulay tame  two-dimensional normal Noetherian local rings, which
are not algebras over a field. It is an interesting problem to introduce
an arithmetic notion of a minimally elliptic singularity and to generalize
methods of the geometric McKay correspondence on the arithmetic case. Existence
of resolutions of surface singularities  in this general situation was proven by Lipman
\cite{Lipman}.

In a forthcoming paper of both authors we are going to show that the
non-isolated surface singularities called \emph{degenerate cusps}
are Cohen-Macaulay tame \cite{BurbanDrozd}.
This includes,  for example,  the
complete intersections $\kk\llbracket x,y,z\rrbracket /(xyz)$ and $\kk\llbracket x,y,z,t\rrbracket /(xy, zt)$
as well as their deformations.   A description of matrix
factorizations corresponding to the Cohen-Macaulay modules
of rank one and two over
the completion of the affine cone of  a nodal cubic curve
$\kk\llbracket x,y,z\rrbracket /(zy^2 + x^3 + x^2 z)$
was obtained by Baciu \cite{Baciu}.

This raises the question about a complete classification
of two-dimension Cohen-Macaulay local rings of tame
Cohen-Macaulay representation type.

\begin{conjecture}\label{C:CMreptype}
Let $(\rA, \idm)$ be a normal surface singularity over a
field $\kk = \rA/\idm$, which is algebraically closed  of characteristic zero.
Then $\rA$ is of tame Cohen-Macaulay representation type if and only
if it has  the form
$\rA = \rB^G$, where $(\rB, \idn)$ is either a simple elliptic  of a cusp singularity
and $G$ is a finite group of automorphisms of $\rB$.
If $\kk = \mathbb{C}$, this  is equivalent to the
condition that the local fundamental group
 $\pi_1(\sX, o)$ is  infinite and solvable, see \cite{Karras}.
\end{conjecture}

Recently, Drozd and Greuel have proven  that a rational normal surface singularity
of tame Cohen-Macaulay representation type is log-canonical \cite{DGrat}.
In a work of
Bondarenko \cite{Bondarenko} it was shown that a hypersurface singularity
$\rA = \kk\llbracket x,y,z\rrbracket /f$ with $f \in (x,y,z)^4$ is always of Cohen-Macaulay wild.  However, the answer on the following question is still unknown.

\begin{conjecture}\label{C:CMtrihot}
Let $(\rA, \idm)$ be a local Cohen-Macaulay ring of Krull dimension two.
Then the
representation type of the category of Cohen-Macaulay $\rA$--modules
 $\CM(\rA)$ is either finite, discrete, tame or wild.
\end{conjecture}

A positive solution of this conjecture was obtained by Drozd and Greuel
in the case of the reduced curve singularities, see \cite{DrozdGreuel}.
Knowing an answer on the following problem   would considerably help
to solve Conjecture \ref{C:CMreptype} and Conjecture \ref{C:CMtrihot}.

\begin{conjecture}\label{C:semicont}
Let $\sX \to \sT$ be a flat family of two-dimensional surface singularities.
Then the set $
\sB =\bigl\{t \in \sT \, \big|  \, \sX_t \,\, \mbox{is of wild Cohen-Macaulay representation
type}\bigr\}
$
is Zariski-closed  is $\sB$.
\end{conjecture}

\noindent
This  conjecture essentially means  that a Cohen-Macaulay tame surface singularity
can not  be locally deformed to a Cohen-Macaulay wild singularity.
Such semi-continuity result is known  in the case of the reduced
 curve singularities, see \cite{Knoerrer1} and \cite{semicont}.
Note that the following proposition is true.

\begin{proposition}
Let $\kk$ be an algebraically closed field of characteristic zero,
$\sX \to \sT$ be a flat family of normal
two-dimensional surface singularities  and  $t_0 \in \sT$ be  a closed point such that
$\sX_{t_0}$ has finite Cohen-Macaulay representation type. Then there exists
an open neighborhood $\sU$ of $t_0$ in $\sB$ such that for all
$t \in \sU$ the singularity  $\sX_t$ is of finite Cohen-Macaulay representation type.
\end{proposition}

\noindent
\emph{Proof}. Indeed, since the singularity $\sX_{t_0}$ has only finitely many
indecomposable Cohen-Macaulay modules, it is a quotient singularity. However, by a result
of Esnault and Viehweg \cite{EsnaultViehweg}
it is known that the quotient surface singularities
deform to quotient surface singularities. This implies the claim.
\qed

\medskip
Another approach to study Cohen-Macaulay modules is provided by the theory
of cluster tilting.

\begin{proposition}[see Theorem 4.1 and Theorem 7.6 in \cite{BIKR}]\label{P:BIKR}
Let $\rA$ be a minimally elliptic curve singularity given by the equality
 $$\rA = T_{p,q}(\lambda)  = \kk\llbracket x,y\rrbracket /(x^p + y^q +  \lambda x^2 y^2),$$
 where
 $\frac{1}{p} + \frac{1}{q} \le \frac{1}{2}$, $\lambda \in \kk^*$
and either both $p$ and $q$ are even, or $p = 3$ and $q$ is even.
Then there exists
a Cohen-Macaulay $\rA$-module $\mM$ such that
\begin{itemize}
\item It is rigid, i.e. $\Ext^1_\rA(\mM, \mM) = 0$ and $\mM \not\cong  0$.
\item If $\mN$ is any indecomposable Cohen-Macaulay $\rA$-module such that
$\Ext^1_\rA(\mN, \mM) = 0$, then $\mN$ is a direct summand of $\mM$.
\end{itemize}
Such an object $\mM$ is called \emph{cluster-tilting}.
Moreover, if $\Gamma = \End_{\underline{\CM}(\rA)}(\mM)$ is the corresponding
\emph{cluster-tilted algebra} then the functor
$$
\Hom_{\underline{\CM}(\rA)}(\mM, \,-\,) : 
\underline{\CM}(\rA)/\tau(\mT) \lar 
\mod-\Gamma 
$$
is an equivalence of categories, where $\underline{\CM}(\rA)$ is the stable category of
Cohen-Macaulay modules and $\tau$ is the Auslander-Reiten translation in
$\underline{\CM}(\rA)$.
\end{proposition}

Moreover, the cluster-tilted algebras $\Gamma$ arising
in Proposition \ref{P:BIKR} are  finite-dimensional and  symmetric. From works of
Erdmann, Bia\l{}kowski-Skowro\'nski and Holm it follows they
are tame, see
\cite{Skowronski}  for an overview  about  tame  self-injective algebras.

\section{Stable category of Cohen-Macaulay modules and computations with \tt{Singular}}

\medskip
Let $(\rA,\idm)$ be a Gorenstein local ring. Since by the definition $\mK_\rA \cong  \rA$,
the exact category $\CM(\rA)$ is Frobenius.

\begin{definition}\label{D:stabcat}
The stable category of Cohen-Macaulay modules
$\underline{\CM}(\rA) = \CM(\rA)/\langle \rA\rangle$ is defined as follows
\begin{itemize}
\item $\Ob\bigl(\underline{\CM}(\rA)\bigr) = \Ob\bigl(\CM(\rA)\bigr)$
\item $\underline{\Hom}_\rA(\mM, \mN) =
\Hom_{\underline{\CM}(\rA)}(\mM, \mN) := \Hom_\rA(\mM, \mN)/\mathfrak{P}(\mM, \mN)$, where
$\mathfrak{P}(\mM, \mN)$ is the submodule of $\Hom_\rA(\mM, \mN)$
generated by those morphisms which factors
through a free $\rA$-module.
\end{itemize}
\end{definition}

The following theorem can be considered as the \emph{raison d'\^etre} for  the  study
of Cohen-Macaulay modules

\begin{theorem}[Buchweitz]\label{T:stabcat}
Let $(\rA, \idm)$ be a complete Gorenstein local ring, then the following
holds,
\begin{itemize}
\item
The  functor  $\Omega = \syz: \underline{\CM}(\rA) \lar \underline{\CM}(\rA)$
is an auto-equivalence of $\underline{\CM}(\rA)$.
\item The category $\underline{\CM}(\rA)$ has a structure of a triangulated category,
where the shift functor is $T:= \Omega^{-1}$.
\item There is an equivalence of triangulated categories
$$
\underline{\CM}(\rA) \to D^b(\rA-\mod)/\Perf(\rA)
$$
induced by the  fully faithful functor $\CM(\rA) \to D^b(\rA-\mod)$
\item
In particular, the  canonical functor $\CM(\rA) \to \underline{\CM}(\rA)$ maps exact sequences
into exact triangles.
\end{itemize}
\end{theorem}

\noindent
This theorem was proven for the first time in \cite{Buchweitz}. Recently, it was
rediscovered in \cite{Orlov1}.

\medskip
\begin{theorem}[Eisenbud \cite{Eisenbud}]\label{T:matrfact}
Let $(\rS,\idn)$ be a regular local ring, $f \in \idn$ and $\rA = \rS/f$. Then any Cohen-Macaulay module
$\mM$ without free direct summands has a 2-periodic minimal free resolution:
$$
 \dots \stackrel{\mu}\lar \mF \stackrel{\nu}\lar \mF \stackrel{\mu}\lar \mF \lar M \lar 0,
$$
where $\mF \cong \rA^n$ is a free $\rA$-module.
In these  terms we can write $\mM = \mM(\mu, \nu)$, where $\mu, \nu \in \Mat_n(\idm)$.
This implies that $\Omega^2 \cong Id$, hence $T^2 \cong Id_{\underline{\CM}(\rA)}$.
Moreover, we have
an equivalence of triangulated categories:
$$
\underline{\CM}(\rA) = \Hot_2(\rA),
$$
where $\Hot_2(\rA)$ is the homotopy category of (unbounded) minimal 2-periodic projective
complexes.
\end{theorem}

\medskip
\begin{proposition}\label{C:Homfin}
Let $(\rA, \idm)$ be a Cohen-Macaulay ring and $\mM$, $\mN$ be Cohen-Macaulay modules.
If either $\mM$ or $\mN$ is locally free on the punctured spectrum, then
the $\rA$-module $\underline\Hom_{\rA}(\mM, \mN)$ has finite length.
\end{proposition}

\noindent
\emph{Proof}. Take any short exact sequence 
 $0 \to \mN' \to \rA^n \to \mN \to 0$ and observe that 
 if $\mN$ is locally free on the punctured spectrum then
 $\mN'$ is locally free on the punctured spectrum, too. 
 From the exact sequence
 $$
 \Hom_\rA(\mM, \rA^n) \to \Hom_\rA(\mM, \mN) \to 
 \Ext^1_\rA(\mM, \mN')
 $$
we obtain an embedding $\underline{\Hom}_\rA(\mM, \mN)
\hookrightarrow \Ext^1_\rA(\mM, \mN')$. If either
$\mM$ or $\mN$ is locally free on the punctured spectrum, then 
$\Ext^1_\rA(\mM, \mN')$ has finite length, hence the claim. 
\qed

\medskip
\noindent
It turns out that the converse statement is also true.

\begin{theorem}\label{T:isolHomfin}
Let $(\rA, \idm)$ be a Cohen-Macaulay local ring. Then it is an isolated
singularity if and only for all Cohen-Macaulay modules $\mM$ and $\mN$ the
module
$\Ext^1_\rA(\mM, \mN)$ has finite length.
\end{theorem}

\noindent
For a proof of this theorem, we refer to \cite[Lemma 3.3]{Yoshino}.

\medskip

\begin{theorem}[Auslander]\label{T:SerreDual}
Let $(\rA, \idm)$ be an isolated  Gorenstein singularity 
of Krull dimension $d$ and  $\kk = \rA/\idm$. Then for any pair of
Cohen-Macaulay  modules $\mM$ and $\mN$ we have a bifunctorial isomorphism of $\rA$-modules
$$
\Ext^1_\rA\bigl(\mM, (\syz^d\Tr(\mN))^*\bigr)  \cong 
\DD\bigl(\underline\Hom_\rA(\mM, \mN)\bigr),
$$
where $\mX^* = \Hom_\rA(\mX, \rA)$ and $\DD = 
\Hom_\rA\bigl(\,-\,,\mE(\kk)\bigr)$ is  the Matlis duality functor.
Moreover, for any Cohen-Macaulay module $\mN$ there is a functorial isomorphism
$
\syz^d\bigl(\Tr(\mM)\bigr)^* \cong \syz^{2-d}(\mM): \underline{\CM}(\rA) \to \underline{\CM}(\rA).
$
\end{theorem}

\noindent
For a proof of this result, see \cite[Proposition 8.8 in  Chapter 1
and Proposition 1.3 in  Chapter 3]{PhilNotes}. 

\medskip
\noindent
If $(\rA, \idm)$ is a Gorenstein $\kk$-algebra, then this theorem can be restated as follows.

\begin{corollary}\label{C:SerreDual}
Let $(\rA, \idm)$ be a Gorenstein $\kk$-algebra ($\kk = \rA/\idm$)
of Krull dimension $d$, which is an isolated singularity.
 Then  $\SS = \syz^{1-d}$ is the Serre functor in the stable category of Cohen-Macaulay modules
$\underline{\CM}(\rA)$. This means that we for any two Cohen-Macaulay modules $\mM$ and $\mN$ we
have a bifunctorial isomorphism
\begin{equation}\label{E:SerreDual}
\underline{\Hom}_\rA(\mM, \mN) \cong \underline{\Hom}_\rA\bigl(\mN, \SS(\mM)\bigr)^*.
\end{equation}
\end{corollary}

\medskip

\begin{remark}\label{R:Buchweitz}
By  \cite[Proposition 10.1.5]{Buchweitz}  the isomorphism (\ref{E:SerreDual}) holds
for an arbitrary Gorenstein ring $d$ and a pair of Cohen-Macaulay modules $\mM$ and $\mN$ such that
$\mM$ is locally free on the punctured spectrum. This means, that the stable category of Cohen-Macaulay
modules over a Gorenstein $\kk$-algebra $\rA$ of Krull dimension $d$, which are locally free on the
punctured spectrum, is  a triangulated $(d-1)$-Calabi-Yau  category.
\end{remark}

\medskip

The following lemma shows that the stable categories of Cohen-Macaulay modules
over an isolated Cohen-Macaulay singularity and its completion are closely related.

\begin{proposition}\label{P:complet}
Let $(\rA,\idm)$ be an isolated Cohen-Macaulay singularity and $\widehat{\rA}$ its completion.
Then the canonical functor
$\underline{\CM}(\rA) \to \underline{\CM}(\widehat{\rA})$ is fully faithful.
\end{proposition}

\noindent
\emph{Proof}. Let $\mM$ and $\mN$ be two Cohen-Macaulay $\rA$-modules. Since they are automatically
free on the punctured spectrum, the $\rA$-module $\underline{\Hom}_\rA(\mM, \mN)$ is annihilated by some
power of the maximal ideal: $\idm^t \cdot  \underline{\Hom}_\rA(\mM, \mN) = 0$ for $t \gg 0$.
Hence it is isomorphic to its completion
$\underline{\Hom}_{\widehat{\rA}}(\widehat{\mM}, \widehat{\mN})$.
\qed

\medskip

In order to compute the dimensions of  $\Hom$ and $\Ext$--spaces in the stable category of
maximal Cohen-Macaulay modules,  one can use
the computer algebra system {\tt Singular}, see
\cite{GP}.
Let $$\rA = \kk[x_1,x_2,\dots,x_n]_{\langle x_1, x_2,\dots, x_n \rangle}/I$$
 be a Cohen-Macaulay local ring, $\mM$ and $\mN$ be a pair of  maximal Cohen-Macaulay modules.

Assume  the vector space  $\Ext^i_\rA(\mM,\mN)$ ($ i \ge 1$) is finite-dimensional over $\kk$. Since
the functor
$\rA-\mod \to \widehat{\rA}-\mod$ is exact, maps the maximal Cohen-Macaulay modules to  maximal Cohen-Macaulay modules and the finite length modules to finite length modules, we can conclude that
$$
\dim_\kk\bigl(\Ext^i_\rA(\mM,\mN)\bigr) = \dim_\kk\bigl(\Ext^i_{\widehat{\rA}}(\widehat{\mM},
\widehat{\mN})\bigr).
$$
Moreover, if $\rA$ is a hypersurface singularity then
by Theorem \ref{T:matrfact} we have:
$$\underline\Hom_\rA(\mM, \mN) \cong \Ext^2_\rA(\mM, \mN).$$

\begin{example}\label{E:Singular}
Let $\rA = \kk\llbracket x,y,z\rrbracket /xyz$.  Then the following modules
$$
\rA^3 \xrightarrow{
\left(
\begin{array}{ccc}
x & z^3 & 0 \\
0 & y & x  \\
y^2 & 0 & z
\end{array}
\right)
}
\rA^3 \lar \mM \to 0,
$$
$$
\rA^3 \xrightarrow{
\left(
\begin{array}{ccc}
x & 0 & 0 \\
z^3 & y & 0  \\
y^2 & x & z
\end{array}
\right)}
\rA^3 \lar \mN \to 0,
$$
and
$$
\rA^2 \xrightarrow{
\left(
\begin{array}{cc}
xy  & -x^2 + y^3 \\
0 & z  \\
\end{array}
\right)} \rA^2 \lar \mK \to 0
$$
are Cohen-Macaulay and locally free on the
punctured spectrum of $\rA$. Let us compute certain $\Hom$ and $\Ext$ spaces between $\mM$, $\mN$ and
$\mK$.

\medskip
\noindent
{\tt > Singular}  (call  the program {\tt ``Singular''})\\
{\tt > LIB ``homolog.lib'';} (call the library of homological algebra)\\
{\tt > ring S = 0,(x,y,z),ds;}   (defines the ring $\rS = \mathbb{Q}[x,y,z]_{\langle x, y, z\rangle}$)\\
{\tt > ideal I = xyz;}     (defines the ideal $xyz$ in $S$) \\
{\tt > qring A  = std(I);} (defines the ring $\mathbb{Q}[x,y,z]_{\langle x, y, z\rangle}/I)$\\
{\tt > module k = [x], [y], [z];} (defines the residue field $\kk$ of $\rA$ as an $\rA$-module)\\
{\tt > module M = [x, 0, y2], [z3, y, 0], [0,x,z];} \\
{\tt > module N = [x,z3,y2], [0,y,x], [0,0,z];} \\
{\tt > module K = [xy,0], [-x2 + y3, z];} \\
{\tt > isCM(M);} \hspace{2cm}(checks, whether $\mM$ is Cohen-Macaulay) \\
{\tt > 1} \hspace{4cm}(yes, $\mM$ is Cohen-Macaulay) \\
{\tt > isCM(N);} \hspace{2cm}(checks, whether $\mN$ is Cohen-Macaulay) \\
{\tt > 1} \hspace{4cm}(yes, $\mN$ is Cohen-Macaulay) \\
{\tt > list l = Ext(1,k,N,1);} \\
{\tt // ** redefining l **} \\
{\tt // dimension of $\Ext^1$:  -1} \hspace{1cm}  $\Ext^1_\rA(\kk, \mN) = 0$ \\
{\tt > isCM(K);} \hspace{2cm} (checks, whether $\mK$ is Cohen-Macaulay) \\
{\tt > 1} \hspace{4cm} (yes, $\mK$ is Cohen-Macaulay) \\
{\tt > list l = Ext(2,M,M,1);} \\
{\tt // ** redefining l **} \\
{\tt // dimension of $\Ext^2$:  0} \hspace{2cm} ($\krdim(\Ext^2_\rA(\mM, \mM)) = 0$) \\
{\tt // vdim of $\Ext^2$:       7} \hspace{1cm}
($\dim_\kk\bigl(\Ext^2_\rA(\mM, \mM)\bigr) =
\dim_\kk\bigl(\underline\Hom_\rA(\mM, \mM)\bigr) = 7$) \\
{\tt > list l = Ext(2, N, K,1);} \\
{\tt // ** redefining l **} \\
{\tt // dimension of $\Ext^2$:  0} \\
{\tt // vdim of $\Ext^2$:       7} \hspace{2cm}
$\dim_\kk \Ext^2_\rA(\mM, \mM) \cong \underline\Hom_\rA(\mM, \mM) = 7$ \\
\end{example}

\end{document}